\documentclass  {thesis}
\usepackage{alltt}  %
\usepackage{amsfonts}

\begin{document}

\newcommand{\eps}{\varepsilon}
\newcommand{\Ztwo}{\mbox{${\bf Z}^2$}}              
\newcommand{\Nzero}{\mbox{${\bf N}_0$}} 
\newcommand{\la}{\mbox{$\longleftrightarrow$}}
\newcommand{\pd}{\mbox{${\cal D}$}}
\newcommand{\f}{\varphi}
\newcommand{\gtwo}{\mbox{$G^{(2)}$}}

\def\qed{\vbox{\hrule\hbox{\vrule\kern3.5pt
  \vbox{\kern7.5pt}\kern4pt\vrule}\hrule}}
\def\fine{\hfill\qed}

\prelimpages


\Title{Discrete Growth Models}
\Author{Dorothea M. Eberz-Wagner}
\Year{1999}
{\Degreetext{A dissertation submitted in partial fulfillment of\\
the requirements for the degree of \\
\\
Doctor of Philosophy (University of Washington).\\
}
%
\setcounter{page}{-1}

\centerline{\bf Discrete Growth Models}
\bigskip
\centerline{\bf Dorothea M. Eberz-Wagner}
\bigskip
\centerline{1999}
\bigskip
A dissertation submitted in partial fulfillment of
the requirements for the degree of 
Doctor of Philosophy (University of Washington).

\tableofcontents

\acknowledgements{I want to express my sincere gratitude to my advisor, Krzysztof Burdzy, for his encouragement and patience throughout the years, and for his invaluable advice. I am also grateful to the members of the reading committee for their work, and to the mathematics department, in particular to the graduate advisors who provided me with various kinds of support. I want to thank the staff of the department and the helpers at Computing and Communications for their time and their willingness to help.

Special thanks go to all of my friends who freely advised me on TeX, to my friend Susan Chin for her tremendous encouragement, to my parents for their
outstanding support, to my husband Mark for his love and tolerance, and to my children for being the way they are. }


\textpages

\newtheorem{lemma}{Lemma}[section]
\newtheorem{proposition}{Proposition}[section]
\newtheorem{corollary}{Corollary}[section]
\newtheorem{theorem}{Theorem}[section]
\newtheorem{notation}{Notation}[section]
\newtheorem{remark}{Remark}[section]
\newtheorem{observation}{Observation}[section]

\chapter{Introduction}

Our thesis is motivated by the celebrated growth model on the hypercubic lattice ${\bf Z}^d$ known as Diffusion Limited Aggregation (DLA). This model was introduced in 1981 by Witten and Sander \cite{Witten&Sander:DLA}
to model aggregates of a condensing metal vapour and since then has attracted much interest both among physicists and mathematicians. The process starts with only the origin ${\bf 0}$ being occupied. Then, one by one, particles diffuse in ``from infinity" according to a random walk conditioned to hit the boundary of the current cluster and attach themselves to the first boundary site visited. (The boundary of a cluster consists of all points of ${\bf Z}^d$ adjacent to the cluster but not belonging to it). Let $A_n$ denote the cluster existing after $n$ steps, $A_0=\{0\}$. Very little is known about the shape of a typical cluster for large values of $n$. Computer simulations suggest that for large $n$ these clusters are highly ramified and of low density and their radius $r(A_n)$ has greater order of magnitude than the $d$-th root of the number of vertices. 

One hard open problem is to prove the existence of $\beta_d= \lim_{n \rightarrow \infty} \log r(A_n)/ \log n$, to show that this random variable is almost surely a constant and to determine its value. So far no lower bound could be rigorously established aside from the trivial, $\liminf_{n \rightarrow \infty} \log r(A_n)/ \log n \geq 1/d$.
As for an upper bound,
Kesten \cite{Kesten:how long,Kesten:growth rate} 
showed that with probability 1, $\limsup_{n \rightarrow \infty} \log r(A_n)/ \log n  \leq 2/{(d+1)}$. A conjectured value of the limit which agrees fairly well with numerical simulations is ${(d+1)}/{(d^2+1)}$ \cite[ Sec. 2.6]{Lawler:Intersections}.

The strong interest in DLA gave rise to a number of variants. Bramson, Griffeath and Lawler \cite{Bramson: Internal DLA} studied a model called Internal Diffusion Limited Aggregation. The growth rule for this model can be described as above, the only difference being that new particles get released at the origin and then diffuse through the interior of the current cluster until they attach themselves at the first boundary site visited. Thus, roughly speaking, these particles tend to stop at sites which are closest to the origin, whereas in standard DLA new particles tend to favor neighbors of ``extreme sites". Due to this difference no fingering of the growing clusters can be expected. Bramson et al. proved in \cite{Bramson: Internal DLA}
that the clusters formed by internal DLA have the limiting shape of a sphere. Kesten \cite{Kesten:caricatures} considered multiple contact DLA and the $\eta$-model, both being variants of DLA in which the distribution of the place where a new particle is added is shifted in a strong way to the points of maximal harmonic measure. He showed that these clusters can grow like ``generalized plus signs", i.e., with positive probability the growing cluster contains only points on the coordinate axes at all times. Barlow, Pemantle and Perkins \cite{Barlow:DLa on a tree}  studied DLA on a regular d-ary tree for tree-valued random walks with transition probabilities determined by assigning conductance $\alpha^{-1}$ to edges from generation $n-1$ to generation $n$, $\alpha<1.$ Thus the cluster starts with only the origin being occupied and new points are added according to harmonic measure on the boundary of the current cluster for the random walk under consideration. This model is mathematically more tractable than DLA on ${\bf Z}^d$  since there is a simple formula for harmonic measure on the boundary of a cluster and the absence of loops implies that disjoint parts of the cluster evolve nearly independently. In \cite{Barlow:DLa on a tree} Barlow et al. rigorously prove that these processes exhibit fingering. The height of the growing clusters increases linearly with their total size. 

In the first part of this thesis we study a one-dimensional growth model on ${\bf Z},$ which starts with only the origin being occupied, and in which growth is obtained by first adding $K$ particles and then deleting $K-1$ particles in each step, $K \geq 2$ being a fixed positive integer. Thus the size of the set of occupied sites increases by one in each step. When adding a new particle, its site is chosen uniformly from the boundary of the present set and when deleting a particle, the site to be removed is chosen uniformly from the set of all occupied sites. Let $A_n$ denote the set of occupied sites after $n$ steps are completed. The main interest focuses on the evolution of the gaps of $A_n$ for large $n$. (A gap of $A_n$ is a nonempty subset of ${\bf Z} \setminus A_n$ of the form $\{z \in {\bf Z}  : z_1 < z < z_2\},$ where $z_1,z_2 \in A_n$). We will prove that with probability one eventually all gaps, with the possible exception of one, have size one, this one exceptional gap (if present) having size two. Infinitely often such an exceptional gap is present and infinitely often there are only gaps of size one. Moreover, letting $L_n$ denote the number of unoccupied sites between the endpoints of $A_n$ and defining $h(n)=\log n/{\log \log n},$ we will show that with probability one, $\liminf_{n \rightarrow \infty} L_n =K-2$ and $\limsup_{n \rightarrow \infty} L_n/h(n) ={(K-1)}/K.$ 

We will start out with straightforward estimates and a supermartingale argument to obtain the first results concerning the evolution of the gaps of $A_n.$ In particular we will show that with probability one infinitely often the set $A_n$ has at most $10K$ gaps. To obtain further results, we will then introduce the distinction between ``new" and ``old" gaps, i.e., we will define a Markov process $(A,I),$ where $I_n$ denotes the index set of new gaps of $A_n.$ 

To illustrate how we will work with this process, let us investigate a little bit closer the evolution of the process $L^I$, where $L^I_n$ denotes the number of gaps sites belonging to new gaps of $A_n$. 
First note that when we start with a sufficiently big set $a$, which has just a few gaps, (all of which are considered as being old), then the next deletions will most often create new gaps of size one. Additions affect both new and old gaps as well as exterior neighbors of the endpoints of $A_n,$  and so we will most likely see that  $L^I_n$ reaches higher and higher values. Yet, once  the number of new gaps is sufficiently big, we observe that an overwhelming majority of additions will (most likely) fill into new gaps and only a small portion will affect old gaps or exterior neighbors of the endpoints. Since we add $K$ particles and delete $K-1$ particles, the number of gaps sites belonging to new gaps will tend to decrease. In other words, if the number of old gaps is bounded above by a fixed constant, then $L^I$ tends to decrease whenever it is sufficiently big, and this holds as long as the number of new gaps is close to the total number of gap sites belonging to new gaps. 

To make use of this fact we will define a stopping time $\sigma$ which has the property that up to this time at most two gap sites belonging to new gaps are adjacent to each other. Moreover, we will define so-called ``excursions" of $L^I.$ Let us think of the path $L^I(\omega)$ as a ``staircase". We now cut off that part of the staircase which is below a certain given constant $C$. What remains are ``separate" pieces each of these being a ``staircase" all the values of which are greater or equal to $C$. These pieces will be called ``excursions" of $L^I(\omega).$ (This description, besides being very informal, does not contain certain technical details, but should convey the main idea. When working out the details we will define the end of an excursion as the first time - after the start of the excursion - when the value of $L^I(\omega)$ is strictly less than $C$). 

We will next use a special construction of the process $(A,I)$. We will construct this process from a sequence $U^{(1)}$ of independent random variables which are uniformly distributed on $(0,1)$. Using the same sequence $U^{(1)}$ we will then construct a random walk ${\hat L}^I$ which has the property, that up to time $\sigma$ the increments of $L^I$ on excursion intervals will be dominated by the corresponding increments of ${\hat L}^I.$ Yet, since the increments of ${\hat L}^I$ may at some places be strictly bigger than the corresponding increments of $L^I$, the part of ${\hat L}^I$ which corresponds to an excursion of $L^I$, may not ``look like an excursion", in fact, the value at its end may be even bigger than the value we start with. We will thus use an additional independent process $U^{(2)}$ to continue those parts until for the first time the value of ${\hat L}^I$ is by an amount $K$ smaller than the value at the beginning of this part.  

Similar algorithmic constructions will be used at several occasions, all of these constructions fitting into one general scheme. To understand this scheme think of a gambler in a casino who may skip games or switch tables according to the outcome of previous games. These constructions will enable us to obtain set inclusions on an enriched probability space
which will then allow us to obtain suitable probability estimates for the sets of interest.

Using these methods we will obtain our two main probability estimates. The first estimate will be valid when we start out with a sufficiently large subset of $\bf Z$ which has just a few gaps (all of these being considered as being old). Using that with probability one infinitely often the sets $A_n$ have at most $10K$ gaps,
this estimate together with the Markov property of $A$ will allow us to show that a certain family ${\cal S}$ of subsets of $\bf Z$ is recurrent for the process $A.$  The second estimate will be valid when we start out with a subset $a$ of $\bf Z$ which already belongs to $\cal S.$ (We now do not distinguish between new and old gaps, or, equivalently, we consider all gaps of $a$ as being new).
Using this second estimate plus a Borel-Cantelli argument we will then be able to prove all, but one, of the remaining claims. This last piece will essentially follow from the strong law of large numbers, since we can use our previous results, and, once again, our special construction of the process $A$.

In the case $K=2$ we will also study the age $O_n$ of the oldest particle of the set $A_n.$ We will show that $(n-O_n)/{\sqrt n}$ converges in distribution to $F$, where $F(x)=(1-\exp(-x^2)) {\bf 1}_{(0, \infty)}.$

The second part of the thesis deals with diffusion limited aggregation on ${\bf Z}^2.$ We attack the question whether the growing clusters will have holes. (A hole of the cluster $A_n$ is a finite maximal connected subset of ${\bf Z}^2 \setminus A_n$). Clearly, with positive probability, a hole may be formed within the first few steps. Yet, is it true that almost surely the growing clusters will eventually have at least one hole? The answer to this question is ``yes". In fact we will show that with probability one the number of holes of the growing cluster tends to infinity. 

The main idea of the proof can roughly be described as follows. Given a cluster $a,$ let us change just a few paths of the next particles by inserting an extra loop such that the cluster obtained with these changed paths will have at least one more hole than the given cluster. More precisely, let us choose a constant $L$ (depending on $a$) and let us start the next $L$ random walks on the boundary of a ``ball" (which surrounds the cluster $a$ and which is sufficiently far away) according to harmonic measure on the boundary of this ``ball". Let $\omega=(\omega_1, \ldots, \omega_L)$ describe these random walks. We will now define changed paths $\varphi(\omega)$ by inserting extra loops into some of the $L$ given paths such that the cluster obtained from $a$ with these changed paths will have at least one more hole than the given cluster $a.$ Thus for any $\omega=(\omega_1, \ldots, \omega_L),$ we will decide which paths will get changed, and if $\omega_i$ is such a path, we will specify at which time an extra loop will get inserted and what this extra loop will look like. In order to make these decisions it will not be necessary to know the full paths, but it will suffice to know these paths up to the hitting time of the set $a.$ If an extra loop will be inserted into path $\omega_i,$ then the insertion will be done before $\omega_i$ hits $a.$ Moreover, the number of inserted loops as well as the length of each loop will be bounded above by a fixed constant and the map $\varphi$ will be one-to-one. We can then conclude that there exists a constant $c>0$ such that for any finite cluster $a$ there exists an integer $L$ with $P(A_L \mbox{ has at least one more hole than } A_0 \mid A_0=a) \geq c.$ The above result will then follow by the Markov property of the process $A.$ 

Yet, the construction of $\varphi$ is not an easy task. Just suppose that we insert a loop into the path of the first particle to make it attach at a designated point $x_1$ of the boundary of $a$ which is different from the previous point of attaching $y_1.$ If the second particle visits a neighbor $x_2$ of $x_1$ strictly before hitting the boundary of $a \cup \{y_1\}$, then $x_2$ will be the new point of attaching for the second particle. Moreover, the point of attaching of the second particle will also be changed, if this particle visits a neighbor of $y_1$ strictly before visiting the boundary of $a \cup \{x_1\}.$ Similar changes may occur for the point of attaching of the third particle and so on. So even a slight change of a single path may cause a chain reaction, possibly changing the positions of all particles attached thereafter. We will thus have to figure out how to deal with this difficulty in order to solve the given task.

\chapter{Gaps in a One-Dimensional Annihilation-Creation Model}

\section {Description of the Model and Statement of Results}

Consider a stochastic process $\{A_n\}_n$ whose state at time $n,$ $A_n
\subset {\bf Z},$ is the set of sites occupied  by a finite set of
particles. Let $A_0=\{0\}$ and fix some integer $K \geq 2.$ Given $A_n,$ the set $A_{n+1}$ is constructed by
first adding $K$ sites one by one and then deleting $K-1$ sites one by
one. To add a particle, its new site is chosen uniformly from the boundary
of the given set of occupied sites, (a site $y$ belongs to the boundary $\partial a$ of a set $a \subset {\bf Z}$ if 
$y \in {\bf Z} \setminus a$ and if there exists $x \in a$ with  $|x-y|=1$). To delete a particle, a site is chosen
uniformly from the set of all occupied sites. Hence the number of elements
of $A_n$ increases by~1 in each period. We are interested in the evolution
of the gaps of~$A_n$, (a gap  of a set $a \subset {\bf Z}$ being defined as a nonempty subset $g$ of ${\bf Z} \setminus a$ with the property that there exist $x_1, x_2
\in a$ with $g=\{y \in {\bf Z}: x_1 < y < x_2 \}$).
 
\vspace{7 mm}
For $a \subset {\bf Z}$ let
$G(a)$ denote the
number of gaps of $a,$ let $G^{(2)}(a)$ denote the number of gaps with
$|g| \geq 2$ and let $L(a)$ denote the number of unoccupied sites between the
endpoints of $a,$ i.e., $L(a)=|g_1|+ \ldots +|g_{G(a)}|,$ where $g_1,
\ldots, g_{G(a)}$ are the different gaps of $a.$
Letting $h(n)=\log n/\log\log n,$ we will prove the following:

\begin{theorem} \label{first theorem}
With probability 1, \newline
1) $\limsup_{n \rightarrow \infty} L(A_n)/h(n) = (K-1)/K,$ \newline 
\nopagebreak
2) $\liminf_{n \rightarrow \infty} L(A_n) = K-2.$ 
\end{theorem}

\begin{theorem} \label{second theorem}
With probability 1, \newline
1) $G^{(2)}(A_n) = 1$ i.o., \newline
2) $G^{(2)}(A_n) = 0$ i.o., \newline
3) $L(A_n) \leq G(A_n)+1$ eventually.
\end{theorem}

Hence eventually all gaps, with the possible exception of one, have size one, this one exceptional gap having size 2. Infinitely often such an exceptional gap is present and infinitely often there are only gaps of size one.

\begin{remark} The above results are valid for any finite initial set $A_0.$
\end{remark}

We will also discuss the age of the oldest particle of the set $A_n.$
For $x \in A_n$ let $Y_n(x)$ be the number of deletions
which happened since the present particle at $x$ was added and let $O_n=\max
\{Y_n(x): x \in A_n \}$ denote the age of the oldest particle of
the set $A_n.$ 
Define a distribution function $F$ by $F(x)=(1-\exp(-x^2)) {\bf 1}_{(0,
\infty)}.$
In the case $K=2$ we have the following result
\begin{theorem} \label{Theorem: age of the oldest particle}
$(n-O_n)/\sqrt n \Longrightarrow F.$
\end{theorem}

\section {Overview}

Using straightforward estimates and a supermartingale argument we will derive the first few results concerning the behavior of $L,$ $G$ and $G^{(2)}.$ We will show that almost surely $\limsup_{n \rightarrow \infty} L(A_n)/h(n) \geq (K-1)/K,$ $L(A_n) \geq K-2$ eventually, $G^{(2)}(A_n) \geq 1$ i.o. and $G(A_n) \leq 10 K$ i.o. When trying to improve upon these results we will have to deal with the fact that the `transition probabilities' for the process $L$ depend on the values of $G$ and $G^{(2)}$ and similar dependencies ooccur when we study the other two processes. The key idea to overcome this difficulty will be to distinguish between `new' and `old' gaps of the set $A_n$ and to observe the evolution of `old' gaps separately from the evolution of `new' gaps. We will thus define a Markov process $(A_n,I_n),$ where $I_n \subset {\bf N}$ will allow to identify gaps as being `new' or `old'. Moreover, stopping times will enable us to study the evolution of one process while controlling other processes within certain bounds. For instance, if, say, at the end of period $n$ the process $G$ is bounded below by 1000, then the probability that $L$ decreases when the first particle gets added in period $n+1$ is at least $1000/1002.$ Using inequalities of this kind, we will compare (on an enriched probability space) the evolution of the process under consideration to the evolution of some other `nicer' process (for instance a random walk on ${\bf Z}$) in a path-by-path way. The proof of our two main lemmas will use this approach and the reader will need some patience to stick with us when going through this algorithmic construction. These two lemmas will allow us to complete the proofs of part (1) of Theorem \ref{first theorem} and of Theorem \ref{second theorem}. Using these results a similar construction will be used to complete the proof of part (2) of Theorem \ref{first theorem}. At the end of this chapter we will investigate the age of the oldest particle for the case $K=2.$ Independently of our previous results we will determine its limiting distribution. 

\vspace{7 mm}

\noindent {\bf Notation.} We will write $P^a$ when we consider the Markov process $A$ with initial set $A_0=a$ and we will write $P$ for $P^{\{0\}}.$ We will write $L_n$ for $L(A_n),$ $G_n$ for $G(A_n)$ etc. Let $\kappa=K+(K-1)$ denote the number of steps in each period and let $A_{n,j},$ $n \in {\bf N},$ $j=0, \ldots, \kappa$ denote the set after $j$ steps in period $n$ are completed. Thus $A_0=A_{1,0}$ and $A_n=A_{n,\kappa}=A_{n+1,0}$ for all $n \in {\bf N}.$ For $a \subset {\bf Z},$ $|a| < \infty,$ any point $x$ with $\min a < x < \max a$ will be called an interior site. Thus an interior site $x$ is occupied if $x \in a,$ otherwise it is unoccupied. 
Let us define the outer boundary of $a$ as ${\bar \partial}a=\{\min a -1, \max a +1\}.$

\section
 {First Results for the Processes $L,G$ and $G^{(2)}$ }

Our first proposition will show that with probability one infinitely often the sets $A_n$ have at least one gap of size bigger than $1.$

\begin{proposition} \label{prop: at least one gap bigger than 1 i.o.}
\qquad \qquad  $P(\gtwo(A_n) \geq 1 \mbox{ i.o.}) =1. $
\end{proposition}
{\bf Proof.} We will show that there exists an index $n_0$ with
\begin{equation} \label{ineq:gaps bigger than 1}
P(\gtwo(A_{n+2}) \geq 1 \mid A_n) \geq  (2/3)^{K+2} \; n^{-1}
\end{equation}
for all $n \geq n_0.$ This will then imply that for all sufficiently large
$n,$ $P(\gtwo(A_l)=0 \mbox{ for all $l \geq n$})$ $\leq$
$\prod_{i=0}^{\infty} (1- (2/3)^{K+2} \; (n+2i)^{-1})=0$ and hence 
$P(\gtwo(A_n)=0 \mbox{ eventually})=0.$
 
So let us prove (\ref{ineq:gaps bigger than 1}).  Note that the set
$A_{n+1}$ has at least one gap if the last deletion in period $n+1$ removes
a particle at an interior site of the present set. Now $A_n$ has size $n+1$
and $K$ sites are added and $K-2$ sites are deleted in period $n+1$ before 
this last deletion occurs. Thus the particle to be removed is chosen from
a set of size $n+3.$ Hence
\[
P(L(A_{n+1}) \geq 1 \mid A_n) \geq (n+1)/(n+3) \geq 2/3
\]
for all sufficiently large $n.$ Next note that if $A_{n+1}$ has at least one gap and if
in period $n+2$ all $K$ particles are added such that at least one interior
site remains unoccupied, then $L(A_{n+2,K}) \geq 1.$ Thus for any 
set $a$ with $|a|=n+2$ and $L(a) \geq 1$ we have 
\begin{equation} \label{(2.2)}
P(L(A_{n+2,K}) \geq 1 \mid A_{n+1}=a) \geq (2/3)^K.
\end{equation}
Moreover, if the next deletion removes a particle which is
adjacent to a gap of $A_{n+1}$ - but not an endpoint - and if the remaining
$K-2$ deletions do not remove any particle which is adjacent to this
particular gap (which is now of size at least 2), then $\gtwo(A_{n+2}) \geq 1.$ Finally
note that if $b$ is a finite set with $L(b) \geq 1$ (and $|b| \geq 3$) then
$b$ has at least one interior site which is adjacent to a gap of $b.$ Thus
for any  set $b$ with $L(b) \geq 1$ and $|b|=n+2+K$ we have 
\begin{eqnarray}
P(\gtwo(A_{n+2}) \geq 1 \mid A_{n+2,K}=b) &\geq& |b|^{-1}
\prod_{i=1}^{K-2} \frac{|b|-i-2}{|b|-i} \nonumber \\
 &\geq& |b|^{-1}  \left( \frac{|b|-K}{|b|-K+2} \right)^{K-2} \nonumber \\
&=& (n+2+K)^{-1} \left( \frac{n+2}{n+4} \right)^{K-2},  \nonumber
\end{eqnarray} 
and thus
\begin{equation}  \label{(2.3)}
P(\gtwo(A_{n+2}) \geq 1 \mid A_{n+2,K}=b)
 \geq  2/3 \; n^{-1}, 
\end{equation} 
if $n$ is sufficiently large. 
Combining (\ref{(2.2)}) and (\ref{(2.3)}) we get for all sufficiently large $n$ and $a \subset {\bf Z}$ with $|a|=n+1,$
$P(\gtwo(A_{n+2}) \geq 1 \mid A_n=a) \geq (2/3)^{K+2} n^{-1} $
and (\ref{ineq:gaps bigger than 1}) thus gets established. 
\newline
\vspace{0 mm}
\fine

Our next goal is to show that with probability one infinitely often the sets $A_n$ have at most $10 K$ gaps.

\begin{proposition}
\begin{equation} \label{upper bound 10K}
P(G(A_n) \leq 10K \mbox{ i.o.})=1.
\end{equation}
\end{proposition}
{\bf Proof.} Let $\rho=\inf\{n \geq 0: G(A_n) \leq 10K \}.$ To establish
(\ref{upper bound 10K}) it suffices to show that 
\begin{equation} \label{ineq in upper bound 10K}
P^a(\rho < \infty)=1
\end{equation}
for all finite $a \subset {\bf Z}.$ Now note that for any $n \in {\bf N}$ 
\begin{eqnarray*}
P^a(\rho < \infty) &=& \sum_{b} P^a(\rho < \infty \mid A_n=b) P^a(A_n=b) \\
&\geq& \sum_{b} P^b(\rho < \infty)  P^a(A_n=b), 
\end{eqnarray*}
where the sum is taken over all sets $b$ with $P^a(A_n=b) > 0.$ Clearly
$P^b(\rho < \infty)=1$ for all sets $b$ with $G(b) \leq 10 K$ and
for any set $b$ with $P^a(A_n=b) > 0$ we have $|b|=|a|+n.$  Thus, if we can show that   
$P^b(\rho < \infty)=1$
for all sufficiently large finite sets $b \subset {\bf Z}$ with $G(b) > 10K,$
then (\ref{ineq in upper bound 10K}) follows by choosing $n$ sufficiently large in the above decomposition.

So let $b$ be a finite set of integers with $G(b) \geq 10 K.$  Let ${\cal F}_n$ denote the $\sigma$-algebra generated by the events $\{A_m=\tilde a\},$ $m=0, \ldots, n,$ ${\tilde a} \subset {\bf Z}$ with $|{\tilde a}| < \infty.$ Our first goal is to  show that, $\{L_{n \wedge \rho},{\cal F}_n
\}$ is a supermartingale under $P^b,$ if the set $b$ is sufficiently large.
  
Let $R_{\alpha}^n$ denote the number of particles which are added during
period $n$ at empty gap sites and let $R_{\delta}^n$ denote the number of
deletions in which we remove an interior site of the present set, i.e., if
$X_{n,i}$ denotes the point which gets added to (or removed from) the set
$A_{n,i-1},$ then $R_{\alpha}^n=\sum_{i=1}^K {\bf 1}_{ \{X_{n,i} \mbox{ belongs to
a gap of } A_{n,i-1}\} }$ and $R_{\delta}^n=\sum_{i=K+1}^\kappa {\bf 1}_{ \{X_{n,i}
\mbox{ is an interior site of } A_{n,i-1} \} }.$ Note that 
\begin{equation} \label{prop "10 K" ineq L vs. R}
L_{n+1} - L_n \leq R_{\delta}^{n+1} - R_{\alpha}^{n+1},
\end{equation}
since $L$ decreases in an addition step if and only if the new site belongs to a gap of the present set and in a deletion step $L$ can only increase when a particle is deleted at an interior site of the present set.
Now if ${\tilde b}$ is a finite set of integers with $G({\tilde b}) \geq 10 K$ then for $1
\leq i \leq K-1$ we have
\begin{eqnarray} \label{prop "10K", first est Rd-Ra}
\lefteqn{P^b(R_{\delta}^{n+1} - R_{\alpha}^{n+1} = i \mid A_n={\tilde b}) } \nonumber \\
   &\leq&  P^b(R_{\delta}^{n+1}=K-1,R_{\alpha}^{n+1}=K-(i+1) \mid A_n={\tilde b}) \nonumber \\
&& \quad + \; P^b(R_{\delta}^{n+1} < K-1 \mid A_n={\tilde b}) \nonumber \\
  & \leq & P^b(K-R_{\alpha}^{n+1}=i+1 \mid A_n={\tilde b}) 
        + P^b(K-1-R_{\delta}^{n+1} \geq 1 \mid A_n={\tilde b}).
\end{eqnarray}
Note that $K-R_{\alpha}^{n+1}$ equals the number of steps in period $n+1$ in
which the new site gets chosen from the outer boundary of
the present set and since ${\tilde b}$ has at least $10K$ gaps each of the sets 
$A_{n+1,i}(\omega)$, $0 \leq i \leq K-1$ has at least $10K+2 - i \geq 9K+3$
boundary sites. Thus
\begin{eqnarray*}
P^b(K-R_{\alpha}^{n+1}=i+1 \mid A_n={\tilde b}) &\leq& \left(\begin{array}{c}
                                                    K  \\
                                                   i+1 
                                           \end{array}
                                     \right)
                                     \left( \frac{2}{9K+3} \right)^{i+1}  
               \leq \frac{K^{i+1}}{(i+1)!} \left(\frac{2}{9K} \right)^{i+1}\\
                      &=& \frac {1}{(i+1)!} (2/9)^{i+1}.
\end{eqnarray*}
Similarly note that $K-1-R_{\delta}^{n+1}$ equals the number of steps in period
$n+1$ in which we remove an endpoint of the present set and for $i=K,
\ldots, \kappa -2$ the set $A_{n+1,i}$ has at least size $|A_n|+2.$ Thus
\[P^b(K-1-R_{\delta}^{n+1} \geq 1 \mid A_n={\tilde b}) \leq (K-1) \; 2/(|{\tilde b}|+2).\]
Using the last two estimates in (\ref{prop "10K", first est Rd-Ra}) we get 
\begin{eqnarray} \label{prop "10K", second est Rd-Ra}
P^b(R_{\delta}^{n+1} - R_{\alpha}^{n+1} = i \mid A_n={\tilde b}) &\leq& \frac {1}{(i+1)!}
(2/9)^{i+1} + (K-1) 2/{|{\tilde b}|} \nonumber \\
& \leq & \frac {1}{(i+1)!}
4^{-(i+1)}, 
\end{eqnarray}
if the set ${\tilde b}$ is sufficiently large, say $|{\tilde b}| \geq N_1$ where $N_1$ depends only on $K.$ 
On the other hand we have
\begin{eqnarray*} 
\lefteqn{ P^b(R_{\delta}^{n+1} - R_{\alpha}^{n+1} \leq -1 \mid A_n={\tilde b})} \\
&\geq& P^b(R_{\alpha}^{n+1}=K, R_{\delta}^{n+1}=K-1 \mid A_n={\tilde b}) \\
&=& P^b(R_{\delta}^{n+1}=K-1 \mid R_{\alpha}^{n+1}=K,
A_n={\tilde b}) P^b(R_{\alpha}^{n+1}=K \mid A_n={\tilde b})  \\
& \geq & 
\left(1-\frac{2}{|{\tilde b}|+2} \right)^{K-1} \; \prod_{i=0}^{K-1} \frac{10K-i}{10K-i+2} \;\; ,
\end{eqnarray*} 
since for $\omega \in \{A_n={\tilde b}\}$ the size of the boundary of $A_{n+1,i}(\omega)$ is bounded below by $10K+2-i$ for $i=0, \ldots, K-1,$ and the size of $A_{n+1,i}(\omega)$ is bounded below by $|{\tilde b}|+2$ for $i=K, \ldots, \kappa-1.$ Now the rightmost product above equals $(9K+2)(9K+1)/((10K+2)(10K+1))$ and thus is bounded below by $(9/10)^2=.81.$ Thus for all sufficiently large sets ${\tilde b},$ say $|{\tilde b}| \geq N_2,$ we get 
\begin{equation} \label{prop "10K", third est Rd-Ra}
P^b(R_{\delta}^{n+1} - R_{\alpha}^{n+1} \leq -1 \mid A_n={\tilde b}) \geq 0.8 .
\end{equation}
Now any set ${\tilde b}$ with $P^b(A_n={\tilde b}) > 0$ has size at least $|b|.$ Thus, letting $N_0=N_1 \vee N_2,$ we can conclude from (\ref{prop "10K", first est Rd-Ra}), (\ref{prop "10K", second est Rd-Ra}) and (\ref{prop "10K", third est Rd-Ra}) that, if $|b| \geq N_0,$ then 
\begin{eqnarray*}
E^b(L_{n+1}-L_n \mid A_n &=& {\tilde b}) \leq \sum_{i=1}^{K-1}  \frac{i}{(i+1)!} 4^{-(i+1)} - 0.8 \leq 4^{-2} \sum_{i=0}^{\infty} 4^{-i} - 0.8 \\
&=& 1/12 - 0.8 < 0 
\end{eqnarray*}
for all finite ${\tilde b} \subset {\bf Z}$ with $G(\tilde b) \geq 10K$ and $P^b(A_n={\tilde b}) > 0.$ Hence $\{L_{n \wedge \rho},{\cal F}_n
\}$ is a supermartingale under $P^b.$ (Here we also use the fact that $(L_{(n+1) \wedge \rho}-L_{n \wedge \rho})(\omega) = 0$  for $\omega \in \{G(A_n)<10K\}$).

Using that  $L_{n \wedge \rho}$ is nonnegative, we can conclude moreover that the supermartingale converges almost surely and thus, since  $L_{n \wedge \rho}$ is integer valued, $L_{n \wedge \rho}$ is almost surely  eventually constant, i.e.,
\begin{equation} \label{prop "10K", stopped L eventually constant}
P^b(\bigcup_{i \in {\bf N}_0} \{L_{n \wedge \rho} =i \mbox{ eventually} \})=1.
\end{equation} 
Yet, for any $i \in {\bf N}_0$ we have 
\begin{equation} \label{prop "10K", L eventually constant}
  P^b(L_n=i \mbox{ eventually }) = 0. 
\end{equation}
To see this it suffices to show that there exists a constant $c>0$ with
$ P^b(L_{n+1} \neq L_n \mid A_n) \geq c$
for all $n \in {\bf N}_0.$ Now, for any set ${\tilde a}
\subset {\bf Z}$ with $P^b(A_n={\tilde a})>0$ we have in the case $G({\tilde a}) \geq 10K,$
$P^b(L_{n+1} < L_n \mid A_n={\tilde a})$ $\geq P(R_{\delta}^{n+1} - R_{\alpha}^{n+1} \leq -1 \mid
A_n={\tilde a})$ $\geq 0.8$ by (\ref{prop "10K", third est Rd-Ra}) and in the case $G({\tilde a}) < 10K$ we have  $P^b(L_{n+1} > L_n \mid
A_n={\tilde a})$ $\geq (2/(2G({\tilde a})+2))^K (1-2/(|{\tilde a}|+2))^{K-1}$ $ \geq  (2/(20K+2))^K (1-2/N_0)^{K-1},$ since $L$ clearly
increases if in period $n+1$ all new particles get added at sites belonging
to the outer boundary of the present set and if all particles which get
removed are at interior sites of $A_{n +1,K}.$ Thus (\ref{prop "10K", L eventually constant}) follows.

Combining (\ref{prop "10K", stopped L eventually constant}) and (\ref{prop "10K", L eventually constant}) we can finally conclude that for all finite sets $b \subset {\bf Z}$ with $|b| \geq N_0$ and $G(b) \geq 10K,$
$P^b(\bigcup_{i \in {\bf N}_0} \{L_{n \wedge \rho}=i \mbox{ eventually and } \rho < \infty\} ) = 1, $
and thus 
$ P^b(\rho < \infty) = 1. $
\newline
\vspace{7 mm}
\fine

Let ${\cal G}$ denote the following family of finite subsets of ${\bf Z}$
\begin{equation}
  {\cal G} = \{ a \subset {\bf Z}: |a| < \infty, G(a) \leq 10K \mbox{ and
} L(a) \leq (K-1) |a| \}. 
\end{equation}

\begin{corollary} \label{cor: recurrence of G}
\qquad \qquad $P(A_n \in {\cal G} \; i.o.) = 1. $
\end{corollary}
{\bf Proof.} Since in each step we add $K$ and delete $K-1$ particles it is
clear that any set $a$ with $P(A_n = a) > 0$ satisfies $|a| <
\infty$ and $L(a) \leq (K-1) n \leq (K-1) |a|.$ Thus the recurrence of
${\cal G}$ follows immediately from the previous proposition.
\newline
\vspace{7 mm}
\fine

In the next proposition we will show that with probability 1 eventually the 
number of gap sites is bigger than $K-3.$

\begin{proposition} \label{prop: lower bound for lim inf of L}
\[ P(L_n \geq K-2 \mbox{ eventually }) = 1. \] 
\end{proposition}
{\bf Proof.} For $a \subset {\bf Z},$ $2 < |a| < \infty,$ let $\alpha(a)$ (resp. $\beta(a)$) denote the second smallest (resp. second biggest) element of $a.$ Define the event ${\cal E}_n=\{$  
there exist at least two indices $i, K < i \leq \kappa$, $\mbox{ with } X_{n,i} \leq \alpha(A_{n,K}) \mbox{ or } X_{n,i} \geq \beta(A_{n,K})  \}.$ Since for any $K < i \leq \kappa,$ $A_{n,i-1}$ is a subset of $A_{n,K},$ we have $| \{ x \in A_{n,i-1}: x \leq \alpha(A_{n,K}) \mbox { or } x \geq \beta(A_{n,K}) \}| \leq 4.$ 
Hence $P(X_{n,i} \leq \alpha(A_{n,K}) \mbox{ or } X_{n,i} \geq \beta(A_{n,K})) \leq 4/{|A_{n,i-1}|} \leq 4/n$ and so $P({\cal E}_n) \leq 2^{-1}K(K-1)\, (4/n)^2.$ By the Borel-Cantelli lemma we can conclude that 
\[ P( \limsup {\cal E}_n) = 0. \]
Now let $\omega \in {\cal E}_n^c.$ Then all but at most one index $i$  
satisfy $\alpha(A_{n,K})(\omega) < X_{n,i}(\omega) < \beta(A_{n,K})(\omega).$ So at least one of the sites $\alpha(A_{n,K})(\omega),$ $\min A_{n,K}(\omega)$ remains occupied and, similarly,  at least one of the sites $\beta(A_{n,K})(\omega), \max A_{n,K}(\omega)$ remains occupied.
Hence $\min A_{n,\kappa}(\omega) \leq \alpha(A_{n,K})(\omega)$ and $\max A_{n,\kappa}(\omega) \geq 
\beta (A_{n,K})(\omega).$ In particular, all but at most one of the sites removed lies between the endpoints of $A_{n,\kappa}(\omega)$ and so $L_n(\omega) \geq K-2.$
Thus we have shown that ${\cal E}_n^c \subset \{L_n \geq K-2 \}$ and so $P(L_n \geq K-2  \mbox{ eventually} ) \geq P( \liminf {\cal E}_n^c)=1.$
\newline
\vspace{7 mm}
\fine

Before proceeding further let us state some useful properties of the
function $h,$ the proof of which is straightforward and will thus be omitted. Recall that $h(x)=\log x/\log\log x.$ 
\begin{proposition} \label{properties of h}
i) $h(2x)/h(x) \rightarrow 1 \mbox{ as } x \rightarrow \infty.$ \newline
ii) Let $c > 0.$ For any $\eta > 0$ there exists $x_0 > 0$ such that for all $x \geq x_0$
\[ x^{1-\eta} \leq (c \, h(x))^{h(x)} \leq x .\]
\end{proposition}

Using these properties we can prove the lower bound in part (1)  of Theorem \ref{first theorem}. 

\begin{proposition} \label{lower bound for limsup of L(A)/h}
\[ P(\limsup_{n \rightarrow \infty} L(A_n)/h(n) \geq (K-1)/K) = 1. \]
\end{proposition}
{\bf Proof.} We need to show that for every $0<c<(K-1)/K,$
\begin{equation} \label{first goal in lower bound for limsup}
P(L_n \geq c \, h(n) \; i.o.  ) = 1.
\end{equation}
So fix $c<(K-1)/K$ and fix a positive integer $N_0 $ to be chosen
later. Let $n_i=2^i N_0 - 1,$ $ i \in \Nzero.$ Thus $|A_{n_i}|=2^i N_0.$
Define the time interval $I_i=(n_{i-1},n_i],$ $i \in {\bf N},$ and
let ${\cal L}_i=\{\exists j \in I_i \mbox{ with } L_j \geq c \, h(j) \}.$ By
the Markov property of the process $A,$ (\ref{first goal in lower bound for limsup}) will be proved when we
show that 
\[P({\cal L}_i \mid A_{n_{i-1}}) \geq 1/2\]
for all $i \in {\bf N}.$ Now note that for $j \in I_i$ and $N_0$
sufficiently large $h(j) \leq h(n_i) \leq h(2^i N_0) = h(2 |A_{n_{i-1}}|).$
So, if $b$ is a set with $P(A_{n_{i-1}}=b) >0$ then 
\begin{eqnarray*}
\lefteqn{P({\cal L}_i \mid A_{n_{i-1}}=b)} \\
& \geq &P(\exists j \in I_i \mbox { with } L_j \geq c \, h(2 |A_{n_{i-1}}| ) \mid  A_{n_{i-1}}=b) \\
&=& P^b(\exists 0 < j \leq |b| \mbox{ with } L_j \geq c \, h(2 |b|)). 
\end{eqnarray*}
Thus it suffices to show that for all sufficiently large finite sets $b \subset {\bf Z}$  
\begin{equation} \label{second goal in lower bound for limsup}
P^b(\exists 0 < j \leq |b| \mbox{ with } L_j \geq c \, h(2 |b|)) \geq 1/2.
\end{equation}     
So let $b$ be a finite set of integers and let $N=|b|.$ Choose $\varepsilon > 0$ with
$(1-\varepsilon)(K-1)/K \geq c$ and let $m_0=m_0(N)=
\lfloor{(1-\varepsilon/2) h(N)/K} \rfloor.$ Each of the intervals
$((j-1)m_0, j m_0],$ $j=1, \ldots, \lfloor{N/m_0}\rfloor,$ is contained in the time interval $(0,N].$ Let 
$\tau_N= \inf \{ n \geq 0: L_n \geq c \, h(2N) \}.$ Our next goal is to show that if $N$ is sufficiently large, then for any finite set $a \subset {\bf Z}$ with $|a| \geq N$ we have  
\begin{equation} \label{third goal in lower bound for limsup}
P^a(\tau_N \leq m_0(N)) \geq N^{-(1-\varepsilon/2)}. 
\end{equation}
Once this inequality gets established we can conclude that if $N=|b|$ is sufficiently large then  
\[P^b(\tau_N > N) \leq
(1-N^{-(1-\varepsilon/2)})^{\lfloor{N/m_0}\rfloor}. \]
Using that $\lfloor{N/m_0}\rfloor \geq
N^{1-\varepsilon/4}$ (for all sufficiently large $N$) the right hand side is bounded above by $(1-N^{-(1-\varepsilon/2)})^{N^{1-\varepsilon/4}}$ and this last expression converges to $0$ as $N \rightarrow \infty.$ Thus the inequality (\ref{second goal in lower bound for limsup}) follows.

So let us now work on establishing (\ref{third goal in lower bound for limsup}). Note that $L$ increases by
$K-1$ in a given period if all new particles attach at the outer
boundary and if only particles at interior sites get deleted. Thus for any
finite set ${\tilde a} \subset {\bf Z}$ with $|{\tilde a}| \geq N $ and $L({\tilde a}) < c h(2 N)$ we
have
$P^{\tilde a}(L_1-L_0 = K-1) \geq (2/(2+L({\tilde a})))^K (1-2/|{\tilde a}|)^{K-1} \geq (2/(c
h(2N)))^K (1-2/N)^{K-1}$ 
and thus by part (i) of Proposition \ref{properties of h}
we get
\[
P^{\tilde a}(L_1-L_0 = K-1) \geq h(N)^{-K},
\]
if $N$ is sufficiently large. Next note that if within each of the first
$m_0(N)$ periods $L$ increases by $K-1$ then the total increase of $L$ is
at least $(K-1)m_0=(K-1)\lfloor{(1-\varepsilon/2) h(N)/K} \rfloor \geq
(K-1)/K \; (1-\varepsilon) h(2N)$ for $N$ sufficiently large. (Here once again we used part (i) of Proposition \ref{properties of h}). Now by choice of $\varepsilon$ this last term can be bounded below by
$c h(2N).$ Thus, if $N$ is sufficiently large, then for any finite set of
integers $a$ with $|a| \geq N$ we have 
\[P^a(\tau_N \leq m_0(N)) 
\geq P^a(L_i-L_{i-1} = K-1 \mbox{ for } i=1, \ldots, \tau_N \wedge m_0(N)) 
\geq h(N)^{-K m_0}.\]
Using that $K m_0 \leq (1- \varepsilon/2) h(N)$ and $h(N)^{-h(N)}
\geq N^{-1}$ by part (ii) of Proposition \ref{properties of h}, we finally get (\ref{third goal in lower bound for limsup}). 
\newline
\vspace{7 mm}
\fine

\section{ New and Old Gaps}

So far we have shown that infinitely often the number of gaps is bounded
above by $10K.$ But gaps may be large and so this result does not yet give
us any upper bound on $L.$ Our next goal will be to show that infinitely
often the set $A_n$ has at most one gap of size two and no bigger gaps, and
for any fixed choice of $M>0$ infinitely often the size of $L(A_n)/h(n)$ will be smaller than $1/M.$ So let $M \in {\bf N}.$ For any finite set of integers $a$ let $C(a)$ (resp. $C^{(2)}(a)$) denote the set of unoccupied   interior sites of $a$ which are within distance $4$ to at least one (resp. two) other unoccupied interior sites of $a$. 
Define
\[ {\cal S}_M=\{ a \subset {\bf Z} : |a| < \infty, L(a) \leq h(|a|)/M,
|C(a)| \leq 3 \mbox{ and } |C^{(2)}(a)| \leq 1. \} . \]
We want to show that the family ${\cal S}_M$ of subsets of ${\bf Z}$ is recurrent for the process $A.$ 
Recall that 
\[ L_{n} - L_{n-1} \leq R_{\delta}^{n} - R_{\alpha}^{n}, \]
where $R_{\alpha}^{n}$ denotes the number of particles which get added in
period $n$ at empty gap sites and $R_{\delta}^{n}$ denotes the number of
deletions in which we remove a particle at an interior site of the present set. Thus, letting $X_{n,i}$ denote the point which gets added to (resp. deleted from) the set $A_{n,i-1},$  we get
\begin{eqnarray*}
L_{n+1} - L_n &\leq& (K-1) - \sum_{i=1}^{K} {\bf 1}_{ \{ X_{n,i} \mbox{ belongs to a gap of $A_{n,i-1}$ } \} } \\
   &=& (K-\sum_{i=1}^{K} {\bf 1}_{ \{ X_{n,i} \mbox{ belongs to a gap of $A_{n,i-1}$ } \} } ) -1 \\
   &=& \sum_{i=1}^{K} {\bf 1}_{ \{ X_{n,i} \in \bar{\partial} A_{n,i-1} \} } -1.
\end{eqnarray*}
Clearly $P(X_{n,i} \in \bar{\partial}A_{n,i} \mid A_{n,i-1}) \leq 2/(2+G(A_{n,i-1}))$ and thus $L$ tends to decrease whenever $G$ is large. So if $L$ and $G$ are not too far apart, then large values of $L$ correspond to large values of $G$ and thus $L$ tends to decrease whenever it gets big. 

Now, why may we hope that $L$ and $G$ will get close together? Let us first observe that if $G(A_n)$ is small compared to the size of $A_n$ then typically the removal of a particle creates a new gap of size 1, since with high probability a particle will get removed at a site which is not adjacent to any existing gap. Next let us suppose a big gap $g$ has formed at time $n_0$ and let us investigate how this gap will change over time. For any $n \geq n_0,$ $1 \leq i \leq K,$ we have $P(X_{n,i} \in g \mid A_{n,i-1}) \geq 1/(2+2G(A_{n,i-1}))$ and for $K < i \leq \kappa$ we have $P(X_{n,i} \mbox{ is adjacent to } g \mid A_{n,i-1}) \leq 2/|A_{n,i-1}|.$ Thus as long as $G$ is small compared to the size of $A$ 
we expect that typically several new particles will fill into this gap before any deletion will increase its length. Yet, the problem is, that we have no good upper bound for the size of a possible increase. Just note that whenever we remove a particle which separates two gaps, these two gaps will form one single enlarged gap the size of which may be considerably larger than the size of the gap $g$ we are considering.

At this point, the idea of distinguishing between `old' and `new' gaps will turn out to be useful. We will be able to show that if we wait for a certain time, chances are high that all `old' gaps will be filled in and we are left with only `new' gaps, each of which is still small. So let us define a Markov process $\{(A,I)_{n,i}: n \in {\bf N} , i=0, \ldots , \kappa \}$ where $I_{n,i} \subset \{1, \ldots, G(A_{n,i} ) \}$ will be interpreted as the set of indices of `new' gaps.

\vspace{7 mm}

For any finite $a \subset {\bf Z}$ let $g_1(a), \ldots ,g_{G(a)}(a)$ denote the gaps of $a$, these gaps being  enumerated from left to right, i.e.,  if $x \in g_i(a)$ and $y \in g_j(a)$ with $i<j$ then $x<y.$  
Let $A_{1,0}$ be a given finite set of integers and let $I_{1,0} \subset \{ 1, \ldots, G(A_{1,0}) \}$ be a given set of indices of `new' gaps for $A_{1,0}.$ 
Given $(A_{n,i},I_{n,i})$ with $n \geq 0,$ $ 0 \leq i < K,$ let $A_{n,i+1}$ be obtained as before by adding a particle at the boundary of $A_{n,i},$ the new site $X_{n,i+1}$ being chosen uniformly from all boundary sites. Define
\[ I_{n,i+1} = \{ 1 \leq j \leq G(A_{n,i+1}): \exists k \in I_{n,i} \mbox{ with } g_j(A_{n,i+1}) \subset g_k(A_{n,i}) \}. \]
Thus, if $X_{n,i+1} \in g_k(A_{n,i})$ and if $|g_k(A_{n,i})| >1,$ (i.e., the new particle does not completely fill the gap $g_k(A_{n,i})$ ) then $g_k(A_{n,i}) \setminus \{X_{n,i+1} \}$ is a gap of $A_{n,i+1}$ and it is considered as being `new' (resp. `old') for $A_{n,i+1}$ if and only if $g_k(A_{n,i})$ was `new' (resp. `old') for $A_{n,i}.$ 
Similarly, given $(A_{n,i},I_{n,i})$ with $n \geq 0,$ $ K \leq i < \kappa,$ let $A_{n,i+1}$ be obtained by deleting a particle from the set $A_{n,i},$ the site of this particle being chosen uniformly from all sites of $A_{n,i}.$ Define
\begin{eqnarray*}
 I_{n,i+1} &=& \{ 1 \leq j \leq G(A_{n,i+1}): \mbox{ there exists no index $k \in \{1, \ldots, G(A_{n,i}) \} \setminus I_{n,i}$ with } \\
  &  & \mbox{ $g_k(A_{n,i}) \subset g_j(A_{n,i+1})$ }  \}. 
\end{eqnarray*}
Thus, if $X_{n,i+1}$ is not an endpoint of $A_{n,i},$ then the gap $g$ of $A_{n,i+1}$ which contains the site $X_{n,i+1}$ is `old' for $A_{n,i+1}$ if and only if $X_{n,i+1}$ is adjacent to an `old' gap of $A_{n,i},$ otherwise $g$ is `new' for $A_{n,i+1}.$ Note in particular that if $X_{n,i+1}$ is adjacent to both an `old' gap $g$ and a `new' gap ${\tilde g}$ of $A_{n,i},$ then the gap $g \cup \{X_{n,i+1}\} \cup {\tilde g}$ of $A_{n,i+1}$ is `old' for $A_{n,i+1}.$ 
To complete the definition let $A_{n+1,0}=A_{n,\kappa}$ and $I_{n+1,0}=I_{n, \kappa}.$ 
The corresponding Markov process of full periods will be denoted by $(A_n,I_n)_{n \geq 0},$ with the convention that $(A_0,I_0)=(A_{1,0},I_{1,0})$ and $(A_n,I_n)=(A_{n,\kappa},I_{n,\kappa})$  for all $n \in {\bf N}.$

We will write $P^{a,\lambda}$ when we consider the Markov process $(A,I)_n$ 
or $(A,I)_{n,i}$ with initial pair $(a, \lambda).$   

A site of ${\bf Z} \setminus A_{n,i}$ which belongs to a `new' gap of $A_{n,i}$ will be called a new gap site for $A_{n,i}$ and similarly we define old gap sites. For any finite set $a$ of integers and $\lambda \subset \{1, \ldots, G(a) \}$ define $L^I(a,\lambda)=\sum_{i \in \lambda} |g_i(a)|$ and let $L_{n,i}^I=L^I(A_{n,i},I_{n,i}).$ Thus $L_{n,i}^I$ counts the number of unoccupied sites belonging to `new' gaps of $A_{n,i}.$ Let $G_{n,i}^I$ denote the number of `new' gaps of $A_{n,i}$ and let $G_{n,i}^{2,I}$ denote the number of `new' gaps of size at least two. Finally let $C_{n,i}^I$ (resp. $C_{n,i}^{I,2}$) denote the set of new gap sites for $A_{n,i}$ which are within distance 4 of at least one (resp. two) other new gap sites of $A_{n,i}.$ Let $L_{n}^I, G_n^I,$ etc. be defined in the obvious way. Whenever we drop the superscript $I$ we consider the set of all gap sites instead of the set of new gap sites. 

\vspace{7 mm}

In our first lemma we will show that if we start with a sufficiently big set $a$ belonging to ${\cal G}$ (all the gaps of which are considered as being old) then with overwhelming probability within the next $|a|^{6/5}$ periods the following holds: (i) all old gaps will get filled in, (ii) during all of these periods  the number of new gap sites will be bounded above by $(1+(K-1)/K) h(|a|),$ and (iii) all new gaps (with the  exception of possibly one) will have size 1. Moreover, for any fixed constant $M>0$ there exists a period $m,$ $|a|^{11/10} < m <|a|^{6/5},$ such that at the start of this period the number of new gap sites is smaller than $h(|a|)/M.$ 

Now, why should this hold? Roughly speaking, as long as all (with the  exception of possibly one) new gaps have size 1, the number of old gap sites  (as long as at least one old gap is present) can be compared to a suitable random walk which `tends to decrease'. Under the same assumption the number of new gap sites can be compared to a similar random walk at those times where at least a certain minimum number of new gap sites is present. Finally, as long as the number of new gap sites is not `too large', new gaps tend to get filled in `pretty soon' after they are formed, and, on the other hand, the chance of selecting a particle for deletion which is close to an existing gap (and hence the chance of enlarging an existing gap) is very small. Thus during the first $|a|^{6/5}$ periods it is very unlikely that a particle gets deleted close to an existing new gap at a time, at which there still exists a new gap of size bigger than $1.$

At this point the reader may get very suspicious about describing a logical circle. The next proposition will provide one of the main tools to overcome this difficulty. This proposition will later enable us to compare parts of the process we are interested in (for instance `excursions' of $L$) to the corresponding parts of some other, `nicer' process (for instance a `K-decrease' of a random walk) in a path-by-path way.

\vspace{7 mm}

Let ${\bf N}^2_{\infty}=\{(i,j): i \in {\bf N}_0 \cup \{ \infty \}, j \in {\bf N}_0 \cup \{ \infty \} \}.$ For $x,y \in {\bf N}^2_{\infty}$ define $x \leq y$ if $x_1 \leq y_1$ and $x_2 \leq y_2.$ Let $(\Omega, {\cal F}, P)$ be a given probability space. A family $\{ {\cal F}_{x}, x \in {\bf N}^2_{\infty}\}$ of sub-$\sigma$-fields of ${\cal F}$ is called a filtration in $\Omega$ if ${\cal F}_{x} \subset {\cal F}_{y}$ for  $x \leq y.$ Given a filtration $\{ {\cal F}_{x}, x \in {\bf N}^2_{\infty}\},$ a random vector $\lambda$ with values in ${\bf N}^2_{\infty}$ is called a stopping time (relative to  
$\{ {\cal F}_{x}, x \in {\bf N}^2_{\infty}\}$), if $\{ \lambda \leq x\} \in  {\cal F}_{x}$ for all $x \in {\bf N}^2_{\infty}.$ For any stopping time $\lambda$ the pre-$\lambda$ field ${\cal F}_{\lambda}$ is defined as $\{\Lambda \in {\cal F}: \Lambda \cap \{\lambda \leq x\} \in {\cal F}_x \mbox{ for all } x \in {\bf N}^2_{\infty}\}.$

\begin{proposition} \label{prop:good selection rules}
Let $F$ be a distribution function and let $\{U_n^{(j)}: n \in {\bf N}, j=1,2
\}$ be a family of independent random variables with $U_n^{(j)} \sim F$ for all $n \in 
{\bf N}, j=1,2.$ Let ${\cal F}_{(i,j)}=\sigma(U_1^{(1)} ,\ldots, U_i^{(1)},U_1^{(2)}
,\ldots, U_j^{(2)})$ (with the convention that ${\cal F}_{(0,0)}=\{\emptyset, \Omega\},$ ${\cal F}_{(\infty,\infty)}=\sigma(U^{(1)},U^{(2)}),$ ${\cal F}_{(0,j)}=\sigma(U^{(2)}_1, \ldots, U^{(2)}_j),$ ${\cal F}_{(\infty,j)}=\sigma(U^{(1)}) \vee \sigma(U^{(2)}_1, \ldots, U^{(2)}_j)$ and similarly for ${\cal F}_{(i,0)}$ and ${\cal F}_{(i,\infty)}$). 
Let $\{\lambda_i\}$ be a family of
stopping times with respect to the filtration $\{{\cal F}_{(i,j)}:(i,j) \in {\bf N}^2_{\infty}\}$. Let
$B_i \in {\cal F}_{\lambda_i}$ and assume that the sequence $\{\lambda_i\}$
satisfies $\{\lambda_{i,1} < \infty\} \subset B_i,$ $\{\lambda_{i,2} < \infty\} \subset B_i^c$ and 
\begin{equation} \label{first property of the selection rule}
{\bf 1}_{B_i}(\lambda_{i,1}+1,\lambda_{i,2})+{\bf
1}_{B_i^c}(\lambda_{i,1},\lambda_{i,2}+1) \leq \lambda_{i+1} 
\end{equation}
for all $i$. Let 
\begin{equation} \label{second property of the selection rule}
\xi_i={\bf 1}_{B_i}U_{\lambda_{i,1}+1}^{(1)}+{\bf 1}_{B_i^c}U_{\lambda_{i,2}+1}^{(2)}.
\end{equation} 
Then $\{\xi_i: i \in {\bf N}
\}$ is a family of independent random variables and $\xi_i \sim F$ for all $i \in {\bf N}.$ 
\end{proposition}
(The reader may think of two independent series of coin tossings. A third series gets constructed from the first two such that the selection rule (which decides when and from which series to select the next element for the new series) depends only on previous coin tossings. Then the new series is
-in a probabilistic sense - no different than the two original ones.)
\vspace{7 mm}

{\bf Proof of the proposition.} Using (\ref{second  property of the selection rule}) and the independence of the random variables $U_n^j$
we have
\begin{eqnarray*}
P(\xi_i \leq t \mid {\cal F}_{\lambda_i}) &=& {\bf 1}_{B_i}P(U_{\lambda_{i,1}+1}^{(1)} \leq t \mid {\cal F}_{\lambda_i}) + {\bf 1}_{B_i^c}P(U_{\lambda_{i,2}+1}^{(2)} \leq t \mid {\cal F}_{\lambda_i}) \\
 &=& {\bf 1}_{B_i} F(t) + {\bf 1}_{B_i^c} F(t) 
= F(t).  
\end{eqnarray*}
Thus $\xi_i$ is independent of ${\cal F}_{\lambda_i}$ and $\xi_i \sim F.$ Since (\ref{first property of the selection rule}) together with (\ref{second property of the selection rule}) implies that for any $i \geq 2$ the random variables $\xi_1, \ldots ,\xi_{i-1}$ are ${\cal F}_{\lambda_i}$-measurable, our proof is completed.
\newline 
\vspace{7 mm}
\fine

Before actually using such an algorithmic scheme in the proof of the next two lemmas, let us state a technical estimate which will be used when investigating a particular random walk occuring in our construction. This estimate will in fact be crucial when proving that with probability 1 we have $\limsup L(A_n)/h(n) \geq (K-1)/K.$ 

\begin{proposition} \label{prop: random walk estimates}
Let $c,N>0$ with $c/h(N) < 1$ and let $Z$ be a random walk starting from $0$ with step law $\mu$ defined by $\mu(1)=c/h(N)$ and $\mu(0)=1-c/h(N).$ Let $\eta > 0$ and define  $S=\inf \{m>0: Z_{mK}-m \geq (\eta+(K-1)/K) h(N) \mbox{ or } Z_{mK}-m \leq -K \}.$  Then for every $\varepsilon > 0$ there exists $N_0$ such that for all $N \geq N_0$ we have
\begin{equation} \label{ineq in the random walk proposition}
P(Z_{SK}-S \geq (\eta+(K-1)/K) h(N)) \leq N^{-(1 +\eta K/(K-1)-\varepsilon)}.
\end{equation}
\end{proposition}
{\bf Proof.} Let $\beta = (K-1)/K+\eta.$ Then 
\[ P(Z_{SK}-S \geq \beta h(N))=\sum_{m \geq 0}P(S=m, Z_{mK} \geq m+ \beta h(N)).\]
Now, if $S(\omega)=m$ and $Z_{mK}(\omega) \geq m+ \beta h(N))$ then, using the obvious bound $Z_{mK} \leq mK,$ we get $mK \geq m + \beta h(N)$ and hence $m \geq (K-1)^{-1} \beta h(N).$ Thus we can restrict the sum on the right hand side to integers $m \geq m_0(N),$ where $m_0(N)= \lceil {(K-1)^{-1} \beta h(N)} \rceil.$ 

Next, if $S(\omega)=m$ and $Z_{mK}(\omega)-m \geq \beta h(N),$ then $Z_{mK}(\omega)=\lceil {m + \beta h(N)} \rceil + i $ for some integer $0 \leq i \leq (K-2).$ (Just recall that if $S(\omega)=m$ then $m$ is the first integer satisfying $Z_{mK}(\omega)-m \geq \beta h(N)$ and hence $\beta h(N)+m \leq Z_{mK}(\omega) \leq Z_{(m-1)K}(\omega) +K < \beta h(N)+ m-1 +K$).  
Therefore
\begin{eqnarray} \label {first ineq}
P(Z_{SK}-S \geq \beta h(N))
\leq \sum_{i=0}^{K-2} \sum_{m \geq m_0(N)} P(Z_{mK} = m+ \lceil{ \beta h(N)} \rceil + i ) \nonumber \\
\leq \sum_{i=0}^{K-2} \sum_{m \geq m_0(N)} 
\left( \begin{array} {cc}
           mK \\
           m+ \lceil{ \beta h(N)} \rceil + i
       \end{array}
\right)
\left(c/h(N)\right)^{m + \lceil{ \beta h(N)} \rceil }, \nonumber
\end{eqnarray}
and so
\begin{eqnarray} \label{sum estimate in technical random walk proposition}
\lefteqn{P(Z_{SK}-S \geq \beta h(N))} \nonumber \\ 
&\leq& g(N) \sum_{i=0}^{K-2} \sum_{l \geq 0} 
\left( \begin{array} {cc}
           (m_0(N) + l)K \\
           m_0(N)+ \lceil{ \beta h(N)} \rceil  +l + i
       \end{array}
\right)
(c/h(N))^{l},
\end{eqnarray}
where $g(N)=(c/h(N))^{m_0(N) + \lceil{ \beta h(N)} \rceil }.$
Now let $r=1+\eta K/(K-1)$ and note that there exists $x$ (depending on $N$), $0 \leq x \leq 1,$ such that  
\[ m_0(N)K =K (K-1)^{-1} \beta h(N) + x K = r h(N) +xK, \]
and there exists $y$ (depending on $N$), $0 \leq y \leq 2$ such that \begin{eqnarray} \label {second numbered equation}
m_0(N)+ \lceil{ \beta h(N)} \rceil  & = & (K-1)^{-1} \beta h(N) + \beta h(N) +y  
\nonumber \\
& = &  
r h(N) + y.
\end{eqnarray}
 Thus 
\begin{eqnarray} \label{eqn for binomial coefficients}
\left( \begin{array} {cc}
           (m_0(N) + l)K \\
            m_0(N)+ \lceil{ \beta h(N)} \rceil +l + i
       \end{array}
\right) &=& 
\left( \begin{array} {cc}
           r h(N) + lK +xK \\
           r h(N) +l + i +y
       \end{array}
\right) \nonumber \\
& = &
\left( \begin{array} {cc}
           r h(N) + lK +xK \\
           l(K-1) + xK - (i + y)
       \end{array}
\right).
\end{eqnarray}
In order to obtain upper bounds for these binomial coefficients we will use that for all integers $n \leq m$ the following estimates hold,
\begin{equation} \label{simple estimates for binomial coefficients}
\left( \begin{array}{c}
           m \\n
       \end{array}
\right) \leq m^n 
\qquad \mbox{and} \qquad  
\left( \begin{array}{c}
           m \\n
       \end{array}
\right) \leq 2^m.
\end{equation}
Now fix $\delta > 0$ and let $l_0=\lfloor {\delta \eta h(N)/{K}} \rfloor.$ For $l \leq l_0$ and sufficiently large $N$ we have
$r h(N) + lK +xK \leq (r + \delta \eta+K/{h(N)}) h(N) \leq  (r + 2 \delta   \eta ) h(N)$ and thus,
by using the first inequality in (\ref{simple estimates for binomial coefficients}),
\begin{eqnarray}
\lefteqn{\left( \begin{array} {cc}
                                 r h(N) + lK +xK \\
                                 l(K-1) + xK - (i + y)
         \end{array}
  \right)
  (c/h(N))^l}  \nonumber \\
& \leq & \bigl( (r + 2 \delta \eta ) h(N) \bigr) ^{l(K-1)+K} \nonumber \\
& \leq & (r + 2 \delta \eta )^K h(N)^{K(l+1)}. \nonumber
\end{eqnarray}
Using (\ref{eqn for binomial coefficients}) and the above estimate we can bound the first part of the sum in (\ref{sum estimate in technical random walk proposition})
for all sufficiently large $N$ as follows, 
\begin{eqnarray} \label{ineq for the first part of the sum}
\lefteqn{\sum_{i=0}^{K-2} \sum_{l = 0}^{l_0} 
\left( \begin{array} {cc}
           (m_0(N) + l)K \\
           m_0(N)+ \lceil{ \beta h(N)} \rceil  +l + i
       \end{array}
\right)
(c/h(N))^{l}} \nonumber \\ 
 & \leq & K (r + 2 \delta \eta)^K  \sum_{l=0}^{l_0}
h(N)^{K(l+1)}  \nonumber \\
& \leq & K (r + 2 \delta \eta)^K  h(N)^{K (l_0+2)}  \nonumber \\
& \leq & K (r + 2 \delta \eta )^K  h(N)^{\delta \eta h(N)+2K} \nonumber \\
& \leq & h(N)^{ h(N) 3 \delta \eta} .
\end{eqnarray} 

Let us now work on the remaining part of the sum in (\ref{sum estimate in technical random walk proposition}). Note that for $l>l_0$ we have $(\delta \eta)^{-1}r \, lK \geq (\delta \eta)^{-1}r \delta \eta h(N) = r h(N)$ and thus $rh(N)+lK+xK \leq (2+(\delta \eta)^{-1} r) lK .$ Hence for $l>l_0$ (and sufficiently large $N$) we get
by using the second inequality in (\ref{simple estimates for binomial coefficients}),
\begin{eqnarray}
 \lefteqn{ \left( \begin{array} {cc}
                                 r h(N) + lK +xK \\
                                 l(K-1) + xK - (i + y)
         \end{array}
  \right)
  (c/h(N))^l } \nonumber \\
& \leq & 2^{(2+(\delta \eta)^{-1}r )lK } \left(c/h(N) \right)^l \nonumber \\
& = & \left(c \; 2^{(2+(\delta \eta)^{-1}r)K }/h(N) \right)^l. \nonumber
\end{eqnarray}
 This last expression can be bounded above by  $2^{-l}/K$ for all sufficiently large $N$. (Just note that $c2^{(2+(\delta \eta)^{-1}r )K }/h(N)$ converges to $0$ as $N \rightarrow \infty$).
So by using this estimate together with (\ref{eqn for binomial coefficients}) we get (for all sufficiently large $N$) the following estimate for the remaining part of the sum in (\ref{sum estimate in technical random walk proposition})
\begin{equation} \label{ineq for the second part of the sum}
\sum_{i=0}^{K-2} \sum_{l = l_0+1}^{\infty} 
\left( \begin{array} {cc}
           (m_0(N) + l)K \\
           m_0(N)+ \lceil{ \beta h(N)} \rceil  +l + i
       \end{array}
\right)
(c/h(N))^{l} 
 \leq \sum_{l = l_0+1}^{\infty} 2^{-l}  \leq 1 .
\end{equation}
Using (\ref{second numbered equation}), (\ref{ineq for the first part of the sum}) and (\ref{ineq for the second part of the sum}) in (\ref{first ineq}) we can conclude that  
\begin{eqnarray*}
P(Z_{SK}-S \geq (\eta+(K-1)/K) h(N)) & \leq & (c/h(N))^{rh(N)}  ( h(N)^{ h(N) 3 \delta \eta} +1)  \\
& \leq & (c/h(N))^{rh(N)} h(N)^{ h(N) 4 \delta \eta}, 
\end{eqnarray*}
if $N$ is sufficiently large.
Now choose $\delta$ such that $4 \delta \eta \leq \varepsilon/2$ and note that by Proposition \ref{properties of h}  we have (for all sufficiently large $N$), $(c/h(N))^{rh(N)} \leq N^{-(r-\varepsilon/2)}$ and $h(N)^{ h(N) 4 \delta \eta} \leq N^{\varepsilon/2}$ . Thus 
\[P(Z_{SK}-S \geq (\eta+(K-1)/K) h(N)) \leq N^{-(r-\varepsilon)}, \]
if $N$ is sufficiently large.
Recalling that $r=1+\eta K/(K-1)$ our claim thus follows.
\newline
\vspace{7 mm}
\fine

In order to state the next two lemmas we will need to define a couple of stopping times. These stopping times will enable us to investigate the evolution of one process while imposing certain restrictions on the evolution of other processes. Lemma \ref{first lemma} together with the recurrence of $\cal{G}$ will give us the recurrence of ${\cal S}_M$ (for any fixed integer $M$) and part (2) of Theorem \ref{second theorem} by using the Markov property of $A$. Once the recurrence of ${\cal S}_M$ gets established,  Lemma \ref{second lemma} will enable us to finish the proof of Theorem \ref{second theorem} and of part (1) of Theorem \ref{first theorem} by a Borel-Cantelli argument.  We will first state the two lemmas plus two corollaries and we will then proceed to proofs. So let us now get started. 

Let 
\[ \tau^1=\inf \{n: L^I_n \geq (K^{-1}(K-1) + 1)h(|A_0|) \} \] 
and let
\[ T^{M,2}=\inf \{n: L_n \geq (K^{-1}(K-1)+ 3 M^{-1} )h(|A_0|) \} .\] 
Thus $\tau^1$ (resp. $T^{M,2}$) catches the first period such that at the end of this period the number of new (resp. all) gap sites exceeds $K^{-1}(K-1) h(|A_0|)$ by at least $h(|A_0|)$ (resp. $3 M^{-1} h(|A_0|)$). 
Let \[ \tau_M = \inf \{ n \geq |A_0|^{11/10}: L^I_{n} \leq h(|A_0|)/M
\} \] and let 
\[ T_{M,2} = \inf \{ n \geq |A_0|: L_n \leq h(|A_0|)/M
\}. \] 
Thus $\tau_M$ (resp. $T_{M,2}$) catches the first
period after period $|A_0|^{11/10}$ (resp. after period $|A_0|$) such that 
at the end of this period the
number of new (resp. all) gaps is `small' compared to $h(|A_0|).$ Let
\[S=\inf \{(n,i): i>K, |C_{n,i-1}| \geq 2 \mbox{ and } X_{n,i} \in
C_{n,i} \}, \]
\[\sigma=\inf \{(n,i): i>K, |C_{n,i-1}^I| \geq 2 \mbox{ and } X_{n,i} \in
C_{n,i}^I \}, \]
and let $\bar{S}$ (resp. $\bar{\sigma}$) denote the period containing $S$ (resp. $\sigma$), i.e. $\bar{S}=\inf \{n: (n,\kappa) \geq S \}$ and $\bar{\sigma}=\inf \{n: (n,\kappa) \geq \sigma \}.$ 
Note that $S$ (resp. $\sigma$) catches the first deletion step in which the following conditions hold: 
i) $X_{n,i}$ is a (new) gap site for $A_{n,i},$  
ii) $X_{n,i}$ is within distance 4 to some other (new) gap site for
$A_{n,i}$ and   
iii) $A_{n,i-1}$ already had at least 2 (new) gap sites which are within
distance 4 to each other. 
Define 
\[ \sigma_0=\inf \{n: n \geq |A_0|^{11/10} \mbox{ and } G_{n}^{I,2}=0
\} ,\]
i.e. $\sigma_0$ denotes the first period after time $|A_0|^{11/10}$ such
that at the end of this period all new gaps have size 1. Finally let 
\[ \theta_0=\inf\{n:L^I_n=L_n\},\]
 i.e., $\theta_0$ denotes the first period such
that at the end of this period there are no more old gaps left.

\begin{lemma} \label{first lemma}
Let $M > K.$ Let $a \in {\cal G}$ and let
$N=|a|.$
Then for all sufficiently large $N$ we have 
\[P^{a, \emptyset}(\tau^1 \wedge \bar{\sigma} >N^{6/5}, \sigma_0 < N^{6/5},
\tau_M < N^{6/5}, \theta_0<N^{11/10}) \geq 1-N^{-1/30}.\]
\end{lemma} 
The above lemma implis part (2) of Theorem \ref{second theorem}  as follows.

\begin{corollary} \label{cor:no gaps of size bigger than two}
\[ P(G^{(2)}(A_n) = 0 \; i.o.) = 1. \]
\end{corollary}
{\bf Proof.} Just observe that by the above lemma for all $ a \in   {\cal G}$ with $N=|a|$  sufficiently big
\[P^{a, \emptyset}(\theta_0<N^{11/10}, N^{11/10} \leq \sigma_0 < N^{6/5})
\geq 1/2. \]
Thus, when starting with $a \in {\cal G}$ (with $|a|$ sufficiently
big) the probability that within the next $|a|^{6/5}$ periods there exists
a period $m$ such that all gaps of $A_m$ have size $1$ is bounded below by
1/2. Since by Corollary \ref{cor: recurrence of G} 
${\cal G}$ is recurrent for the process $A,$ 
the corollary now follows immediately from the Markov property of $A.$
\newline
\vspace{7 mm}
\fine

The next corollary will give us the recurrence of ${\cal S}_M$ for any fixed integer $M.$ 
\begin{corollary} \label{second corollary to Lemma1}
For any $M \in {\bf N}$ we have
\[ P(A_n \in {\cal S}_M   \; i.o.) = 1. \]
\end{corollary}      
{\bf Proof.} Clearly ${\cal S}_M \subset {\cal S}_{\tilde M}$ for $M \geq {\tilde M}.$  Hence we only need to show the above claim for integers $M > K.$  So fix $M > K.$ Note that ${\cal G}$ is recurrent for the process $A$ and $|A_n|
\rightarrow \infty$ as $n \rightarrow \infty.$ It thus suffices to show that for all $a \in {\cal G}$ (the size of which is sufficiently big)
\begin{equation} \label{lower bound for P(enter S within ...}
P^a(A_m \in {\cal S}_M \mbox{ for some } m \leq |a|^{6/5}) \geq 1/2. 
\end{equation}  
Now fix $a \in {\cal G}$ and let $N=|a|.$ Consider the Markov process $(A,I)$ with initial pair $(a, \emptyset),$ i.e., all gaps of $a$ are considered as being old. Note that $|{ C}^I|$ or $|{C}^{I,2}|$ can only increase when a particle gets deleted and the site
$X_{n,i}$ which gets unoccupied belongs to ${C}^I_{n,i}.$ Now, if
${ C}^I_{n,i-1}(\omega) = \emptyset$ and $X_{n,i}(\omega) \in {C}^I_{n,i}(\omega)$ then   
$|{ C}_{n,i}^I|(\omega) \leq 3$ and $|{C}_{n,i}^{2,I}| \leq
1.$ 
(For instance, if
$0$ and $5$ are the only new gap sites for $({ A}_{n,i-1},{ I}_{n,i-1})(\omega)$ and if $(n,i)$ is a deletion time and 
$X_{n,i}(\omega)=2$ then $ C_{n,i}^I=\{0,2,5\}$ and ${
C}_{n,i}^{2,I}=\{2\}.$) 
Hence, if ${\sigma}(\omega) > (n,i)$ then $|{ C}_{n,i}^I|(\omega) \leq 3$ and $|{ C}_{n,i}^{2,I}|(\omega) \leq
1.$ 

Thus, if $\omega \in \Omega$ satisfies $\theta_0(\omega) <
N^{11/10},$ $N^{11/10} \leq { \tau}_M(\omega) < N^{6/5}$ and ${
{\bar \sigma}}(\omega) > N^{6/5}$ then at the end of period $m=\tau_M(\omega)$  the following hold:   
(i) all gaps of $({A}_{m},{I}_{m})(\omega)$ are
new,  
(ii) $L^I({A}_{m},{I}_{m})(\omega) \leq h(N)/M,$ and 
(iii) $|{C}_{m}^I|(\omega) \leq 3$ and $|{C}_{m}^{2,I}|(\omega) \leq
1$.  
Note that (ii) and (iii) are in fact not only valid for new gap sites but, by using (i), for all gap sites.  Hence for $\omega$ as above we have ${A}_{\tau_M}(\omega) \in {\cal S}_M.$ 

We can thus conclude that 
$P^{a, \emptyset}({ A}_m \in {\cal S}_M \mbox{ for some } m \leq N^{6/5}) $ $\geq$   
$P^{a, \emptyset}( \theta_0 < N^{11/10}, N^{11/10} \leq {
\tau}_M < N^{6/5}, { {\bar \sigma}} > N^{6/5}).  $ 
Now by Lemma \ref{first lemma} (and recalling the definition of $\tau_M$) this last probability is bounded below by $1/2$ for all sufficiently large $N$ and observing that 
$P^{a, \emptyset}({ A}_m \in {\cal S}_M \mbox{ for some } m \leq N^{6/5}) = P^a({A}_m \in {\cal S}_M \mbox{ for some } m \leq N^{6/5})$,
 (\ref{lower bound for P(enter S within ...})
gets established. 
\newline
\vspace{7 mm}
\fine

The next lemma will be stated in two different versions. The first version will be used to finish  the proof of Theorem \ref{second theorem} and of part (1) of Theorem \ref{first theorem}. The second version is stated for technical reasons only. It will allow us to do jointly major parts of the proof of the two lemmas.

\begin{lemma} {\rm (First Version)} \label{second lemma}
Let $M \in {\bf N},$ $M > K$,  $a \in {\cal S}_M$ and let $N=|a|.$ Then for all sufficiently large $N$ we have 
\begin{equation} \label{inequality in the second lemma}
P^a(T^{M,2} \wedge \bar{S} > 2N, T_{M,2} < 2N) \geq 1-N^{-1/(2M)}. 
\end{equation} 
\end{lemma}

Note that if the Markov process $(A,I)$ starts with an initial pair $(a,\lambda),$ where $\lambda=\{1, \ldots, G(a)\},$ (i.e., all gaps of $a$ are considered as being new), then also at any later time the set of new gap sites coincides with the set of all gap sites. Thus defining \[ \tau^{M,2}=\inf \{n: L^I(A_n,I_n) \geq ((K-1)/K + 3/M )h(|A_0|) \}  \] 
and
\[ \tau_{M,2} = \inf \{ n \geq |A_0|: L^I(A_n,I_n) \leq h(|A_0|)/M
\}, \] 
Lemma \ref{second lemma} can be reformulated as follows

\vspace{7 mm}

\noindent {\bf Lemma 2.4.2 }(Second Version)  {\it
Let $M \in {\bf N}, M >K $,  $a \in {\cal S}_M$ and $\lambda=\{1, \ldots, G(a)\}.$ Let $N=|a|.$ Then for all sufficiently large $N$ we have }
\begin{equation}
 P^{a,\lambda}(\tau^{M,2} \wedge \bar{\sigma} > 2N, \tau_{M,2} < 2N) \geq 1-N^{-1/(2M)}. 
\end{equation}
The proof of the two lemmas will stretch over the next couple of sections. Note that in Lemma \ref{first lemma} we consider the Markov process $(A,I)$ with initial pair $(a, \lambda),$ where $a \in {\cal G}$ and $\lambda= \emptyset,$ (i.e., all gaps of $a$ are considered as being old) and in the second version of Lemma \ref{second lemma} we consider the Markov process $(A,I)$ with initial pair $(a, \lambda),$ where $a \in {\cal      
S}_M$ and $\lambda=\{1,\ldots,G(a)\},$ (i.e., all gaps of $a$ are considered as being new.) For all pairs $(a, \lambda)$ ($a$ a finite set of integers, $\lambda \subset \{1, \ldots ,G(a)\}$) we will now construct the Markov process $(A,I)$ with initial pair $(a,\lambda)$ on one and the same probability space $(\Omega,{\cal F},P)$ such that for every $\omega \in \Omega,$ $(A,I)(\omega)$ depends explicitly on $(a,\lambda)$ and such that for every 
${\tilde a} \subset {\bf Z},$ $|{\tilde a}| < \infty$ and $
{\tilde \lambda} \subset \{1, \ldots ,G({\tilde a})\}$ we have 
\[P^{a,\lambda}(A_{n,i}={\tilde a},I_{n,i}={\tilde \lambda})
 = P(A_{n,i}(a,\lambda)={\tilde a},I_{n,i}(a,\lambda)={\tilde \lambda}).\]
Note that the notation on the left does not refer to any specific construction of $(A,I)$ or of the underlying probability space, while the notation on the right refers to the explicit construction of $(A,I)$ as $(A,I)(a, \lambda)$ on $(\Omega,{\cal F},P).$ For ease of notation we will continue to write $(A,I)$ instead of $(A,I)(a,\lambda),$ but the reader is warned that not only $(A,I),$ but also the stopping times $\tau^1,$ $\tau_M$ etc. will depend pathwise on the initial pair $(a,\lambda).$   

\vspace{7 mm}

 For any pair $(\tilde a, \tilde{\lambda}),$ where $\tilde a$ is a finite set of integers and $\tilde{\lambda} \subset \{1, \ldots, G(\tilde a)\}$ denotes the indices of new gaps of $\tilde a,$ let us partition the inner boundary of ${\tilde a}$ into  
$\partial^o (\tilde a,\tilde{\lambda})$ and  
$\partial^n (\tilde a,\tilde{\lambda}),$ 
where $\partial^o (\tilde a,\tilde{\lambda})$ (resp. $\partial^n (\tilde a,\tilde{\lambda})$) consists of the boundary points of $\tilde a$  which belong to old (resp. new) gaps. Let $V^o(\tilde a,\tilde{\lambda})$ denote the subset of  $\tilde a$ consisting of all points which are adjacent to an old gap but not an endpoint of $\tilde a$ and let the set $V^n(\tilde a,\tilde{\lambda})$ consist of all points of $\tilde a$ which are within distance 4 to a new gap site but which are neither adjacent to an old gap nor an endpoint of $\tilde a.$

\section{An Explicit Construction of the Markov Process $(A,I)$}

Let $(\Omega, {\cal F}, P)$ be a probability space and let $\{U_j\}$ be a sequence  of independent random variables on $(\Omega, {\cal F}, P),$ which are uniformly
distributed on $(0,1).$ Let $(a,\lambda)$ with ${a} \subset {\bf Z},$ $|{a}| < \infty$ and $
{\lambda} \subset \{1, \ldots ,G({a})\}$ be our initial pair. Define $A_0=A_0(a,\lambda)=a$ and  $I_0=I_0(a,\lambda)=\lambda.$ Now suppose that $(A,I)_{n,i}=(A,I)_{n,i}(a,\lambda)$ has already been defined. In order to define $(A,I)_{n,i+1}=(A,I)_{n,i+1}(a,\lambda)$ it suffices to specify how 
the point $X_{n,i+1}$ which gets added to
(resp. deleted from) the set $A_{n,i}$ gets defined.

If $(n,i+1)$ is an addition time let us enumerate the boundary sites 
of $A_{n,i}$ (depending on the properties of those sites) as follows: 
Let $x_1 <  \ldots <
x_{|\partial^o (A_{n,i},I_{n,i})|}$ be the boundary points of
$A_{n,i}$ which belong to old gaps and let $x_{|\partial^o (A_{n,i},I_{n,i})|+1}=\min A_{n,i}-1$ and $x_{|\partial^o (A_{n,i},I_{n,i})|+2}=\max(A_{n,i})+1.$
Let
$x_{|\partial^o (A_{n,i},I_{n,i})|+3} <  \ldots <
x_{|\partial A_{n,i}|-|C^I|(A_{n,i},I_{n,i})}$ denote the
boundary points of $A_{n,i}$ which belong to new gaps and which do not
belong to $C^{I}_{n,i}$ and let  
$x_{|\partial A_{n,i}|-|C^I|(A_{n,i},I_{n,i})+1}< \ldots
< x_{| \partial A_{n,i}|}$  denote the boundary points of $A_{n,i}$ 
which belong to $C^{I}(A_{n,i},I_{n,i}).$
Define 
\begin{equation} \label{choice of the new point to be added}
 X_{n,i+1}(\omega)=x_j(\omega)  \mbox{ if } (j-1)/{| \partial A_{n,i}|(\omega)} <
U_{(n-1)\kappa+i+1}(\omega) \leq j/{| \partial A_{n,i}|(\omega)} .
\end{equation} 
Similarly, if $(n,i+1)$ is a deletion time let us enumerate the points of $A_{n,i}$ as follows:
Let the points $y_1 <  \ldots
<y_{|V^n|(A_{n,i},I_{n,i})}$ of $V^n(A_{n,i},I_{n,i})$ be followed by those points of $A_{n,i}$ which are neither endpoints nor adjacent to an old gap nor within distance 4 to a new
gap site (once again enumerated in increasing order) and then by the endpoints $y_{| A_{n,i}|-|V^{o}|(A_{n,i},I_{n,i})-1}= \min A_{n,i}$ and $y_{| A_{n,i}|-|V^{o}|(A_{n,i},I_{n,i})} = \max A_{n,i}.$ 
Finally let $y_{| A_{n,i}|-|V^{o}|(A_{n,i},I_{n,i})+1} <
\ldots < y_{| A_{n,i}|}$ enumerate the points of $V^{o}(A_{n,i},I_{n,i}).$ 
Define
\begin{equation} \label{choice of the point to be deleted}
X_{n,i+1}(\omega)=y_j(\omega)  \mbox{ if } (j-1)/{|A_{n,i}|(\omega)} < U_{(n-1)\kappa+i+1}(\omega) \leq j/{| A_{n,i}|(\omega)} . 
\end{equation}
Before proceeding further let us introduce one more piece of notation. For $B,C \in {\cal F}$ define $B \sqsubset C$ (res. $B \sqsupset C$) if there exist ${\tilde B}, {\tilde C} \in {\cal F}$ with $P(B \bigtriangleup {\tilde B}) = 0,$ $P(C \bigtriangleup {\tilde C}) = 0$ and  ${\tilde B} \subset {\tilde C}$ (resp. $B \supset C$).

\section{Excursions of $L$ - Construction of $\hat L^I$}

Let $M$ be a fixed positive integer.
Note that when we add a particle, the number of new gap sites can never
increase and it decreases if and only if the newly added site was a new
gap site for $(A_{n,i-1},I_{n,i-1}).$ When we delete a particle, the
number of new gap sites can increase at most by 1. Thus  the increments of $L^I=L^I(A,I)$ can be estimated as follows, 
\begin{eqnarray}
L^I_{n}-L^I_{n-1} &\leq& - \sum_{i=1}^{K} {\bf 1}_{ \{X_{n+1,i} \in
\partial^n (A_{n,i-1},I_{n,i-1}) \} } + (K-1) \nonumber \\
                  &=&  -\sum_{i=1}^{K} (1-{\bf 1}_{ \{X_{n,i} \in \bar{\partial}
A_{n,i-1} \cup \partial^o(A_{n,i-1},I_{n,i-1}) \} })+ (K-1) \nonumber \\
                  &=& -1 + \sum_{i=1}^{K} {\bf 1}_{ \{X_{n,i} \in \bar{\partial}
A_{n,i-1} \cup \partial^o(A_{n,i-1},I_{n,i-1}) \} } \nonumber \\
                  &=&  -1 + \sum_{i=1}^{K} {\bf 1}_{ \{U_{(n-1)\kappa+i} \leq (2 + |\partial^o(A_{n,i-1},I_{n,i-1})|)
/|\partial A_{n,i-1}| \} }. \label{general ineq for increments of LI}
\end{eqnarray}
Now, if we have a starting pair $(a,\emptyset)$ with $a \in \cal {G}$ then $|\partial^o(A_{n,i-1},I_{n,i-1})| \leq 20K$ for all time points $(n,i)$ since $G(a) \leq 10K$ and the number of old gaps never increases. Moreover, if $\sigma(\omega) > (n,i-1)$ then $(A_{n,i-1},I_{n,i-1})(\omega)$ has no new gaps of size bigger than 2 and at most one new gap of size 2 and so $(L^I_{n,i-1} - G^I_{n,i-1})(\omega) \leq 1.$
Similarly for a starting pair $(a,\lambda)$ with $a \in {\cal S}_M$ and $\lambda=\{1, \ldots, G(a)\}$ we have  $|\partial^o(A_{n,i-1},I_{n,i-1})| = 0 $  and, if $\sigma(\omega) > (n,i-1)$ then $(L^I_{n,i-1} - G^I_{n,i-1})(\omega) \leq 1.$ (Here we use that $L(a)-G(a) \leq 1$ since $a \in {\cal S}_M$). Finally note that if 
$L^I(A_{n,i-1},I_{n,i-1})(\omega) \geq h(N)/M$ and $(L^I_{n,i-1} - G^I_{n,i-1})(\omega) \leq 1$ then $|\partial A_{n,i-1}| \geq 
h(N)/M.$  
Thus for any starting pair $(a,\lambda)$ as in Lemma \ref{first lemma}, (i.e., $a \in {\cal G}$ and $\lambda=\emptyset$)  or in the second version of Lemma \ref{second lemma} ,
(i.e., $a \in {\cal S}_M$ and $\lambda=\{1,\ldots,G(a)\}$), we get the following: 

If $L^I(A_{n-1},I_{n-1})(\omega) \geq h(N)/M+K$ and 
$\bar{\sigma}(\omega) > n-1$ then 
\begin{equation} \label{ineq: essential for the construction of L-hat}
(L^I_{n}-L^I_{n-1})(\omega) \leq  -1 + \sum_{i=1}^{K} {\bf 1}_{ \{U_{(n-1)\kappa+i} \leq  c/h(N) \} }(\omega),
\end{equation}
where $c=22KM.$
(Here we also use that if $L^I(A_{n-1},I_{n-1})(\omega) \geq h(N)/M+K$ then $L^I(A_{n,i-1},I_{n,i-1})(\omega) \geq h(N)/M$ for $i=1 ,\ldots, K$ and if $\bar{\sigma}(\omega) > n-1$ then $\sigma(\omega) > (n,K)$ since obviously no new empty site can be created when we add a particle). We will use the above estimate to compare in a path-by-path way `excursions' of $L^I$ (above $h(N)/M+K$), which start before $\sigma$, to the corresponding
parts of a random walk. So let us proceed to define these `excursions'.

\vspace{7 mm}

If $L^I_{n}(\omega) < K+h(N)/M$ for all $n < \bar{\sigma}(\omega),$ let ${\cal N}(\omega)=0,$ i.e., no excursion of $L^I(\omega)$ above $K+h(N)/M$ starts before $\sigma(\omega).$ Otherwise let us recursively define those excursions as follows: 
Let $f_0(\omega) = \inf \{n < \bar{\sigma}(\omega):
L^I_{n} \geq K+h(N)/M \}$   
 denote the start of the first `excursion' of $L^I(\omega)$ above $K+h(N)/M$. 
Suppose that $f_{i-1}(\omega),$ the start of the $i$-th excursion of $L^I,$
is already defined with $f_{i-1}(\omega) < \bar{\sigma}(\omega).$ Let us define the length $\rho_i(\omega)$ of the $i$-th
`excursion' as
\[ \rho_i(\omega)= \left\{ \begin{array}{ll}
                             \bar{\sigma}(\omega)-f_{i-1}(\omega) & \mbox{if $L^I_{n}(\omega) \geq K+h(N)/M$} \\
                        & \mbox{  for $n=f_{i-1}(\omega), \ldots, \bar{\sigma}(\omega)-1$} \\
                             \inf\{n:
L^I_{f_{i-1}(\omega)+n}(\omega) < K+h(N)/M \} & \mbox{otherwise}
                           \end{array}
                   \right. ,  \]
and let
\[e^i(\omega)=(L^I_{f_{i-1}+j})_{j=0 ,\ldots, \rho_i}(\omega)\]
denote the $i$-th excursion of $L^I(\omega).$ (The reader may note that we consider the present excursion as being terminated when we reach period $\bar{\sigma}(\omega)$ even if $L^I(\omega)$ is still bounded below by $K+h(N)/M$). Let us define the start of the next excursion of $L^I$ (if there is any) as
$f_i(\omega) = \inf \{n:
f_{i-1}(\omega)+\rho_i(\omega)<n<\bar{\sigma}(\omega) \mbox{
and } L^I_{n}(\omega) \geq K+h(N)/M \} .$
 
We have thus recursively defined all the excursions of $L^I(\omega)$ above $K+h(N)/M$, which start before period $\bar{\sigma}(\omega).$ Define ${\cal N}(\omega)$ to be the number of such excursions.

Now let us assume that on our probability space $(\Omega, {\cal F},P)$ a process $U^{(1)}$ is defined, which is independent of $U$ and such that the random variables $U^{(1)}_j$ are independent and uniformly distributed on $(0,1).$
We will recursively construct a process $U^L$ and we define $\hat{L}^I$ as
\begin{eqnarray*}
\hat{L}^I_0   &=&  0, \\
\hat{L}^I_n   &=&  -n + \sum_{l=1}^n \sum_{j=1}^{K} {\bf 1}_{ \{U^L_{(l-1)K+j}
\leq c/h(N) \} }, \\
\end{eqnarray*}
where $c=22KM.$
We are interested in successive $K$-decreases of ${\hat L}^I,$ a $K$-decrease being defined as follows. Let $\hat{f}_0  = 0$ and for each $n \in {\bf N}$ let $\hat{f}_n$ be recursively defined as $\hat{f}_n(\omega)  = \inf \{j > \hat{f}_{n-1}(\omega): (\hat{L}^I_j
-  \hat{L}^I_{\hat{f}_{n-1}})(\omega) \leq -K \}.$ Then  $({\hat L}^I_j)_{{\hat f}_{n-1} \leq j \leq {\hat f}_n}$ denotes the $n$-th $K$-decrease of ${\hat L}^I.$  Let $\hat{\rho}_n=\hat{f}_n-\hat{f}_{n-1}$ denote its length.

\vspace{7 mm}

For each $\omega \in \Omega$ we construct the sequence $U^L(\omega)$ by pasting together parts of $U(\omega)$ and of $U^{(1)}(\omega),$ the decision, whether to choose the next $K$ elements for $U^L(\omega)$ from $U(\omega)$ or from  $U^{(1)}(\omega),$ depending in each step only on those elements of $U(\omega)$ and of $U^{(1)}(\omega)$ which have been investigated so far.
   
If no excursion of $L^I(\omega)$ above $K + h(N)/M$ starts before period 
$\bar{\sigma}(\omega)$ 
 we put $U^L(\omega)=U^{(1)}(\omega).$ Otherwise, starting at the end of period $f_0(\omega),$ we select successive elements of the sequence $U(\omega)$ (but only those which correspond to addition times) until we reach the end of the first excursion of $L^I(\omega).$ If by then $\hat{L}^I(\omega)$ has not yet terminated its first K-decrease, we fill in with elements from the sequence $U^{(1)}(\omega)$ until $\hat{L}^I(\omega)$ terminates its first K-decrease. (We will explain below, why $\hat{L}^I(\omega)$ can not terminate its first $K$-decrease before all elements of the sequence $U(\omega),$ which correspond to addition times and which contribute to the first excursion of $L^I(\omega),$ are selected for $U^L(\omega)$).
We keep iterating this procedure (replacing the previous $f_{i-1}(\omega)$ by $f_{i}(\omega)$ in each new iteration) as long as this is possible. Thus, if $\hat{f}_{i-1}(\omega) < \infty$ then all elements of the sequence $U(\omega),$ which correspond to addition times and which contribute to the $i$-th excursion of $L^I(\omega)$ ($i < {\cal N}(\omega)+1$), are selected for $U^L(\omega).$ If by then $\hat{L}^I(\omega)$ has not yet terminated its $i$-th $K$-decrease then we switch to the sequence $U^{(1)}(\omega)$ and we keep selecting from this sequence until this $i$-th $K$-decrease of $\hat{L}^I(\omega)$ is completed.

If ${\cal N}(\omega) = \infty,$ the whole sequence $U^L(\omega)$ is thus defined and if  ${\cal N}(\omega) < \infty,$  $U^L(\omega)$ has been constructed in this way up to (including) ${\cal N}(\omega)$ $K$-decreases for ${\hat L}^I(\omega).$
So suppose that ${\cal N}(\omega)$ $K$-decreases for  ${\hat L}^I(\omega)$ have been completed in this way.  
After this time all further elements for $U^L(\omega)$ are selected from the sequence $U^{(1)}(\omega).$  

Thus, $U^L$ satisfies the following: 
For $1 \leq i < {\cal N}(\omega)+1$ we have that if $\hat{f}_{i-1}(\omega) < \infty$ then 
\[U^L_{(\hat{f}_{i-1}+l-1)K+j}(\omega) = U_{({f}_{i-1}+l-1)\kappa+j}(\omega) \qquad \qquad \mbox{ for $1 \leq l \leq \rho_i,$ $1 \leq j \leq K,$}  \]
and
\[U^L_{(\hat{f}_{i-1}+\rho_i)K+m}(\omega) = U^{(1)}_{\sum_{k=1}^{i-1}(\hat{\rho}_k-\rho_k)K + m}(\omega) \qquad \qquad \mbox{for $1 \leq m \leq (\hat{\rho}_i-\rho_i)K$}. \]
If ${\cal N}(\omega) < \infty$ and $\hat{f}_{\cal N}(\omega) < \infty,$
then
\[ U^L_{\hat{f}_{\cal N}K+m}(\omega) = U^{(1)}_{\sum_{k=1}^{\cal N}(\hat{\rho}_k-\rho_k)K + m}(\omega) \qquad \qquad \mbox{for all $m \in {\bf N}$}.
\]

\vspace{7 mm}

Let us now verify that for any $\omega$ satisfying $\hat{f}_{i-1}(\omega) < \infty,$  $\hat{L}^I(\omega)$ can not terminate its $i$-th $K$-decrease ($i < {\cal N}(\omega)+1$) before all elements of the sequence $U(\omega),$ which correspond to addition times and which contribute to the $i$-th excursion of $L^I(\omega),$ are selected for $U^L(\omega)$.
 
If ${\cal N}(\omega) > 0,$ then the excursion $e^1(\omega)$ is defined, and satisfies $L^I_{f_0}(\omega) \leq (K-1)+h(N)/M+(K-1)=2K-2+h(N)/M.$  
By definition of the first excursion we have  $L^I_{f_{0}+j}(\omega) \geq K+h(N)/M$  and $f_0+j < \bar{\sigma}(\omega)$ for any $0 \leq j < \rho_1(\omega).$ Thus by (\ref {ineq: essential for the construction of L-hat})
for any starting pair $(a,\lambda)$ as previously described,   
the increment of $L^I(\omega)$ satisfies
$(L^I_{f_0+j+1}-L^I_{f_0+j})(\omega) \leq  -1 + \sum_{i=1}^{K} {\bf 1}_{ \{U_{(f_0+j)\kappa+i} \leq  c/h(N) \} }(\omega).$ Now using that the elements of $U(\omega),$ which correspond to addition times and which contribute to the first excursion for $L^I(\omega),$ are the first ones to be selected  for $U^L(\omega)$ we can conclude that $(L^I_{f_0+j+1}-L^I_{f_0+j})(\omega) \leq  -1 + \sum_{i=1}^{K} {\bf 1}_{ \{U^L_{jK+i} \leq  c/h(N)\} }(\omega) $ and hence
\[ (L^I_{f_0+j+1}-L^I_{f_0+j})(\omega) \leq ({\hat L}^I_{j+1}-{\hat L}^I_{j})(\omega) \]
for any $0 \leq j < \rho_1(\omega).$
Thus
\begin{eqnarray*} 
({\hat L}^I_{j}-{\hat L}^I_{0})(\omega) &\geq& (L^I_{f_{0}+j}-L^I_{f_{0}})(\omega) \\
&\geq& K+h(N)/M-(2K-2+h(N)/M) \\
&=& -(K-2)
\end{eqnarray*}
for all $j<\rho_1(\omega)$ and so the first $K$-decrease for ${\hat L}^I(\omega)$ can not be completed before all elements of $U(\omega),$ which correspond to addition times and which contribute to the first excursion for $L^I(\omega),$ have been selected.
In particular we can conclude that
$\rho_1(\omega) \leq  \hat{\rho}_1(\omega)$
and
for $j=1 ,\ldots, \rho_1(\omega),$
$
 (e^1_j-e^1_0)(\omega) \leq (\hat{L}^I_{j}-\hat{L}^I_{0})(\omega). $
By our selection rule we can repeat this argument for all other excursions and so our claim follows. Moreover, we can conclude that for all $i, \omega$ with $i < {\cal N}(\omega)+1$ and $\hat{f}_{i-1}(\omega) < \infty$ we have 
\begin{equation} \label{comparison for the length of excursions}
\rho_i(\omega) \leq  \hat{\rho}_i(\omega)
\end{equation}
and 
\begin{equation} \label{comparison for the height of excursions}
 (e^i_j-e^i_0)(\omega) \leq (\hat{L}^I_{\hat{f}_{i-1}+j}-\hat{L}^I_{\hat{f}_{i-1}})(\omega) 
\end{equation}
for $j=1, \ldots, \rho_i(\omega).$ 

The following property of $U^L$ will be most useful in the sequel.

\begin{proposition}
$U^L$ is a process of independent random variables which are uniformly distributed
on $(0,1).$ 
\end{proposition}

{\bf Proof.} Let us rewrite the procedure in order to show that the construction fits into the scheme described in Proposition \ref{prop:good selection rules}. 
Define $S_1(\omega)   =   f_0(\omega) \wedge \bar{\sigma}(\omega)$ and $T_1(\omega)   =   0.$ 
Let $i \geq 1$ and let us suppose that the first $i-1$ periods of the
sequence $U^L(\omega)$ were already constructed and that we defined  $S_i(\omega)$ and $T_i(\omega),$ $S_i(\omega)$ representing the number of periods of length $\kappa$ of the sequence $U(\omega)$ that we investigated so far and $T_i(\omega)$ representing the number of periods of length $K$ of the sequence $U^{(1)}(\omega)$ that we investigated so far. In order to define the $i$-th period of length $K$ of $U^L(\omega)$ and  $S_{i+1}(\omega)$ and $T_{i+1}(\omega),$ let us consider the following three cases:

(i) If $S_i(\omega)< \bar{\sigma}(\omega)$ and $L^I_{S_i}(\omega) \geq K+h(N)/M,$ let $U^L_{(i-1)K+j}(\omega)= U_{S_i \kappa +j}(\omega)$ for $j=1 ,\ldots, K, \; T_{i+1}(\omega)=T_i(\omega)$ and define $S_{i+1}(\omega)$ as follows: \newline
(a) If $L^I_{S_i+1}(\omega) \geq K+h(N)/M,$ or if (b) $L^I_{S_i+1}(\omega) < K+h(N)/M$ and $\hat{L}^I(\omega)$ does not finish a $K$-decrease in step $i,$ let
$S_{i+1}(\omega)=S_{i}(\omega)+1.$
Otherwise, (i.e., if (c) $L^I_{S_i+1}(\omega) < K+h(N)/M$ and $\hat{L}^I(\omega)$ finishes its next $K$-decrease in step $i$), let
$ S_{i+1}(\omega)=\inf \{n>S_{i}(\omega):L^I_{n}(\omega) \geq
K+h(N)/M \} \wedge \bar{\sigma}(\omega) .$ 

(ii) If $S_i(\omega)< \bar{\sigma}(\omega)$ and $L^I_{S_i}(\omega) < K+h(N)/M,$ let $U^L_{(i-1)K+j}(\omega)= U^{(1)}_{T_i K +j}(\omega)$ for $j=1 ,\ldots, K, \; T_{i+1}(\omega)=T_i(\omega)+1$ and define $S_{i+1}(\omega)$ as follows: \newline
If $\hat{L}^I(\omega)$ does not finish a $K$-decrease in step $i,$ let
$S_{i+1}(\omega)=S_{i}(\omega). $
Otherwise let 
$ S_{i+1}(\omega)=\inf \{n>S_{i}(\omega):L^I_{n}(\omega) \geq
K+h(N)/M \} \wedge \bar{\sigma}(\omega) .$
 
(iii) If $S_i(\omega) = \bar{\sigma}(\omega),$ let $U^L_{(i-1)K+j}(\omega)= U^{(1)}_{T_i K +j}(\omega)$ for $j=1 ,\ldots, K, \; T_{i+1}(\omega)=T_i(\omega)+1$ and 
$S_{i+1}(\omega)=S_{i}(\omega). $

\vspace{7 mm}

Letting ${\cal F}_{ij}=\sigma(U_1 ,\ldots, U_i,U_1^1 ,\ldots, U_j^1),$ it is not hard to check that for each $n \in {\bf N},$ $(S_n \kappa, T_n K)$ is a stopping time with respect to the filtration  $\{{\cal F}_{(i,j)}:(i,j) \in {\bf N}^2_{\infty}\}.$ Moreover, there exist sets $B_n^1,B_n^2 \in  {\cal F}_{S_n \kappa, T_n K}$ with $U^L_{(n-1)K+j}={\bf 1}_{B_n^1} U_{S_n \kappa + j} + {\bf 1}_{B_n^2} U^{(1)}_{T_n K + j}$ for $j=1 ,\ldots, K$ and such that $(S_{n+1},T_{n+1}) \geq  {\bf 1}_{B_n^1} (S_{n}+1,T_{n+1}) +
{\bf 1}_{B_n^2} (S_{n},T_{n}+1).$ Thus by Proposition \ref{prop:good selection rules} our claim follows.
\newline
\vspace{7 mm}
\fine

\begin{remark}
The previous proposition and our definition of ${\hat L}^I$ imply in particular that for all $i \in {\bf N}$, ${\hat f}_i < \infty$ a.s.
\end{remark}  

\section{Set Relations and Probability Estimates
     Related to the Process $\hat L^I$}

Let $M>K$ be a fixed integer and let us fix an initial pair $(a, \lambda)$ as in Lemma \ref{first lemma}, (i.e., $a \in \cal{G}$ and $\lambda=\emptyset$), or in the second version of Lemma \ref{second lemma}, (i.e., $a \in {\cal S}_M$ and $\lambda=\{1, \ldots, G(a)\}$).  Let $N=|a|.$ Let $(A,I)=(A,I)(a,\lambda)$ be the Markov process with initial pair $(a, \lambda)$ constructed from the pocess $U$ and let ${\hat L}^I$ be defined as in the previous section. (We continue to use the definitions and the set-up of the previous section). 
Let us define the events 
\[ {\cal E}^{1} = \bigcup_{i=1}^{\lfloor{N^{6/5}}\rfloor}\{ \max_{\displaystyle{j \! \leq \hat{\rho}_i}}({\hat L}^I_{\hat{f}_{i-1}+j} - {\hat L}^I_{\hat{f}_{i-1}})  
 \geq 2 (1-1/K)h(N) \} \]
and 
\[ {\cal E}^{2} = \bigcup_{i=1}^{2N} \{ \max_{\displaystyle{j \! \leq \hat{\rho}_i}}({\hat L}^I_{\hat{f}_{i-1}+j} -
{\hat L}^I_{\hat{f}_{i-1}})  
 \geq ((K-1)/K + 1/M)h(N) \}. \]
Then the following holds.

\begin{proposition} 
For all initial pairs $(a,\lambda)$ with $a \in \cal{G},$ $\lambda=\emptyset$ and $N=|a|$ sufficiently large, we have
\begin{equation} \label{set inclusion for E1}
   \{ \bar{\sigma} \wedge \tau^1 \leq N^{6/5},
 \tau^1 \leq \bar{\sigma} \} \sqsubset {\cal E}^{1},
\end{equation}
and for all initial pairs $(a,\lambda)$ with $a \in {\cal S}_M,$   $\lambda=\{1, \ldots, G(a) \}$ and $N=|a|$ sufficiently large, 
  we have
\begin{equation} \label{set inclusion for E2}
  \{ \bar{\sigma} \wedge \tau^{2,M} \leq 2N,
\tau^{2,M} \leq \bar{\sigma} \}  \sqsubset {\cal E}^{2}.
\end{equation}
\end{proposition}

(Note that the sets ${\cal E}^{1}$ and ${\cal E}^{2}$ and the stopping times in (\ref{set inclusion for E1}) and (\ref{set inclusion for E2}) depend on the initial pair $(a,\lambda)$).

\vspace{7 mm}

{\bf Proof.} Let us first assume that $a \in \cal{G}$ and  $\lambda=\emptyset$. Clearly, any $\omega \in \{ \bar{\sigma} \wedge \tau^1 \leq N^{6/5},
 \tau^1 \leq \bar{\sigma} \}$ satisfies $ \tau^1(\omega) \leq \bar{\sigma}(\omega)$  and  $\tau^1(\omega) \leq
N^{6/5}.$ Hence
there exists an excursion of $L^I(\omega)$ above $K + h(N)/M$   
which during the first $\lfloor{N^{6/5}}\rfloor$ periods gets higher than $(1+(K-1)/K)h(N)$ and thus (assuming that $N$ is sufficiently large)
 there exists an index $i$, $1 \leq
i \leq \lfloor{N^{6/5}}\rfloor \wedge {\cal N}(\omega),$ with
\begin{eqnarray*}
\max_{\displaystyle{j \!  \leq \rho_i(\omega)}}(e^i_j -e^i_0)(\omega) & \geq & (1+(K-1)/K)h(N) - (h(N)/M+K+(K-2)) \nonumber\\
&  \geq & ( 2 - 1/K - 1/M - (2K-2)/h(N)) h(N)   
\nonumber \\
&  \geq & 2(1-1/K) h(N). 
\end{eqnarray*}
Now ${\hat f}_{i-1} < \infty$ a.s. and so by  (\ref{comparison for the length of excursions}) and (\ref{comparison for the height of excursions}) for almost every $\omega$ as above we have
\begin{eqnarray*}
\max_{\displaystyle{j \! \leq \rho_i(\omega)}} (e^i_j -e^i_0)(\omega) 
& \leq & \max_{\displaystyle{j \! \leq \rho_i(\omega)}}({\hat L}^I_{\hat{f}_{i-1}+j} -
{\hat L}^I_{\hat{f}_{i-1}})(\omega) \nonumber \\
& \leq & \max_{\displaystyle{j \! \leq \hat{\rho}_i(\omega)}}({\hat L}^I_{\hat{f}_{i-1}+j} -
{\hat L}^I_{\hat{f}_{i-1}})(\omega). 
\end{eqnarray*}
Hence there exists $i \leq \lfloor{N^{6/5}}\rfloor$ with 
$ \max_{\displaystyle{j \! \leq \hat{\rho}_i(\omega)}}({\hat L}^I_{\hat{f}_{i-1}+j} -
{\hat L}^I_{\hat{f}_{i-1}})(\omega)  
 \geq  2(1-1/K)h(N)$ and so (\ref{set inclusion for E1}) follows.

Similarly, for an initial pair $(a,\lambda)$ with $a \in {\cal S}_M$ and  $\lambda=\{1, \ldots, G(a) \}$ we have that if $\tau^{2,M}(\omega) \leq \bar{\sigma}(\omega)$ and $\tau^{2,M}(\omega) \leq 2N,$ then there exists an excursion of $L^I$ which during the first $2N$ periods gets higher than $((K-1)/K) + 3/M)h(N).$ Hence
(assuming that $N$ is sufficiently large)
there exists an index $i$, $1 \leq
i \leq 2N,$ with
\begin{eqnarray*}
\max_{\displaystyle{j \! \leq \rho_i(\omega)}}(e^i_j -e^i_0)(\omega) & \geq & ((K-1)/K) + 3/M)h(N) - (h(N)/M+K+(K-2)) \nonumber\\
&  \geq & ((K-1)/K + 1/M) h(N), 
\end{eqnarray*}
and so 
$ \max_{\displaystyle{j \! \leq \hat{\rho}_i(\omega)}}({\hat L}^I_{\hat{f}_{i-1}+j} -
{\hat L}^I_{\hat{f}_{i-1}})(\omega)  
 \geq  ((K-1)/K + 1/M) h(N) $ for almost every $\omega$ satisfying the conditions above. Thus (\ref{set inclusion for E2}) gets established for this choice of $(a, \lambda)$. 
\newline
\vspace{7 mm}
\fine

The next proposition will give us useful upper bounds
for the probability of the events ${\cal E}^{1}$ and ${\cal E}^{2}$.

\begin{proposition}
Let $M > K.$ There exists $N_0$ such that for all initial pairs $(a,\emptyset)$ with $a \in \cal{G}$ and $N=|a| \geq N_0$ we have 
\begin{equation} \label{estimate for calE1}
P({\cal E}^1) \leq N^{-3/5} 
\end{equation}
and such that for all initial pairs $(a,\lambda)$ with $a \in {\cal S}_M$ , $\lambda=\{1, \ldots, G(a) \}$  and $N=|a| \geq N_0$ we have
\begin{equation} \label{estimate for calE2}
P({\cal E}^{2}) \leq N^{-5/(8M)}. 
\end{equation}
\end{proposition}
{\bf Proof.} Fix an initial pair $(a,\emptyset)$ with $a \in \cal{G}$ and let $N=|a|.$ For any $i \in {\bf N}$ define the event $E_i^{1}= \{ \max_
{j  \leq \hat{\rho}_i}(\hat{L}_{\hat{f}_{i-1}+j} - \hat{L}_{\hat{f}_{i-1}})  
 \geq 2 (1-1/K)h(N) \}.$ In order to show (\ref{estimate for calE1}) it will suffice to show that 
\[ P(E_i^1)=P(E_1^1) \; \mbox{  for all $i$} \qquad \mbox{and} \qquad P(E_1^1) \leq
N^{-9/5}.\]
Let $Z^i$ denote the process $Z_n^i=\hat{L}^I_{\hat{f}_{i-1}+n} - 
\hat{L}^I_{\hat{f}_{i-1}}$, thus  
\[Z_n^i=-n+\sum_{l=1}^{n}\sum_{j=1}^{K} {\bf
1}_{\{U_{(\hat{f}_{i-1}+l-1)K+j}^L \leq c/h(N) \}},\]
where $c=22KM.$
Let $\hat{e}^i$ be obtained by stopping $Z^i$
at its first $K$-drop, i.e., $\hat{e}_n^i = Z_{n \wedge \hat{\rho}}^i,$
where $\hat{\rho}(\omega)=\inf\{n:Z_n^i(\omega) \leq -K \}.$
 Since the random variables
$U_n^L$ are independent and identically distributed and $\hat{f}_{i-1}K$ is
a stopping time for the process $U^L,$ we have $Z^i \sim Z^1$ for
all $i$ and hence also $\hat{e}^i \sim \hat{e}^1.$ Noting that $E_i^1= \{
\sup_{n} \hat{e}_n^i \geq 2(1-1/K)h(N) \}$ it follows that 
$P(E_i^1) = P(E_1^1) .$

Next let us define a process $\tilde{Z}$ by 
\[ \tilde{Z}_{nK+k}= \sum_{l=1}^{n} \sum_{j=1}^{K} {\bf
1}_{\{U_{(l-1)\kappa+j}^L \leq c/h(N) \}} + \sum_{j=1}^{k} {\bf
1}_{\{U_{n \kappa+j}^L \leq c/h(N) \}} ,\]
$n \in \Nzero,$ $0 \leq k \leq K,$
and let
\[ \tilde{\tau}=\inf \{m \geq 0 : \tilde{Z}_{mK}-m \geq 2(1-1/K)h(N)
\mbox{ or } \tilde{Z}_{mK}-m \leq -K \} .\]
Then $P(E_1) = P(\tilde{Z}_{\tilde{\tau}K}-{\tilde \tau} \geq 2(1-1/K)h(N)).$ Note that
$\tilde{Z}$ is a random walk with step law $\mu$ defined by $\mu(1)=c/h(N)$
and $\mu(0)=1-c/h(N).$ Thus by Proposition \ref{prop: random walk estimates} (with $\eta=(K-1)/K$ and
$\varepsilon=1/5$) we can conclude that for all sufficiently large $N$ \[ P(E_1^1) \leq N^{-9/5}.\]
Let us now consider an initial pair $(a,\lambda)$ with $a \in {\cal S}_M$ ,$\lambda=\{1, \ldots, G(a) \}$  and  $N=|a|.$ For any $i \in {\bf N}$ define the event 
$E_i^{2}= \{ \max_{j \! \leq \hat{\rho}_i}(\hat{L}_{\hat{f}_{i-1}+j} -
\hat{L}_{\hat{f}_{i-1}})  
 \geq ((K-1)/K + 1/M)h(N) \}.$ Proceeding as above we can show that $P(E_i^2)=P(E_1^2)$ for all $i$ and using Proposition \ref{prop: random walk estimates} (with $\eta=1/M$ and
$\varepsilon=1/(4M)$) we obtain $P(E_1^2) \leq N^{-(1+3/(4M))}.$ Thus (\ref{estimate for calE2}) follows.
\newline
\vspace{7 mm}
\fine

\section{Investigation of $\{\bar{\sigma} \wedge \tau^1 > N^{6/5}, \sigma_0 < N^{6/5} \}$ (resp. $\{\bar{\sigma} \wedge \tau^{2,M} > 2N \}$) } 

Let $M>K$.   Let $(A,I)=(A,I)(a,\lambda)$ be the Markov process with initial pair $(a, \lambda)$ constructed from the pocess $U.$ (We continue to use the set-up and the definitions of the previous section. Thus, in particular, $(\Omega, {\cal F},P)$ is a probability space on which two independent processes $U$ and $U^{(1)}$ are defined, the random variables $U_n, n \in {\bf N},$ (resp. $U^{(1)}_n, n \in {\bf N}$), being independent and uniformly distributed on $(0,1)$   ).
 
We will investigate the event $ \{ \bar{\sigma} < \tau^1;  \bar{\sigma} \leq N^{6/5} \},$  if the initial pair $(a, \lambda)$ satisfies  $a \in \cal{G}$ and $\lambda=\emptyset$, and we will investigate the event $ \{ \bar{\sigma} < \tau^{2,M};  \bar{\sigma} \leq 2N \}$, if $a \in \cal{S}_M$ and $\lambda=\{1, \ldots, G(a)\}.$

So let $(a, \lambda)$ satisfy any of the above sets of conditions . 
Note that if  $\omega \in \Omega$ satisfies  
$\bar{\sigma}(\omega,a,\lambda) \leq N^{6/5}$ then there exists a period $n \leq N^{6/5}$ and a deletion time $(n,i)$ such that at the end of step $(n,i-1)$ there are at least 2 new gap sites which are within distance to each other and such that the site $X_{n,i}(\omega)$ which is deleted in step $(n,i)$
is within distance 4 to a new gap site of $(A_{n,i-1},I_{n,i-1})(\omega).$ Moreover, $X_{n,i}(\omega)$ is a new gap site for $(A_{n,i},I_{n,i})(\omega).$ 

\vspace{7 mm}

Let us now define a function $f$ which describes the transition from numbering steps in the form $0,(1,1),(1,2),\ldots,(1,\kappa),(2,1),(2,2),\ldots,(2,\kappa), \ldots$
to enumerating them as  $0,1,2,3, \ldots,$ i.e., $f: \{(0,0)\} \cup {\bf N} \times \{1, \ldots, \kappa\} \rightarrow {\bf N}_0$ is defined by $f(0,0)=0$ and $f(n,j)=(n-1) \kappa +j$ for all $n \geq 1, 1 \leq j \leq \kappa.$ Note that if $f^{-1}=(f^{-1}_1,f^{-1}_2)$ denotes the inverse function, then for each $m \in {\bf N},$ $f^{-1}_1(m)$ denotes the period, which the step $m$ belongs to.

\vspace{7 mm}
  
Let the process $\Gamma$ be defined as 
\[ \Gamma_m=|C^I|(A_{f^{-1}(m)},I_{f^{-1}(m)}),\] 
i.e., $\Gamma_m$ denotes the number of new gap sites which are within distance 4 to some other new  gap site after $m$ steps are completed.
Thus the path $\Gamma(\omega)$ increases at time $f(\sigma(\omega))$ if $\sigma(\omega)< \infty.$ Let us define $V_0(\omega) = \inf \{m \geq 0: \Gamma_m(\omega) > 0\}$ and let  
$V_i(\omega) = \inf \{m>V_{i-1}(\omega): \Gamma_m(\omega) >
\Gamma_{m-1}(\omega) \}$ denote the $i$-th point of increase of $\Gamma(\omega)$ after time $V_0(\omega).$ In order to investigate the above sets we will look at the number of deletions between successive increases of $\Gamma$ and we will look at the number of additions between an increase and the next return to 0 for $\Gamma$.  

First note that if $\Gamma(\omega)$ increases in step $m$ then $m$ corresponds to a deletion time and by (\ref{choice of the point to be deleted}) we have $U_{m}(\omega) \leq 8(L^I_{f^{-1}(m-1)}/|A_{f^{-1}(m-1)}|)(\omega)$ since the point $X_{f^{-1}(m)}$ has to be within distance $4$ to a new gap site, but can not be adjacent to an old gap nor be an endpoint of $A_{f^{-1}(m-1)}.$ Clearly $|A_{f^{-1}(m-1)}|(\omega) \geq |a| =N$ and if $m$ belongs to a period n such that $\tau^1(\omega) > n-1$  then $L^I_{f^{-1}(m-1)}(\omega) < (K^{-1}(K-1)+1)h(N)+K-1 \leq \ln N /8$ for all sufficiently large $N$. We thus get the following (assuming that $N$ is sufficiently large).

\begin{remark} \label{remark; points of increase} If $\Gamma(\omega)$ increases at time $m$ and if $\tau^1(\omega) > f^{-1}_1(m)-1$ then  $U_{m}(\omega) \leq \ln N/N.$
\end{remark}
Roughly speaking, the number of deletions between
two successive increases for $\Gamma$ is bounded below by a
geometric random variable with parameter $\ln N/N$ as long as we have not yet
reached period $\tau^1.$

On the other hand, being interested in points of decrease of $\Gamma,$ we can conclude from (\ref{choice of the new point to be added}) that if $m$ is an addition time with $U_m(\omega) > 1-(\Gamma(m-1)/|\partial A_{f^{-1}(m-1)}|)(\omega),$ then $\Gamma(\omega)$ decreases at time $m.$ (The reader may notice that $\Gamma(\omega)$ may also decrease in a deletion step. This happens for instance when a particle gets removed which separates an old and a new gap belonging to $C^I(A_{f^{-1}(m-1)})(\omega)$). Now, clearly  $|\partial A_{f^{-1}(m-1)}|(\omega) \leq 2+2(G(a)-G^I(a))+L^I_{f^{-1}(m-1)}.$ So if $\tau^1(\omega) > f^{-1}_1(m)-1$  then $|\partial A_{f^{-1}(m-1)}|(\omega) \leq 2+20K+(K^{-1}(K-1)+1)h(N)+K-1 \leq  \ln N $ for all sufficiently large $N$. (Here we also used that, given our initial pair $(a,\lambda)$ as above, the number of old gaps stays always bounded above by $10K$). We thus get the following (assuming that $N$ is sufficiently large).

\begin{remark} If $m$ is an addition time with $\Gamma_{m-1}(\omega) > 0$ and $\tau^1(\omega) > f^{-1}_1(m)-1$ and if $U_m(\omega) > 1-1/{\ln N}$ then $\Gamma(\omega)$ decreases at time $m.$
\end{remark}

Next let us observe,that if $\Gamma_{V_i-1}(\omega)=0$ then we need at most two points of decrease for $\Gamma(\omega)$ in the interval
$(V_{i}(\omega),V_{i+1}(\omega))$ to return to the value $0.$ This holds because $\Gamma_{V_i}(\omega) \leq 3$ if $\Gamma_{V_i-1}(\omega)=0.$ So, roughly speaking, the number of additions between $V_i$ and the next return to $0$ for $\Gamma$ is
bounded above by the sum of two geometric
random variables with parameter $1/{\ln N}$ as long as we have not yet reached $\sigma$ or period $\tau^1.$

In the next subsection we will be working on making this statement and the corresponding statement for the waiting times between two successive increases of $\Gamma$ more precise.

\subsection{The Processes $U^X$ and $U^Y$ and the Random Times $V^X_i$ and $V^Y_i$}
 
Let us assume that on our probability space $(\Omega, {\cal F}, P)$ processes $U^{(2)}$ and $U^{(3)}$ are defined which are independent of $(U,U^{(1)})$  and such that the
random variables $U_n^{(i)},$ $i=2,3, n \in {\bf N},$ are independent and uniformly distributed on $(0,1).$ Using $U$ and $U^{(2)}$ we will construct a process $U^X$ and using $U$ and $U^{(3)}$ we will construct a process $U^Y$ in a path by path way. So let us first construct $U^X.$

\vspace{7 mm}

Let  ${\tilde \zeta}=\sigma \wedge (\tau^1 ,\kappa)$ and let $\zeta=f({\tilde \zeta}).$ Let us call a  deletion time $m$ with $U_{m}(\omega) 
\leq \ln N/N$  a d-success for $U(\omega).$  For each $\omega \in \Omega$ we construct the sequence $U^X(\omega)$ by pasting together parts of $U(\omega)$ and of $U^{(2)}(\omega).$ Roughly speaking, starting at $V_0(\omega)$ we select successive elements from the sequence $U(\omega)$ (but only those which correspond to
deletions) until we reach the next d-success, assuming that up to this time we have not yet reached $\zeta(\omega)$. We then ignore all the elements of $U(\omega)$ which lie between this d-success and $V_1(\omega)$, assuming that $V_1(\omega)<\zeta(\omega).$ 
(Note that by Remark \ref{remark; points of increase},
$V_1(\omega)$ is a d-success for $U(\omega)$ if $V_1(\omega) \leq (\tau^1(\omega),\kappa).$ Yet, not every d-success for $U(\omega)$ is a point of increase of $\Gamma(\omega)$. So there may exist a d-success for $U(\omega)$ which lies strictly prior to $V_1(\omega)$).
Restarting at $V_1(\omega)$ we continue selecting elements from $U(\omega)$
until we reach the next following d-success, still assuming that this time also lies before $\zeta(\omega)$. (To be a little bit more specific: If the first d-success after $V_0(\omega)$ lies strictly before $V_1(\omega),$ say at time $m$, then $U_{m}(\omega)$ gets included into $U^X(\omega)$ but
not $U_{V_1}(\omega).$ The next value to be selected for the sequence
$U^X(\omega)$ is $U_{d(V_1)}(\omega),$ where $d(V_1)(\omega)$ corresponds to
the first deletion after $V_1(\omega),$ once again assuming that up to this time $\zeta(\omega)$ has not yet been reached). We continue in this way until
we reach $\zeta(\omega).$ At this time we are no longer interested in
the sequence $U(\omega).$ All further elements of the sequence get selected from  
the
sequence $U^{(2)}(\omega).$  

\vspace{7 mm}

More precisely, and to show that the above construction fits into the scheme described in Proposition \ref{prop:good selection rules} let us rewrite the above procedure. For any random time R define $t(R)(\omega)=\infty,$ if $R(\omega)=\infty,$ and $t(R)(\omega)=\min \{m \geq R(\omega): \mbox{$(m+1)$ is a deletion time} \},$ if $R(\omega)<\infty$. 
Let $S_1(\omega)=\zeta(\omega) \wedge t(V_0)(\omega)$.  

Now let $n \geq 1$ and let us assume that the first $n-1$ elements of the sequence $U^X(\omega)$ have already been constructed and that we defined integers $S_n(\omega)$ and $T_n(\omega),$ these integers denoting the number of elements of $U(\omega)$ (resp. $U^{(2)}(\omega)$) we investigated so far. 

If  
$S_n(\omega)=\zeta(\omega)$ let $U^X_n(\omega)=U^{(2)}_{T_n+1}(\omega),$ $T_{n+1}(\omega)=T_n(\omega)+1$ and $S_{n+1}(\omega)=S_n(\omega).$ 

Otherwise, choose an integer $l \geq 0$ with $V_l(\omega) \leq S_n(\omega) < V_{l+1}(\omega),$
let $U^X_n(\omega)=U_{S_n+1}(\omega),$ $T_{n+1}(\omega)=T_n(\omega)$ and define $S_{n+1}(\omega)$ as follows: 

If $S_n(\omega)+1$ is a d-success for $U(\omega),$ let $S_{n+1}(\omega)=\zeta(\omega) \wedge t(V_{l+1})(\omega),$  and if  $S_n(\omega)+1$ is no d-success for $U(\omega),$ let
\[S_{n+1}(\omega) = \left \{ \begin{array}{ll}
                              S_n(\omega)+1 & \mbox{ if                         $S_n(\omega)+2$ is a deletion time} \\
                              \zeta(\omega) \wedge (S_n(\omega)+1+K) &   \mbox{ otherwise.}                                              
                           \end{array}
                     \right. \]
\vspace{7 mm}

Note that the above procedure fits into the scheme described in Proposition \ref{prop:good selection rules}. Moreover, by construction we have a useful comparison between waiting times in the original process and the corresponding waiting times for the new process $U^X.$ Define $V^X_0=0$ and let $V^X_i(\omega) = \inf \{n > V^X_{i-1}(\omega): U^X_{n}(\omega) \leq \ln N/N \},$ i.e., $V^X_i$ denotes the $i-$th d-success for $U^X.$  Let ${\cal M}(\omega)=\sup \{i \geq 0:V_i(\omega) < \zeta(\omega) \},$ (with the convention that $\sup \emptyset = 0$), thus ${\cal M}$ denotes the number of
points of increase of $\Gamma$ in the time interval $(V_0(\omega),\zeta(\omega))$. Finally, for any time interval $\cal{C}$ let $\varphi_d(\cal{C})$ denote the number of deletion times in $\cal{C}.$ (Here we refer to our partitioning of time into addition and deletion times). We then get the following (assuming that $N$ is sufficiently large). 

\begin{proposition} \label{comparison for the number of deletions between successive increases of Gamma}
The random variables $U^X_n$ are independent and uniformly distributed on $(0,1)$ and the waiting times $V^X_i-V^X_{i-1}, i \geq 1,$ are independent geometric random variables with parameter $\ln N/N.$
Moreover, for all $i < {\cal M}(\omega)+1$ we have
\begin{equation}
(V^X_i-V^X_{i-1})(\omega) \leq \varphi_d((V_{i-1}(\omega),V_i(\omega)]),
\end{equation}
and if $V_{{\cal M}+1}(\omega)=\zeta(\omega)<\infty$ then 
\begin{equation}
(V^X_{{\cal M}+1}-V^X_{\cal M})(\omega) \leq \varphi_d((V_{\cal M}(\omega),V_{{\cal M}+1}(\omega)]).
\end{equation}

\end{proposition}

Let us now construct $U^Y$ from $U$ and $U^{(3)}$. Let us define an a-success for a sequence $u$ as a time $k$ with $u_k > 1-1/{\ln N}.$ For each $\omega \in \Omega$ we construct the sequence $U^Y(\omega)$ by pasting together parts of $U(\omega)$ and of $U^{(3)}(\omega).$ Our procedure can essentially be
described as follows. 

Starting at $V_0(\omega)$ we select successive
elements from the sequence $U(\omega)$ (but only those which correspond to
addition times) until we reach the first point of decrease for
$\Gamma(\omega)$ (assuming that this time lies before $\zeta(\omega)$).
(The reader may recall that this time may be an addition time or a deletion time). If
this decrease of $\Gamma(\omega)$ corresponds to an addition time which is an
a-success of $U(\omega),$ we now already have constructed $U^Y(\omega)$ up to its first
a-sucess. If not, we select succesive elements from the sequence
$U^{(3)}(\omega)$ until we reach its first a-success. Thus, in either case,
$U^Y(\omega)$ has now been constructed up to its first a-success. We now proceed to
construct $U^Y(\omega)$ up to its second a-success. So, say, $U_{m}(\omega)$ was
the last element of the sequence $U(\omega)$ which was investigated so far.
If $\Gamma_{m}(\omega) \neq 0$ we iterate the above procedure,
i.e., we once again select elements from $U(\omega)$ until $\Gamma(\omega)$
decreases for the next time - assuming that this time lies before
$\zeta(\omega)$ - and, if necessary, we fill in with elements from
$U^{(3)}(\omega)$ until we reach the next a-success for $U^Y(\omega).$
If $\Gamma_m(\omega) = 0$ we immediately select elements from
$U^{(3)}(\omega)$ until we reach the next a-success. 
Thus, in either case, $U^Y(\omega)$ has now been constructed up to its second a-success.
Replacing $V_0(\omega)$ by $V_1(\omega)$ we now iterate the whole procedure
and we continue in this way until we reach $\zeta(\omega).$ (The reader
may notice that if $\Gamma(\omega)$ does not drop down to $0$ in the
interval $(V_{i}(\omega),V_{i+1}(\omega))$ then $f(\sigma(\omega))$ and hence 
$\zeta(\omega)$ happens at the latest at time $V_{i+1}(\omega)$).
After time $\zeta(\omega)$ we are no longer interested in
the sequence $U(\omega).$ All further elements of the sequence $U^Y(\omega)$ are
selected successively from the
sequence $U^{(3)}(\omega).$

\vspace{7 mm}

More precisely, and to show that the above construction fits into the scheme described in Proposition \ref{prop:good selection rules} let us once again rewrite the above procedure. For any random time R define $s(R)(\omega)=\infty,$ if $R(\omega)=\infty,$ and $s(R)(\omega)=\min \{m \geq R(\omega): \mbox{$m$ is a point of decrease of $\Gamma(\omega)$ or  $(m+1)$ is an addition time} \},$ if $R(\omega)<\infty$. 
Let $S_1(\omega)=\zeta(\omega) \wedge s(V_0)(\omega),$ $T_1(\omega)=0.$
Now let $n \geq 1$ and let us assume that the first $n-1$ elements of the sequence $U^Y(\omega)$ have already been constructed and that we defined integers $S_n(\omega)$ and $T_n(\omega),$ these integers denoting the number of elements of $U(\omega)$ (resp. $U^{(3)}(\omega)$) we investigated so far.

If  $S_n(\omega)=\zeta(\omega)$ let $U^Y_n(\omega)=U^{(3)}_{T_n+1}(\omega),$ $T_{n+1}(\omega)=T_n(\omega)+1$ and $S_{n+1}(\omega)=S_n(\omega).$ 

Otherwise, choose an integer $l \geq 0$ with $V_l(\omega) \leq S_n(\omega) < V_{l+1}(\omega)$ and proceed as follows: \newline
1) If $\Gamma(\omega)$ decreases at time $S_n(\omega)$ and ($n=1$ or $\Gamma_{S_n}(\omega) = 0$ or $U^Y(\omega)$ has at most $2l$ a-successes among times $1, \ldots, n-1$, where  $n \geq 1$),
let $U^Y_n(\omega)=U^{(3)}_{T_n+1}(\omega),$ $T_{n+1}(\omega)=T_n(\omega)+1$ and define $S_{n+1}(\omega)$ as follows:
 
a) If $n$ is no a-success for $U^Y(\omega)$ let $S_{n+1}(\omega)=S_n(\omega);$ 

b) If $n$ is an a-success for $U^Y(\omega)$ and if $\Gamma_{S_n}(\omega) \neq 0$ let 
\[S_{n+1}(\omega) = \left \{ \begin{array}{ll}
                                S_n(\omega) & \mbox{ if 
$S_n(\omega)+1$ is an addition time} \\                                                                                  
                                  \zeta(\omega) \wedge s(S_n+1)(\omega) & \mbox{ otherwise}
                             \end{array}
                    \right.  ;\]
 
c) If $n$ is an a-success for $U^Y(\omega)$ and if $\Gamma_{S_n}(\omega) = 0$ let 
\[S_{n+1}(\omega) = \left \{ \begin{array}{ll}
                                S_n(\omega) & \mbox{ if                                                                                   $U^Y(\omega)$ has at most $2l+1$ a-successes} \\                                                    & \; \; \; \; \; \; \; \; \; \; \; \; \; \; \; \; \; \; \;  \mbox{among times $1, \ldots, n$ } \\ 
                                \zeta(\omega) \wedge                                     s(V_{l+1})(\omega) & \mbox{ otherwise}
                             \end{array}
                    \right.  .\]
2) Otherwise, (i.e., (1) does not apply),   
let $U^Y_n(\omega)=U_{S_n+1}(\omega),$ $T_{n+1}(\omega)=T_n(\omega)$ and define $S_{n+1}(\omega)$ as follows:

a) If $\Gamma(\omega)$ does not decrease at time $S_n(\omega)+1$ let $S_{n+1}(\omega)=  \zeta(\omega) \wedge s(S_n+1)(\omega)$;

b) If $\Gamma(\omega)$ decreases at time $S_n(\omega)+1$ and if $S_n(\omega)+1$ is no a-success for $U(\omega)$ 

$ \; \; \; \;$ let $S_{n+1}(\omega)=S_n(\omega)+1;$ 
 
c) If $\Gamma(\omega)$ decreases at time $S_n(\omega)+1$ and if $S_n(\omega)+1$ is an a-success for $U(\omega)$

$\; \; \; \; $ consider the following cases:

$ \; \; \; \alpha)$ If $\Gamma_{S_n+1}(\omega) \neq 0$ let $S_{n+1}(\omega)=  \zeta(\omega) \wedge s(S_n+1)(\omega);$
 
$ \; \; \; \beta)$ If $\Gamma_{S_n+1}(\omega) = 0$ let 
\[S_{n+1}(\omega) = \left \{ \begin{array}{ll}
                                S_n(\omega)+1 & \mbox{ if                                                                                  
$U^Y(\omega)$ has at most $2l+1$ a-successes} \\                                                    & \; \; \; \; \; \; \; \; \; \; \; \; \; \; \; \; \; \; \;  \mbox{among times $1, \ldots, n$ } \\ 
                                \zeta(\omega) \wedge                                     s(V_{l+1})(\omega) & \mbox{ otherwise}
                             \end{array}
                    \right.  .\]

Note that the above procedure fits into the scheme described in Proposition \ref{prop:good selection rules}. Moreover, once again we have a useful comparison between waiting times in the original process and the corresponding waiting times for the new process $U^Y.$ Let $V^Y_0(\omega)=0$ and let $V^Y_i(\omega)$ denote the $i-$th a-success of $U^Y(\omega).$ 
For $i=1, \ldots, {\cal M}(\omega)$ let $W^{(1)}_i(\omega)$ denote the first time of decrease of $\Gamma(\omega)$ after time $V_{i}(\omega)$ and let $W^{(2)}_i(\omega)$ denote the smallest integer $m >V_{i}(\omega)$  such that $\Gamma_m(\omega) =0.$ Finally, for any time interval $\cal{C}$ let $\varphi_a(\cal{C})$ denote the number of addition times in $\cal{C}.$ (Here we refer once again to our partitioning of time into addition and deletion times). We then get the following (assuming that $N$ is sufficiently large). 

\begin{proposition} \label{comparison for the number of additions between V_i and the next return to 0 for Gamma}
The random variables $U^Y_n$ are independent and uniformly distributed on $(0,1)$ and the waiting times $V^Y_i-V^Y_{i-1}, i \geq 1,$ are independent geometric random variables with parameter $1/{\ln N}.$
Moreover, we have
\begin{eqnarray}
(V^Y_{2i+1}-V^Y_{2i})(\omega)  &\geq& \varphi_a((V_{i}(\omega),W^{(1)}_{i}(\omega)]), \\
(V^Y_{2i+2}-V^Y_{2i+1})(\omega) & \geq&    \varphi_a((W^{(1)}_{i}(\omega),W^{(2)}_{i}(\omega)]),
\end{eqnarray}
for $0 \leq i < {\cal M}(\omega),$ (assuming that the left hand side is well-defined), and if ${\cal M}(\omega) < \infty$ and $V_0(\omega) < \zeta(\omega)$ then 
\begin{eqnarray}
(V^Y_{2{\cal M}+1}-V^Y_{2 {\cal M}})(\omega) & \geq&  \varphi_a((V_{{\cal M}}(\omega),(\zeta \wedge W^{(1)}_{{\cal M}})(\omega)]),  \\
(V^Y_{2{\cal M}+2}-V^Y_{2{\cal M}+1})(\omega) & \geq&  \varphi_a(((\zeta \wedge W^{(1)}_{{\cal M}})(\omega),(\zeta \wedge W^{(2)}_{{\cal M}})(\omega)]),
\end{eqnarray}
(assuming that the left hand side is well-defined).
\end{proposition}

\begin{remark} It can be shown that $U^X$ and $U^Y$ are independent since $U^{(2)}$ and $U^{(3)}$ are independent, but this fact will not be needed for our estimates.
\end{remark}

\subsection{Set Relations and Probability Estimates Related to $V^X$ and $V^Y$}

Let us continue to use the set-up and the definitions of Section 2.8.1. Thus, in particular, $(\Omega, {\cal F},P)$ is a probability space on which independent processes $U$ and $U^{(i)}, i = 1,2,3,$ are defined, the random variables $U_n, n \in {\bf N},$ and $U^{(i)}_n, i = 1,2,3, \; n \in {\bf N}$, being independent and uniformly distributed on $(0,1)$.  
$(A,I)=(A,I)(a,\lambda)$ denotes the Markov process with initial pair $(a, \lambda)$  constructed from $U$ as in Section 2.5 and the processes $U^X$ and $U^Y,$ and thus the random times $V^X_i, i \in {\bf N},$ and $V^Y_i, i \in {\bf N},$
are as previously defined. The reader may recall that these processes depend pathwise on the initial pair $(a,\lambda)$. As before, let $N=|a|.$

Let us now define the set   
\begin{eqnarray*}
{\cal D} &=& \{V^X_{\lfloor{N^{2/5}}\rfloor} \geq \kappa N^{6/5}; V^X_i-V^X_{i-1} > 4K
\lceil{N^{1/10}}\rceil, i=1 ,\ldots, \lfloor{N^{2/5}}\rfloor; \\
& & \qquad V^Y_i-V^Y_{i-1} \leq  N^{1/10},i=1 ,\ldots,
2 \lfloor{N^{2/5}}\rfloor \}, 
\end{eqnarray*}
and let ${\cal E}^1$ and ${\cal E}^2$ be defined as in Section 2.7.
Then the following set inclusions hold.  

\begin{proposition}  There exists an integer $N_0$ such that for all initial pairs $(a,\lambda)$ with $a \in \cal{G},$ $\lambda=\emptyset$ and $N=|a| \geq N_0$ we have  
\begin{equation} \label{main set inclusion for tau1 in subsection 2.8.2}
 \{\bar{\sigma} \wedge \tau^1 > N^{6/5}, \sigma_0 < N^{6/5} \} \sqsupset ({\cal E}^1)^c \cap {\cal D} 
\end{equation}
and such that for all initial pairs $(a,\lambda)$ with $a \in {\cal S}_M,$ $\lambda=\{1, \ldots, G(a)\}$ and $N=|a| \geq N_0$ we have
\begin{equation} \label{main set inclusion for tauM,2 in subsection 2.8.2}
 \{\bar{\sigma} \wedge \tau^{M,2} > 2N \} \sqsupset ({\cal E}^2)^c \cap {\cal D}. 
\end{equation}
\end{proposition}
{\bf Proof.} Let $(a,\lambda)$ satisfy $a \in \cal{G}$ and $\lambda=\emptyset,$ or $a \in {\cal S}_M$ and $\lambda=\{1, \ldots, G(a)\}.$ Let $N=|a|.$
Our first goal is to show that if 
$\omega \in {\cal D} $ satisfies 
$\bar{\sigma}(\omega) \leq \tau^1(\omega)$ then $\bar{\sigma}(\omega) >
N^{6/5},$ if $N$ is sufficiently large. So let $\omega  \in {\cal D} $ satisfy $\bar{\sigma}(\omega) \leq
\tau^1(\omega)$ and assume that $\sigma(\omega)$ lies within the first
$\lfloor{N^{6/5}}\rfloor$ periods. Then ${\tilde \zeta}(\omega)=\sigma(\omega)$ and 
$V_{{\cal
M}+1}(\omega)=\zeta(\omega),$ i.e., $V_{{\cal M}}(\omega)$ is a point of
increase for $\Gamma(\omega)$ ($\Gamma(\omega)$ increases from 0 to the value
2 or 3), $\Gamma(\omega)$ increases again at time  $V_{{\cal M}+1}(\omega),$
but $\Gamma(\omega)$ does not drop down to 0 in the interval $(V_{{\cal
M}}(\omega),V_{{\cal M}+1}(\omega)].$ Let $n_1$ denote the period containing $V_{{\cal
M}}(\omega)$ and let $n_2$ denote the period containing $V_{{\cal
M}+1}(\omega).$ Let us first suppose that ${\cal M}(\omega)+1 \leq N^{2/5}.$
Then by Proposition \ref{comparison for the number of deletions between successive increases of Gamma} we have $(n_2-n_1+1)(K-1) \geq \varphi_d
((V_{{\cal M}}(\omega),V_{{\cal M}+1}(\omega)]) \geq (V^X_{{\cal M}+1}-V^X_{{\cal M}})(\omega) > 4K
N^{1/10}$ and hence
$ n_2-n_1 > 3 N^{1/10}. $
On the other hand, using Proposition \ref{comparison for the number of additions between V_i and the next return to 0 for Gamma} we have 
$2 N^{1/10} \geq (V^Y_{2{\cal M}+2}-V^Y_{2{\cal
M}})(\omega) \geq  \varphi_a( (V_{{\cal M}}(\omega),V_{{\cal M}+1}(\omega)]) =
(n_2-n_1)K $
 and thus we get a contradiction. So let us now suppose that ${\cal M}(\omega) +1 > N^{2/5}.$ Then 
$ f(\sigma(\omega)) = V_{{\cal M}+1}(\omega) > V_{\lfloor{N^{2/5}}\rfloor}(\omega)$ and by Proposition \ref{comparison for the number of deletions between successive increases of Gamma} and our definition of ${\cal D}$ we have
$ \varphi_d( (0,V_{\lfloor{N^{2/5}}\rfloor}(\omega)]) \geq
V^X_{\lfloor{N^{2/5}}\rfloor}(\omega) \geq \kappa N^{6/5}$ in contradiction to our
assumption that $\sigma(\omega)$ lies within the first  $\lfloor{N^{6/5}}\rfloor$ periods.
Thus, in either case we have now shown that any  $\omega \in {\cal D} $
with $\bar{\sigma}(\omega) \leq \tau^1(\omega)$ satisfies $\bar{\sigma}(\omega) > N^{6/5}$ and thus $\{\bar{\sigma} \leq N^{6/5},\bar{\sigma} \leq \tau^1 \} \cap {\cal D} = \emptyset,$ or, equivalently,
\begin{equation} \label{first set inclusion for tau1}
\{\bar{\sigma} \wedge \tau^1 \leq N^{6/5},\bar{\sigma} \leq \tau^1 \} \subset {\cal D}^c.
\end{equation}
Noting that
$\{\bar{\sigma} \leq 2N, \bar{\sigma} \leq \tau^{M,2} \}  \subset  
  \{\bar{\sigma} \leq N^{6/5}, \bar{\sigma} \leq \tau^{1} \}$
we can also conclude that 
\begin{equation} \label{first set inclusion for tauM,2}
\{\bar{\sigma} \wedge \tau^{M,2} \leq 2N, \bar{\sigma} \leq \tau^{2,M} \} \subset {\cal D}^c.
\end{equation}

For $(a,\lambda)$ with $a \in {\cal G},$ $\lambda = \emptyset$ and $N=|a|$ sufficiently large, we can now combine (\ref{first set inclusion for tau1}) with (\ref{set inclusion for E1}) and obtain for this choice of $(a,\lambda),$
$ \{\bar{\sigma} \wedge \tau^1 \leq N^{6/5} \} \sqsubset {\cal E}^1 \cup {\cal D}^c, $
or, equivalently,
\begin{equation} \label{intermediate set inclusion for tau1}
 \{\bar{\sigma} \wedge \tau^1 > N^{6/5} \} \sqsupset ({\cal E}^1)^c \cap {\cal D}. 
\end{equation}
Similarly, using (\ref{first set inclusion for tauM,2}) and (\ref{set inclusion for E2}) we immediately obtain (\ref{main set inclusion for tauM,2 in subsection 2.8.2}) for all
$(a,\lambda)$ with $a \in {\cal S}_M,$ $\lambda=\{1, \ldots, G(a) \}$ and $N=|a|$ sufficiently large.

Our final goal is to show that if $\omega \in {\cal D}$ satisfies 
$(\bar{\sigma} \wedge \tau^1)(\omega) > N^{6/5}$ then $\sigma_0(\omega) <
N^{6/5}.$ Combining this fact with (\ref{intermediate set inclusion for tau1}), the set inclusion (\ref{main set inclusion for tau1 in subsection 2.8.2})
then gets established.

So let $\omega \in  {\cal D}$  satisfy $(\bar{\sigma} \wedge
\tau^1)(\omega) > N^{6/5}$ and thus $ \zeta(\omega) > N^{6/5} \kappa.$ If $C^I_{\lceil{N^{11/10}}\rceil}(\omega)=\emptyset$ we obviously have
$\sigma_0(\omega) = \lceil{N^{11/10}}\rceil < N^{6/5}.$ So suppose that
$C^I_{\lceil{N^{11/10}}\rceil}(\omega) \neq \emptyset.$ Then there exists $i$, $i <
\lfloor{N^{2/5}}\rfloor,$ with $V_i(\omega) \leq \lceil{N^{11/10}}\rceil \kappa < V_{i+1}(\omega).$ 
(The reader may notice that we can exclude the case $i \geq \lfloor{N^{2/5}}\rfloor$ since in
this case we would have $V_{\lfloor{N^{2/5}}\rfloor}(\omega) \leq \lceil{N^{11/10}}\rceil \kappa$ and - using that
$\zeta(\omega)>N^{6/5} \kappa$ -
we can thus conclude that ${\cal M}(\omega) \geq \lfloor{N^{2/5}}\rfloor.$ Hence by Proposition \ref{comparison for the number of deletions between successive increases of Gamma} and our definition of ${\cal D}$ , 
$\varphi_d((0, V_{\lfloor{N^{2/5}}\rfloor}])(\omega)
\geq V^X_{\lfloor{N^{2/5}}\rfloor}(\omega) \geq \kappa N^{6/5}$
 in contradiction to
$V_{\lfloor{N^{2/5}}\rfloor}(\omega)) \leq \lceil{N^{11/10}}\rceil \kappa$). 
Let $m$ denote the period containing
$V_{i}(\omega).$ 

We now claim that $V_{i+1}(\omega)$ does not belong to any of the periods $m, \ldots, m+3 \lceil{N^{1/10}}\rceil.$ Suppose not. Then $V_{i+1}(\omega) \leq
(\lceil{N^{11/10}}\rceil+3\lceil{N^{1/10}}\rceil) \kappa < \lfloor{N^{6/5}}\rfloor \kappa < \zeta(\omega)$ and hence
$i+1 \leq {\cal M}(\omega).$ Thus by Proposition \ref{comparison for the number of deletions between successive increases of Gamma} we obtain $\varphi_d((V_i,V_{i+1}])(\omega)  \geq
(V_{i+1}^X-V_{i}^X)(\omega) > 4K N^{1/10}$
 which contradicts our assumption
that $V_{i+1}(\omega)$ belongs to one of the periods $m, \ldots, m+3 \lceil{N^{1/10}}\rceil.$ 

Next we claim that $W_i^{(2)}(\omega)$ happens 
within periods $m, \ldots, m+\lceil{N^{1/10}}\rceil.$ To see this note that if $W_i^{(2)}(\omega)$ does not happen 
within these periods, then, since $\zeta(\omega) > \lfloor{N^{6/5}}\rfloor \kappa,$ all additions in periods $m+1 ,\ldots, m+\lceil{N^{1/10}}\rceil$ contribute to $U^Y(\omega)$ and
there is at most one a-success among these additions. Hence 
$(V^Y_{2(i+1)}-V^Y_{2i})(\omega) > 2 \lceil{N^{1/10}}\rceil.$ Yet, since $\omega \in {\cal D}$ and $i+1 \leq N^{2/5}$ we have $(V^Y_{2(i+1)}-V^Y_{2i})(\omega) \leq 2N^{1/10}$ and we thus arrive at a contradiction.

Now recall that $\Gamma(\omega)$ drops
down to the value 0 at time $W_i^2(\omega)$ and $V_{i+1}(\omega)$ is the
first point of increase of $\Gamma(\omega)$ after time $V_{i}(\omega).$ Thus our two claims above imply that $C^I(\omega)$ is empty in periods $m+\lceil{N^{1/10}}\rceil+1 ,\ldots,
m+3\lceil{N^{1/10}}\rceil.$ Finally using that $m \leq \lceil{N^{11/10}}\rceil \leq W_i^2(\omega) \leq
m+\lceil{N^{1/10}}\rceil,$ (the second inequality holds by our choice of
$V_i(\omega)$ and since $C^I_{\lceil{N^{11/10}}\rceil}(\omega) \neq \emptyset$), we can conclude that period
$m+\lceil{N^{1/10}}\rceil$ lies strictly before period $\lfloor{N^{6/5}}\rfloor$ for all sufficiently large $N$ and $m+\lceil{N^{1/10}}\rceil \geq \lceil{N^{11/10}}\rceil.$ Hence $\sigma_0(\omega) < N^{6/5}.$ 
\newline
\vspace{7 mm}
\fine

In Section 2.7 we proved upper bounds for the probability of the events ${\cal E}^1$ and ${\cal E}^2$ for those  initial pairs $(a,\lambda)$ that we are interested in. We will now derive a suitable upper bound for $P({\cal D}^c(a,\lambda))$. 

\begin{proposition} 
Let $(a,\lambda)$ be an initial pair with $a \in {\cal G}, \lambda=\emptyset$ or $a \in {\cal S}_M, \lambda=\{1, \ldots, G(a)\}$ and let $N=|a|.$ For all sufficiently large $N$ we have
\begin{equation} \label{upper bound for the probability of the complement of D}
P^{a,\lambda}({\cal D}^c) \leq 7N^{-2/5}. 
\end{equation}
\end{proposition}
{\bf Proof.} Let $X_i=V_i^X-V_{i-1}^X.$ Since by Proposition \ref{comparison for the number of deletions between successive increases of Gamma}  the
random variables $X_i$ are independent geometric random variables with
parameter $p={\ln N}/N$ we have  $E X_1=1/p$ and $\mbox{Var}(X_1)=(1-p)/{p^2} \leq
1/{p^2}.$ Hence for all sufficiently large $N$ we have
\begin{eqnarray*}
P(V^X_{\lfloor{N^{2/5}}\rfloor} < \kappa N^{6/5}) & \leq & P(\lfloor{N^{2/5}}\rfloor^{-1} \sum_{i=1}^{\lfloor{N^{2/5}}\rfloor} X_i
< \kappa N^{4/5})\\ 
& = & P(\lfloor{N^{2/5}}\rfloor^{-1}  \sum_{i=1}^{\lfloor{N^{2/5}}\rfloor} X_i - E X_1 
< \kappa N^{4/5} - N/{\ln N})\\
& \leq & P( | \lfloor{N^{2/5}}\rfloor^{-1} \sum_{i=1}^{\lfloor{N^{2/5}}\rfloor} X_i - E X_1  | 
 > N/{(2 \ln N)}). 
\end{eqnarray*} 
Now $E(\lfloor{N^{2/5}}\rfloor^{-1}
\sum_{i=1}^{\lfloor{N^{2/5}}\rfloor} X_i)=E X_1$ and $\mbox{ Var}(\lfloor{N^{2/5}}\rfloor^{-1}
\sum_{i=1}^{\lfloor{N^{2/5}}\rfloor} X_i) =\lfloor{N^{2/5}}\rfloor^{-1} \mbox{ Var}(X_1).$ Thus by Chebychev's inequality 
we get
\begin{eqnarray}
P(V^X_{\lfloor{N^{2/5}}\rfloor} < \kappa N^{6/5}) & \leq & \lfloor{N^{2/5}}\rfloor^{-1} \mbox{Var}(X_1) \;\; (N/{(2 \ln N)})^{-2} \nonumber \\
& \leq & 5 N^{-2/5} \nonumber.
\end{eqnarray} 
Next note that
\begin{eqnarray*}
  P(V_i^X-V_{i-1}^X \leq 4K \lceil{N^{1/10}}\rceil) &=& 1- (1-{\ln N}/N)^{4K \lceil{N^{1/10}}\rceil}\\
&=& 1-\left(
(1-{\ln N}/N)^{N/{\ln N}} \right ) ^{4K \lceil{N^{1/10}}\rceil{\ln N}/N}, 
\end{eqnarray*}
and for all sufficiently large $N$ we have $(1-{\ln N}/N)^{N/{\ln N}} \geq
\exp(-2).$ Thus, for all sufficiently large $N$, 
\[ P(V_i^X-V_{i-1}^X \leq 4K \lceil{N^{1/10}}\rceil ) \leq 1-\exp(-N^{-8/10}) \leq
N^{-8/10}. \]
Similarly,using that $V^Y_i-V^Y_{i-1}$ are geometric random variables with parameter $1/{\ln N},$ we have 
 \[ P(V_i^Y-V_{i-1}^Y > N^{1/10}) = (1-1/{\ln N})^{\lfloor{N^{1/10}}\rfloor} \leq
\exp(-N^{1/20}) \]
for all sufficiently large $N.$ Since 
${\cal D}^c = \{ V^X_{N^{2/5}} <
\kappa N^{6/5}\} \cup \bigcup_{i=1}^{\lfloor{N^{2/5}}\rfloor} \{ V^X_i-V^X_{i-1} \leq 4K
N^{1/10} \}
\cup \bigcup_{i=1}^{2\lfloor{N^{2/5}}\rfloor} \{V^Y_i-V^Y_{i-1} >  N^{1/10} \}, $
the proposition follows immediately from the previous inequalities.
\newline
\vspace{7 mm}
\fine

\section{Proof of Lemmas \ref{first lemma} and \ref{second lemma}}

We will first prove Lemma \ref{second lemma} and we will then proceed to a proof of Lemma \ref{first lemma}. In the following we continue to work on the same probability space as in Section 2.8.2. 

\vspace{7 mm}

{\bf Proof of Lemma \ref{second lemma}}.
Let $M >K$ be a fixed positive integer  and let our initial pair  $(a, \lambda)$ satisfy $a \in {\cal S}_M$ and $\lambda=\{1, \ldots, G(a)\}$. Let $N=|a|.$  
Let $(A,I)=(A,I)(a,\lambda)$ be the Markov process with initial pair $(a,\lambda)$ constructed from the process $U$ as in Section 2.5. 
  
Clearly, $\{\bar{\sigma} \wedge \tau^{M,2} > 2N,  \tau_{M,2}< 2N \} = \{ \bar{\sigma} \wedge \tau^{M,2} > 2N \} \setminus \{ \bar{\sigma} \wedge \tau^{M,2} > 2N,
\tau_{M,2} \geq 2N \}$ and combining (\ref{main set inclusion for tauM,2 in subsection 2.8.2}) with (\ref{estimate for calE2}) and (\ref{upper bound for the probability of the complement of D}) yields (for sufficiently large $N$) $P^{a,\lambda}( \bar{\sigma} \wedge \tau^{M,2} > 2N ) \geq 1-N^{-5/(8M)}-7N^{-2/5}.$ So the proof of Lemma \ref{second lemma} gets completed if we can show that (for sufficiently large $N$) 
\[ 
P^{a,\lambda}(\bar{\sigma} \wedge \tau^{M,2} > 2N, \tau_{M,2} \geq 2N ) \leq
N^{-1}.\]
Let us emphasize that for our special construction the above stopping times depend pathwise on the initial pair $(a,\lambda)$ and hence this inequality can be rewritten as   
\begin{equation} \label{ineq in the proof or Lemma 2}
P((\bar{\sigma} \wedge \tau^{M,2})(a,\lambda) > 2N, \tau_{M,2}(a,\lambda) \geq 2N ) \leq
N^{-1}.
\end{equation}
To prove this inequality we first claim that,  given the initial pair $(a,\lambda),$ the following set inclusion holds
\begin{equation} \label{set inclusion, proof of Lemma 2}
\{\bar{\sigma} \wedge \tau^{M,2} > 2N, \tau_{M,2} \geq 2N \} \subset \{\sum_{n=N}^{2N-1} \sum_{i=1}^K {\bf 1}_{ \{ U_{n \kappa +i} \leq c_1/h(N) \}} \geq N/2 \},
\end{equation}
where $c_1=4M.$
So let $\omega \in \{ \bar{\sigma} \wedge \tau^{M,2} > 2N, \tau_{M,2} \geq
2N \} .$ Then $L^I_{N}(\omega) \leq (3/M+(K-1)/K)h(N) \leq \ln N$ and hence 
\begin{equation} \label{lower bound for the increment of LI}
(L^I_{2N} - L^I_{N})(\omega) \geq -\ln N.
\end{equation}
Let us now use (\ref{general ineq for increments of LI}) to bound the increments $(L^I_{n}-L^I_{n-1})(\omega)$ from above. 
By choice of the initial pair $(a, \lambda)$ we have $|\partial^o(A_{n,i-1},I_{n,i-1})|(\omega)=0$ for all $n \in {\bf N}, i=1, \ldots, K.$ Moreover, if $\sigma(\omega) > (n,i-1)$ then $|\partial A_{n,i-1}|(\omega) \geq L^I(A_{n,i-1},I_{n,i-1})(\omega)$ and if $\tau_{M,2}(\omega) \geq 2N$ then for $N<n \leq 2N, 1 \leq i \leq K,$ we have $L^I(A_{n,i-1},I_{n,i-1})(\omega) > h(N)/M-K \geq h(N)/(2M),$ if $N$ is sufficiently large. Hence, if $N$ is sufficiently large, then  for $\omega \in \{\bar{\sigma} \wedge \tau^{M,2} > 2N, \tau_{M,2} \geq 2N\}$ and $N \leq n < 2N$ we have 
\[(L^I_{n+1} - L^I_{n})(\omega) \leq -1
+ \sum_{i=1}^K {\bf 1}_{ \{ U_{n \kappa + i} \leq c_1/h(N) \} } \]
(with $c_1$ defined as above), and  summing over $n$ we get  
\[ (L^I_{2N}(\omega) - L^I_{N})(\omega) \leq - N +
\sum_{n=N}^{2N-1} \sum_{i=1}^K {\bf 1}_{ \{ U_{n \kappa + i} \leq c_1/h(N) \}}. \]
Combining this last inequality with (\ref{lower bound for the increment of LI}) we can conclude that (for $N$ sufficiently large) any $\omega \in \{ \bar{\sigma} \wedge \tau^{M,2} > 2N, \tau_{M,2} \geq
2N \}$ satisfies 
\[\sum_{n=N}^{2N-1} \sum_{i=1}^K {\bf 1}_{ \{ U_{n \kappa + i} \leq c_1/h(N) \}}(\omega) \geq N/2, \]
i.e., (\ref{set inclusion, proof of Lemma 2}) holds.
To prove (\ref{ineq in the proof or Lemma 2}) it thus suffices to bound the probability of the event on the right hand side of (\ref{set inclusion, proof of Lemma 2}) from above by $N^{-1}$. Clearly, $\mbox{ Var}({\bf 1}_{ \{ U_{n \kappa + i} \leq c_1/h(N) \}}) \leq c_1/h(N)$ since the random variables $U_n, n \in {\bf N},$ are uniformly distributed on $(0,1)$. Thus, using the independence of the random variables $U_n, n \in {\bf N},$  
and  Chebychev's inequality we get 
\begin{eqnarray*}
\lefteqn{P ( \sum_{n=N}^{2N-1} \sum_{i=1}^K  {\bf 1}_{ \{ U_{n \kappa + i
} \leq  c_1/h(N) \}}  \geq  N/2 ) \leq } \\
& \leq & P \left( | (KN)^{-1} \sum_{n=N}^{2N-1} \sum_{i=1}^K {\bf 1}_{
\{ U_{n \kappa + i }\leq c_1/h(N) \}} - c_1/h(N) | \geq 1/(4K) \right) \\
&  \leq &  c_1 (4K)^2 /(KN h(N)) \\
&  \leq & N^{-1},
\end{eqnarray*}
which completes the proof of (\ref{ineq in the proof or Lemma 2}).
\newline
\vspace{7 mm}
\fine

\vspace{7 mm}

{\bf Proof of Lemma \ref{first lemma}}. 
Let $M>K.$ Let $a \in \cal{G}$ and let $N=|a|.$ Let $(A,I)=(A,I)(a,\emptyset)$ be the Markov process with initial pair $(a,\emptyset)$ constructed from the process $U$ as in Section 2.5. 

Our first goal is to show that for all sufficiently large $N$ we have
\begin{equation} \label{first inequality in 'Section 2.10', only theta0 is still missing}
 P^{a,\emptyset}(\bar{\sigma} \wedge \tau^1 > N^{6/5}, \sigma_0 < N^{6/5}, \tau_M < N^{6/5}) \geq 1-8N^{-2/5}.
\end{equation}
To see this, note that $\{\bar{\sigma} \wedge \tau^1 > N^{6/5}, \sigma_0
<N^{6/5}, \tau_M < N^{6/5} \} = \{ \bar{\sigma} \wedge \tau^1 > N^{6/5},
\sigma_0 <N^{6/5} \} \setminus \{ \bar{\sigma} \wedge \tau^1 > N^{6/5},
\tau_M \geq N^{6/5} \}$  and  (for sufficiently large $N$) $P^{a,\emptyset}( \bar{\sigma} \wedge \tau^1 > N^{6/5},
\sigma_0 <N^{6/5} ) \geq 1-N^{-3/5}-7N^{-2/5}$ by using the set inclusion (\ref{main set inclusion for tau1 in subsection 2.8.2}) and the inequalities (\ref{estimate for calE1}) and (\ref{upper bound for the probability of the complement of D}). Thus (\ref{first inequality in 'Section 2.10', only theta0 is still missing}) follows if we can show that for all sufficiently large $N$
\begin{equation} \label{second inequality in 'Section 2.10' involving sigma,tau1 and tauM}
 P^{a,\emptyset}( \bar{\sigma} \wedge \tau^1 > N^{6/5}, \tau_M \geq N^{6/5}) \leq N^{-1}. 
\end{equation} 
Now, similar to the proof of Lemma \ref{second lemma} it can be shown that 
\[ \{ (\bar{\sigma} \wedge \tau^1)(a,\emptyset) > N^{6/5}, \tau_M(a,\emptyset) \geq N^{6/5} \} \subset \{ \sum_{n=N^{11/10}+1}^{N^{11/10}+N} \sum_{i=1}^K  {\bf 1}_{ \{ U_{n,i
} \leq  c_1^*/h(N) \}}  \geq  N/2 \} ,\]
where $c_1^*=44KM.$ 
By the independence of the random variables $U_n,$ the probability of this last event is bounded above by $N^{-1}.$ 
The inequality (\ref{second inequality in 'Section 2.10' involving sigma,tau1 and tauM}), and hence (\ref{first inequality in 'Section 2.10', only theta0 is still missing}), thus gets established.

We are now ready to investigate  the event $\{ \bar{\sigma}
\wedge \tau^1 > N^{6/5}, \sigma_0 <N^{6/5},\tau_M < N^{6/5}, \theta_0 \leq N^{11/10} \}$ for the pocess $(A,I)(a,\emptyset)$ . Obviously we have $ \{ \bar{\sigma}
\wedge \tau^1 > N^{6/5}, \sigma_0 <N^{6/5}, \tau_M < N^{6/5}, \theta_0 \leq N^{11/10} \} = 
  \{ \bar{\sigma} \wedge
\tau^1 > N^{6/5}, \sigma_0 < N^{6/5}, \tau_M < N^{6/5} \}  
\setminus  \{ \bar{\sigma}
\wedge \tau^1 > N^{6/5}, \theta_0 > N^{11/10} \}.$
Thus by (\ref{first inequality in 'Section 2.10', only theta0 is still missing}) the proof of the lemma gets completed if
we can show
that for all $a \in {\cal G}$ with $N=|a|$ sufficiently large we have 
\begin{equation} \label{third inequality in 'Section 2.10'  involving sigmabar, tau1 and theta0} 
P^{a, \emptyset}(\{ \bar{\sigma}
\wedge \tau^1 > N^{6/5}, \theta_0 > N^{11/10} \}) \leq K N^{-1/20}.
\end{equation}

Let $L_{n,i}^o$ denote the number of old gap sites of
$(A_{n,i},I_{n,i}).$ Note that $L^o$ can only increase when a
particle is deleted which is adjacent to an old gap but not an endpoint,  and $L^o$ decreases
when the site of a particle which is added belongs to an old gap. Thus,
by construction of our Markov chain, $L^o(\omega)$ can only increase in
step $(n,i),$ if $(n,i)$ is a deletion step and 
$U_{(n-1) \kappa+i}(\omega) > 1 - (|V^{o}|(A_{n,i-1},I_{n,i-1})/|A_{n,i-1}|)(\omega).$
Using  that the number of old gaps is nonincreasing and $|A_{n,i-1}|
\geq |a| = N,$ this last term can be bounded below by $1-20K/N.$  Moreover, if all new gaps of $A_{n,i-1}(\omega)$
have at most size two, then in any deletion step an increase of $L^o(\omega)$ is bounded above by 3. 

On the other hand, $L^o(\omega)$ will decrease in step $(n,i),$ if $(n,i)$ is an addition
time and if the new particle is added at a boundary point belonging to an
old gap, i.e., if  
$U_{(n-1) \kappa + i}(\omega) \leq (|\partial^o(A_{n,i-1},I_{n,i-1})|/|\partial
A_{n,i-1}|)(\omega).$
Now if the number of new gap sites is
bounded above by $(1+(K-1)/K)h(N)+K-1$ then the total number of boundary points
of  $A_{n,i-1}(\omega)$ is
bounded above by $\ln N$ for all sufficiently large  $N$. (Here we use that the set $a$ has at most $10K$ gaps). Hence in this
case $(|\partial^o(A_{n,i-1},I_{n,i-1})|/|\partial
A_{n,i-1}|)(\omega)$
is bounded below by $1/{\ln N}$ if
$A_{n,i-1}(\omega)$ has at 
least one old gap. 
So let us define the process $\hat{L}^o$ as
\[\hat{L}^o_n = \sum_{m=1}^n \left( -\sum_{i=1}^K {\bf 1}_{ \{ U_{(m-1) \kappa + i} \leq
1/{\ln N} \} } + 3 \sum_{i=K+1}^{\kappa} {\bf 1}_{ \{ U_{(m-1) \kappa + i} > 1-20K/N \}
} \right), \]
and let us define the random time
\[\hat{\tau} = \inf \{ n: \hat{L}^o_n \leq - (K-1) N \}.\]  Recalling that
$L(a) \leq (K-1)N$ the previous observations imply that (for $N$ sufficiently large)
\[ \{ \bar{\sigma}
\wedge \tau^1 > N^{6/5}, \theta_0 > N^{11/10} \} \subset \{ \hat{\tau} >
N^{11/10} \}. \]
Thus we have reduced the problem to bounding the probability of this last
event. So let ${\cal U}_m= \cup_{i=K+1}^{\kappa} \{U_{(m-1) \kappa +i} > 1-20K/N \},$
${\cal D}_m= \{U_{m,1} \leq 1/{\ln N} \} \cap \bigcap_{i=K+1}^{\kappa} \{U_{(m-1) \kappa +i}
\leq  1-20K/N \}$ and define $\tilde{L}_n^o = \sum _{m=1}^n \left( 3K {\bf
1}_{{\cal U}_m } -  {\bf 1}_{{\cal D}_m } \right).$ Clearly $\hat{L}^o_n
\leq \tilde{L}_n^o$ and thus 
\[ \{ \hat{\tau} >N^{11/10} \} \subset \{ \tilde{\tau} >N^{11/10} \},\]
 where $\tilde{\tau}= \inf \{ n: \tilde{L}^o_n
\leq - (K-1) N \}.$ Let $\sigma_1, \sigma_2, \ldots $ denote the jump times
of the process $\tilde{L}_n^o$ and let $Y_n=\tilde{L}_{\sigma_n}^o.$ Let 
$\tau^* = \inf \{ n: Y_n \leq - (K-1) N \}.$ Then
\begin{eqnarray*}
  \{ \tilde{\tau} > N^{11/10} \} & \subset &  \{ \tau^* > N^{21/20} \} \cup \{
\tau^* \leq N^{21/20}, \sigma_{N^{21/20}} > N^{11/10} \} \\
& \subset & \{ \tau^*
> N^{21/20} \} \cup \{\sigma_{N^{21/20}} > N^{11/10} \}.
\end{eqnarray*}
The waiting times $\sigma_i-\sigma_{i-1}$ are independent geometric random
variables with parameter $p_w$ satisfying $1/{\ln N} \leq p_w.$ Thus
$E \sigma_1 \leq {\ln N}$ and $\mbox{Var} (\sigma_1) \leq ({\ln N})^2.$ Hence for
all sufficiently large $N$ we have by Chebychev's inequality
\begin{eqnarray*}
P( \sigma_{N^{21/20}} > N^{11/10} ) &\leq& P(|N^{-21/20}
\sigma_{N^{21/20}}- E \sigma_1| > 2^{-1} {N^{1/20}}) \\ 
& \leq & 4 (\ln N)^2/(N^{21/20} N^{1/10})\\
& \leq & N^{-1}. 
\end{eqnarray*}
So (\ref{third inequality in 'Section 2.10'  involving sigmabar, tau1 and theta0}) gets established if we can show that 
\begin{equation} \label{fourth numbered inequality in 'Section 2.10'  involving the stopping time tau*}
P(\tau^* > N^{21/20}) \leq (K-1/2) N^{-1/20}
\end{equation}
for all sufficiently large $N.$ 
Note that $Y_n$ is a random walk with step law $\mu$ (depending on $N$)
satisfying $\mu(3K)+\mu(-1) = 1$ and $\mu(3K) \neq 0 \neq \mu(-1).$ Thus we
can choose $p,$ $0<p<1,$ such that any Bernoulli random walk $Z$ with
parameter $p$ satisfies 
\[P(Z_S = 3K) = \mu(3K), \]
where $S=\inf \{ n>0: Z_n \in \{-1,3K\} \}$ and we assume that $Z$ starts
from $0$. So let $Z$ be such a random walk. Define the stopping times
$S_n$ by $S_1= \inf \{n>0: Z_n \in \{-1,3K\} \}, S_{n+1}= \inf \{n>0:
Z_m-Z_{S_n} \in \{-1,3K\} \}$ and let $T=\inf \{ n: Z_{n} \leq  - (K-1) N
\}$ and $T^*=\inf \{ n: Z_{S_n} \leq  - (K-1) N
\}.$ Then 
$ \{T^* > N^{21/20} \} \subset \{T > N^{21/20} \} ,$ 
since $Z_T(\omega) =Z_{S_n}(\omega)$ for some index $n$ (depending on
$\omega$) with $n \leq T(\omega).$  Hence 
\[
P(\tau^* > N^{21/20}) = P(T^* > N^{21/20}) \leq  P(T > N^{21/20}) \leq ET/{N^{21/20}}.
\] 
Now recalling the definition of $\mu$ and using that ${\ln N}/N$ converges
to $0$ as $N \rightarrow \infty$ it is clear that $\mu(3K)$ - and hence $p$
- converges to 
$0$ as $N \rightarrow \infty.$ In particular $p<1/2$ for all sufficiently
large $N$ and so $ET=(K-1)N/{(1-2p)}$ (see \cite[page 318]{Feller}). Choosing $N$ sufficiently large, (\ref{fourth numbered inequality in 'Section 2.10'  involving the stopping time tau*})  follows and hence the proof of the lemma gets completed.
\newline
\vspace{7 mm}
\fine

\section{Completion of the Proof of Theorems \ref{first theorem} and \ref{second theorem} }

Note that if the process $A$ starts with an initial set $a$ belonging to ${\cal S}_M,$ ($M$ a fixed positive integer),   then for $\; \omega \in \{T^{M,2} \wedge
S \geq 2|a|, T_{M,2}<2|a|\}$ we have 
$A_{T_{M,2}}(\omega) \in {\cal S}_M \; , \; | A_{T
  _{M,2}}(\omega) | \geq 2 |a|,$ and up to time
$T_{M,2}(\omega)$ the number of unoccupied sites has never become too
big and the positions of these have remained fairly spread out. Since for any $M \in {\bf N}$ the family ${\cal S}_M$ of subsets of ${\bf Z}$  is recurrent we can now use a Borel-Cantelli argument to prove the following. 

\begin{proposition}
Let $A_0=\{0\}.$ With probability $1$ we have \newline
1) $\limsup_{n \rightarrow \infty} L(A_n)/h(n) \leq (K-1)/K,$ \newline
2) $G^{(2)}(A_n)=1$ i.o., \newline
3) $L(A_n) \leq G(A_n)+1$ eventually.
\end{proposition}
{\bf Proof.} First note that if $L(A_n)(\omega) \leq G(A_n)(\omega)+1$ then $A_n(\omega)$ has at most one gap of size bigger than $1,$ and if there is such an exceptional gap, then it has size $2.$ 
Thus (2) follows immediately from (3) and Proposition \ref{prop: at least one gap bigger than 1 i.o.}.
So it suffices to prove (1) and (3).
Fix $c>(K-1)/K$ and choose an integer $M\geq 3$ with $(K-1)/K+3/M<c$. Let ${\cal S} = {\cal S}_M$. Define 
$${\cal R} = \{ a \subset {\bf Z} : |a| < \infty ,
L(a)>G(a)+1 \; \mbox{or} \; L(a)>c \, h(|a|) \}.$$
We want to show that eventually $A$ stays outside of $\cal R$ . By Lemma \ref{second lemma}
we can choose  $N_0 \in {\bf N}$ such that for all $a \in {\cal S}_M$ with $|a|=N \geq N_0$ the inequality (\ref{inequality in the second lemma}) holds.
Define 
$T_0=\inf\{n \geq N_0: A_n \in {\cal S}\}$ and 
$T_i=\inf\{n \geq 2 T_{i-1}+1: A_n \in {\cal S}\},$ $i \in {\bf N}.$
By Corollary \ref{second corollary to Lemma1} the stopping times $T_i$ are finite with probability 1. 
Let 
$R_i=\inf\{n \geq T_{i-1}:A_n \in \cal R$$
 \}, i= 1,2,\ldots $ and let 
$T(N)=\inf\{n \geq N: A_n \in \cal S \}$ and $R = \inf\{n \geq 0: A_n \in \cal R \}$. We want to bound the probability of the event $\{R_i \leq T_i \}.$
By the Markov property of $A$ and using that $|A_{T_{i-1}}|=T_{i-1}+1$
we have $P(R_i \leq T_i \mid A_{T_{i-1}}=a)=P^a(R \leq T(|a|))$ for all $a \subset {\bf Z}$ with $P(A_{T_{i-1}}=a) > 0.$ Now observe that if the Markov process $A$ starts from an initial set $a \in {\cal S}$ of size $N$ then $\{T^{M,2} \wedge \bar S > 2N, T_{M,2} < 2N \} \subset \{R >T(|a|)\}$ and thus by Lemma \ref{second lemma}, $P^a(R  \leq T(|a|)) \leq |a|^{-1/(2M)}.$ Combining this estimate with the fact that $|A_{T_{i-1}}| \geq  2^{i-1} N_0$ we get
\begin{eqnarray*}
P(R_i \leq T_i) &=&  \sum_{a \subset {\bf Z}: P(A_{T_{i-1}}=a) > 0}
P^a( R \leq T(|a|)) P(A_{T_{i-1}}=a)\\
& \leq& \sum_{a \subset {\bf Z}: P(A_{T_{i-1}}=a) > 0}
(2^{i-1}N_0)^{-1/(2M)} P(A_{T_{i-1}}=a)\\
&=& (2^{i-1}N_0)^{-1/(2M)}.
\end{eqnarray*} 
Hence by the Borel-Cantelli lemma $\, P(R_i \leq T_i$ for infinitely many
i$)=0$, which implies (1) and (3) above.
\newline
\vspace{7 mm}
\fine

The previous proposition together with Corollary \ref{cor:no gaps of size bigger than two} proves Theorem \ref{second theorem}. As for Theorem \ref{first theorem}, the first part of this theorem is implied by part (1) of the previous proposition together with Proposition \ref{lower bound for limsup of L(A)/h}. Moreover, in Proposition \ref{prop: lower bound for lim inf of L} we established a lower bound for $\liminf_{n \rightarrow \infty} L(A_n)$. Thus the proof of Theroem \ref{first theorem} gets completed once we show the following.

\begin{proposition} 
\[P(L_n \leq K-2 \; i.o.) = 1. \]
\end{proposition}
{\bf Proof.} Let us construct the Markov process $A$ explicitly. (The construction corresponds to the one given in the proof of our two lemmas). Let $(\Omega, {\cal F}, P)$ be a probability space on which is defined a sequence  $\{U_j\}$ of independent random variables, which are uniformly
distributed on $(0,1).$ Define $A_0=\{0\}.$ Now suppose that $A_{n,i}$ has already been defined. In order to define $A_{n,i+1}$ it suffices to specify how 
the point $X_{n,i+1}$ which gets added to
(resp. deleted from) the set $A_{n,i}$ gets defined.

If $(n,i+1)$ is an addition time let us enumerate the boundary points 
of $A_{n,i}$ such that the outer boundary points $x_{1}=\min A_{n,i}-1$ and $x_{2}=\max A_{n,i} +1$ are followed by all other boundary points, those being enumerated in increasing order. Define
\[ X_{n,i+1}(\omega)=x_j(\omega)  \mbox{ if } (j-1)/{| \partial A_{n,i}|(\omega)} <
U_{(n-1)\kappa+i+1}(\omega) < j/{| \partial A_{n,i}|(\omega)} .\] 
Similarly, if $(n,i+1)$ is a deletion time let us first enumerate in increasing order the points of $A_{n,i}$ which are not endpoints and let us define $y_{| A_{n,i}|-1}= \min A_{n,i}$ and $y_{| A_{n,i}|} = \max A_{n,i}.$ 
Define
\[X_{n,i+1}(\omega)=y_j(\omega)  \mbox{ if } (j-1)/{|A_{n,i}|(\omega)} < U_{(n-1)\kappa+i+1}(\omega) < j/{| A_{n,i}|(\omega)} .\] 

Let us fix an integer $C$ , $C \geq 10K+1.$ Let us define the stopping times $\alpha_0=0$ and $\alpha_i(\omega)=\inf \{n \geq \alpha_{i-1}(\omega)+C: L_n(\omega) \leq C \},$ $i \geq 1.$ Note that by Proposition \ref{upper bound 10K}
and (3) of Theorem \ref{second theorem}, each of these stopping times is finite with probability $1.$ We now claim that if $l$ is sufficiently large and $U_{\alpha_l+k,i}(\omega)  > 2/3$ for $k=1, \ldots, C,$ $i=1, \ldots, K,$ and $U_{\alpha_l+C,\kappa}(\omega)  > 1-1/{\alpha_l(\omega)}$ then 
\begin{equation} \label{right chain of U's to get the value of L down}
L_{\alpha_l+C}(\omega)   \leq K-2. 
\end{equation}
So let us assume that $U(\omega)$ satisfies the above conditions for the first $C$ periods following period $\alpha_l(\omega)$. Let $n=\alpha_l(\omega).$ Let us consider period $n+k$ for $1 \leq k \leq C.$  Note that if at the beginning of the $i$-th addition step there is at least one gap present, (i.e., if $L_{n+k,i-1}(\omega) > 0$), then the condition $U_{n+k,i}(\omega) > 2/3,$ implies that 
a particle gets added at a site which does not belong to the outer boundary of $A_{n+k,i-1}(\omega)$ and hence $L(\omega)$ decreases in this step. Once $L(\omega)$ reaches the value $0,$ (i.e., all gaps got removed), it stays at this value when we add another particle. On the other hand, in each single deletion step $L(\omega)$ can increase at most by $1.$ Thus, considering now the whole period, we can conclude that if $L_{n+k-1}(\omega) > K-1$ then $L_{n+k}(\omega) \leq L_{n+k-1}(\omega)-1,$ and if $L_{n+k-1}(\omega) \leq K-1$ then $L_{n+k,K}(\omega) =0$ and $L_{n+k}(\omega) \leq K-1.$ Now by definition of $n=\alpha_l(\omega)$ we have $L_n(\omega) \leq C.$ Thus under the above assumptions on $U(\omega)$, $L(\omega)$ will reach a value less than $K$ within the first $C-K+1$ periods following period $n$ and so $L_{n+C-1}(\omega) \leq K-1$ and $L_{n+C,K}(\omega)=0$. Thus after the next $K-2$ deletions we have $L_{n+C,\kappa-1}(\omega) \leq K-2.$ Now note that $|A_{n+C,\kappa-1}|(\omega) = n+C+3 \leq 2n$ (if $n$ is sufficiently large). So $1-1/n \geq 1-2/{|A_{n+C,\kappa-1}|(\omega)}$ and thus under the above assumption on  $U_{n+C,\kappa}(\omega),$  the last deletion in period $n+C$ removes an endpoint of $A_{n+C,\kappa-1}(\omega).$ Hence $L(\omega)$ does not increase in this step and our claim thus follows. 

We next claim that for any sufficiently large $C$ there exists an integer $D>0$ such that 
\begin{equation} \label{alpha less than Dl eventually}
P(\alpha_l < D l \mbox{ eventually}) = 1.
\end{equation}
Assuming that this claim holds, let us choose $C \geq 10K+1$ and $D>0$ such that  (\ref{alpha less than Dl eventually}) holds. Let ${\cal E}_l=\{U_{\alpha_l+k,i} > 2/3,$  $k=1, \ldots, C,$ $i=1, \ldots, K;$ $U_{\alpha_l+C,\kappa} > 1-1/{D l} \}.$ The events ${\cal E}_l$ are independent with $P({\cal E}_l) = 3^{-CK} (D l)^{-1}.$ Hence by the Borel-Cantelli lemma
$P(\limsup {\cal E}_l) = 1.$ Using that by our choice of $C$ and $D$ we have $P(\alpha_l < D l \mbox{ eventually}) = 1,$ we can conclude that for almost every $\omega$ there exist infinitely many indices $l$ with $U_{\alpha_l+k,i}(\omega)  > 2/3$ for $k=1, \ldots, C,$ $i=1, \ldots, K$ and $U_{\alpha_l+C,\kappa}(\omega)  > 1-1/{\alpha_l(\omega)}$ and thus by (\ref{right chain of U's to get the value of L down}),  $L_{\alpha_l+C}(\omega)   \leq K-2$ for infinitely many indices $l$. Hence the proposition follows. 

\vspace{7 mm} 

So it remains to show that we can choose $C \geq 10K+1$ and $D>0$ such that  (\ref{alpha less than Dl eventually}) holds. Let $\beta_l(\omega) = \alpha_l(\omega)-(\alpha_{l-1}(\omega) +C).$ Let us first observe that if $G_{n-1}(\omega) \geq C$ then 
for $i=0, \ldots, K-1$ we have $G_{n,i}(\omega) \geq  G_{n-1}(\omega)-i > C-K$ and hence $|\partial A_{n,i}|(\omega) > C-K.$ Now, $L_n-L_{n-1} \leq \sum_{i=1}^{K} {\bf 1}_{ \{X_{n,i} \in \bar{\partial} A_{n,i-1} \} } -1$ and 
$\{X_{n,i} \in \bar{\partial} A_{n,i-1} \} = \{U_{{(n-1) \kappa +i}} \leq 2/{|\partial A_{n,i-1}|} \}.$  
Thus, if $G_{n-1}(\omega) \geq C$ then 
\[(L_n-L_{n-1})(\omega) \leq \sum_{i=1}^{K} {\bf 1}_{ \{U_{{(n-1) \kappa +i}} \leq \gamma \} } (\omega) -1, \]
where $\gamma=2 (C-K)^{-1}.$
Next note that if $\beta_l(\omega)>0$ then $L_{\alpha_{l-1}+C+m}(\omega) > C$ for $m=0, \ldots, \beta_l(\omega)-1.$ 
So, if $L_n(\omega) \leq G_n(\omega)+1$ for all $n \geq \alpha_{l-1}(\omega)+C$ then for $m=0, \ldots, \beta_{l}(\omega)$ we have 
\[
(L_{\alpha_{l-1}+C+m}-L_{\alpha_{l-1}+C})(\omega) \leq \sum_{n=\alpha_{l-1}(\omega)+C+1}^{\alpha_{l-1}(\omega)+C+m} \sum_{i=1}^{K} {\bf 1}_{ \{U_{(n-1) \kappa +i} \leq \gamma \} }(\omega) - m .
\]
Clearly, $L_{\alpha_{l-1} +C} \leq L_{\alpha_{l-1}}+C(K-1) \leq C+C(K-1)= CK$ and so the left hand side can be bounded below by $-CK.$ Thus we get the following: If 
$L_n(\omega) \leq G_n(\omega)+1$ for all $n \geq \alpha_{l-1}(\omega)$ then for $m=0, \ldots, \beta_{l}(\omega)$ we have
\begin{equation} \label{ineq: relating alpha l to the CK-decrease of a random walk}
\sum_{n=\alpha_{l-1}(\omega)+C+1}^{\alpha_{l-1}(\omega)+C+m} \sum_{i=1}^{K} {\bf 1}_{ \{U_{{(n-1) \kappa +i}} \leq \gamma \} }(\omega) - m > -CK.
\end{equation}
So in view of part (3) of Theorem \ref{second theorem} we can say that, roughly speaking, in a path-by-path way the waiting time $\beta_l(\omega)$ is eventually dominated by the time a corresponding random walk will need to decrease by the amount $CK.$ Let us now work on making this statement precise.
 
Let us assume that on our probability space $(\Omega, {\cal F},P)$ a process $U^4$ is defined, which is independent of $U$ and such that the random variables $U^4_j$ are independent and uniformly distributed on $(0,1).$
We will recursively construct a process $U^*$ and we define the random walk  $L^*$ as
\begin{eqnarray*}
L^*_0 &=& 0 \\
L^*_n &=& -n+ \sum_{k=1}^{n} \sum_{i=1}^K {\bf 1}_{ \{U^*_{k,i} \leq \gamma \} }.
\end{eqnarray*}
Let $f_n^*$ denote the end of the $n$-th CK-drop for $L^*(\omega),$ i.e., let $f^*_0(\omega)=0,$ $f_n^*(\omega)=\inf \{j>f_{n-1}^*(\omega):(L^*_j-L^*_{f_{n-1}^*})(\omega) \leq -CK \}$ and let $\beta^*_n=f_n^*-f_{n-1}^*$ denote its length.
 
For each $\omega \in \Omega$ we construct the sequence $U^*(\omega)$ by pasting together parts of $U(\omega)$ and of $U^4(\omega)$ as follows:
Starting at the end of period $\alpha_1(\omega)+C$ we select successive elements from the sequence $U(\omega)$ until we reach $\alpha_2(\omega)$ or until the first CK-drop of $L^*(\omega)$ is completed (whatever happens first). If at this time the first CK-drop of $L^*(\omega)$ has not yet been completed we continue with elements from the sequence $U^4(\omega)$ to complete this CK-drop. (In particular, if $\alpha_2(\omega)=\alpha_1(\omega)+C$ then all elements of $U^*(\omega)$ which have been selected so far, were chosen from $U^4(\omega)$). $U^*(\omega)$ thus is constructed up to its first CK-drop. We then restart the procedure at the end of period $\alpha_2(\omega)+C,$  filling in (if necessary) with elements from the sequence $U^4(\omega)$ to complete the next CK-drop for $L^*(\omega).$ We keep iterating this procedure replacing the current $\alpha_n(\omega)$ by $\alpha_{n+1}(\omega)$ at each new start.

It is not hard to show that the above construction of $U^*$ fits into the scheme described in Proposition \ref{prop:good selection rules}. Hence we can conclude that the random variables $U^*_{n}$ are independent and uniformly distributed on $(0,1).$ Thus $L^*$ is a random walk and the random variables $\beta^*_l$ are independent and identically distributed. We now claim that 
\begin{equation} \label{beta leq beta*}
P(\beta_l \leq \beta^*_l \mbox{ eventually})=1.
\end{equation}
First note that by part (3) of Theorem 2 for almost every $\omega$  there exists an index $m$ (depending on $\omega$) such that $L_n(\omega) \leq G_n(\omega)+1$ for all $n \geq \alpha_{m-1}(\omega).$ So, when in the construction of $U^*(\omega)$ we restart our iterative procedure at time $\alpha_{l-1}(\omega)+C,$ $l \geq m,$ we keep selecting from the sequence $U(\omega)$ until we reach $\alpha_l(\omega).$ (Here we use that by (\ref{ineq: relating alpha l to the CK-decrease of a random walk}) no CK-decrease gets completed in this time interval). Thus $\beta_l(\omega) \leq \beta^*_l(\omega)$ for all $l,\omega$ as above. 

Next we claim that for all sufficiently large $C$ we have 
\[ E \beta^*_1 < \infty. \]
So let us define ${\cal D}_n = \cup_{i=1}^K \{U^*_{n,i} \leq \gamma \}$ and let $W^*_n=\sum_{i=1}^{n} (-{\bf 1}_{ {\cal D}_i^c } + K  {\bf 1}_{ {\cal D}_i }).$ Then $L^*_n \leq W^*_n,$ and so, defining $R^*=\inf \{m : W^*_m = -CK \},$ we have $\beta^*_1 \leq R^*.$ So it suffices to show that for all sufficiently large $C$ we have $E R^* < \infty.$ Define $p^*=P({\cal D}_1).$ Let $Z$ be a Bernoulli random walk with $P( Z_1=1)=p,$ $p \in (0,1)$ to be chosen right after the next definition. Define the random times $\xi_n$ as $\xi_0=0,$ $\xi_n(\omega)=\inf \{m \geq \xi_{n-1}(\omega): (Z_m-Z_{\xi_{n-1}})(\omega) \in \{-1,K\} \}$ and let $\psi=\inf \{m : Z_m = -CK \}.$ Now let us choose $p$ such that $P(Z_{\xi_{1}}=K)=p^*$ and let us consider the random walk $W_n=Z_{\xi_n}.$ Let $R=\inf \{m : W_m = -CK \}.$ Note that by choice of $p$ the random variable $R$ has the same distribution as $R^*.$ Moreover, $R \leq \psi$ and $\psi$ has finite expectation if and only if $p<1/2.$ Now $p$ (which depends on $p^*$) converges to $0$ as  
$p^*$ converges to $0$ and $p^*$ (which by definition depends on $\gamma$ and hence on $C$) converges to $0$ as $C$ converges to infinity. Thus for all sufficiently large $C$ we have $p<1/2$ and so $\psi,$ and hence $R$ and $\beta^*_1,$ have finite expectation for any such choice of $C.$

We are now ready to choose $C$ and $D$ such that (\ref{alpha less than Dl eventually}) holds. Let us fix a constant $C \geq 10K+1$ such that $\beta^*_1$  has finite expectation. Let $r=E \beta^*_1.$  By the strong law of large numbers we have $P(\sum_{k=1}^{l} \beta^*_k < (r+1) l \mbox{ eventually}) = 1$ and combining this with (\ref{beta leq beta*}) we can conclude that  $P(\sum_{k=1}^{l} \beta_k < (r+2) l \mbox{ eventually}) = 1.$  Now recall that $\alpha_l=Cl+\sum_{k=1}^{l} \beta_k.$ Thus, if we choose $C$ as above and $D=C+r+2$ then (\ref{alpha less than Dl eventually}) holds. 
\newline
\vspace{7 mm}
\fine

\section{Age of the Oldest Particle}

Let us now fix $K=2.$ Recall that for $x \in A_n,$ $Y_n(x)$ denotes the number of deletions which happened since the present particle at $x$ was added and  $O_n=\max \{Y_n(x): x \in A_n \}$ denotes the age of the oldest particle of the set $A_n.$ Thus $O_1=1$ since both particles of $A_1$ survived one deletion, $O_2=2,$ since at least one particle of $A_1$ does not get deleted in period 2 and $O_3 \in \{2,3\},$ both values having positive probability since  $O_3(\omega)=2$ if and only if both particles of $A_1(\omega)$ are deleted in periods 2 and 3. Note that for all $n,$ $O_n \geq (n+1)/2$ since $A_n$ has $n+1$ particles and two particles are added in each period. Let the distribution function $F$ be defined by  $F(x)=(1-\exp(-x^2)) {\bf 1}_{(0,\infty)}.$ We will now prove Theorem \ref{Theorem: age of the oldest particle},
i.e., we will show that as $n \rightarrow \infty$, $(n-O_n)/\sqrt n$ has a limiting distribution whose distribution function is $F$. The proof will be independent of our previous work.

{\bf Proof of Theorem \ref{Theorem: age of the oldest particle} }.
For $i,j,k \in {\bf N},$ $\min(j,k) \geq i,$ let $p(i,j,k)$ denote the probability that for a fixed subset $B$ of size $i$ of the set $A_{j-1}$ all the particles at sites belonging to $B$ will get deleted within the first $k$ periods following period $j-1,$ i.e., for any subset $B$ of $A_{j-1}$ of size $i,$ $p(i,j,k)$ denotes the probability that for each $y \in B$ there exists a period $j \leq l < j+k$ with $y \not \in A_l.$
For $k \in {\bf N},$ $1 \leq i \leq k$ let $(k)_i$ denote the product $\prod_{j=1}^i (k-i+1).$ We will prove by induction that 
\begin{equation}\label{induction claim}
p(i,j,k)=\frac{(k)_i}{(j+k+1)_i} \; .
\end{equation}
Let $i=1$ and $y \in A_{j-1}.$ Note that $1-p(1,j,k)$ equals the probability that in each of the first $k$ periods following period $j-1$ a particle gets deleted at a site different from $y.$ Now $|A_{j,K}|=j+2$ and $|A_{j,K} \setminus \{y\}|=j+1$  and thus
\[ 1-p(1,j,k)=\prod_{l=1}^k \frac{j+l}{j+l+1} = \frac{j+1}{j+k+1} \; . \]
Hence 
\[ p(1,j,k) = \frac{k}{j+k+1} \; . \]
Next note that for $j \geq i$ 
\[ p(i,j,i)=\prod_{l=1}^{i} \frac{i-l+1}{j+l+1} = \frac{(i)_i}{(j+i+1)_i} \; , \]
since in order to delete $i$ particles within $i$ deletions we need to 
delete one of these in each single deletion step.
So fix $n,j,m \in {\bf N}$ with $n \geq 2,$ $j \geq n$ and $m>n$ and suppose that (\ref{induction claim}) holds for all $(i,j,k) \in {\bf N}^3$ with  $\min(j,k) \geq i$ and  ($i \leq n-1$ or  $i=n$ and $k \leq m-1$). Let $B$ be a subset of size $n$ of $A_{j-1}.$ Using that in period $j$ the probability of deleting a particle at a site belonging to $B$ equals $n/(j+2)$ we get
\[ p(n,j,m)= \frac{n}{j+2} p(n-1,j+1,m-1) + \frac{j+2-n}{j+2}p(n,j+1,m-1), \]
and thus by our induction hypothesis
\begin{eqnarray*}
 p(n,j,m) &=& \frac{n}{j+2} \; \frac{(m-1)_{n-1}} {(j+m+1)_{n-1}} + \frac{j+2-n}{j+2} \;  \frac{(m-1)_{n}} {(j+m+1)_{n}}  \\
          &=& \frac{(m-1)_{n-1}} {(j+m+1)_{n-1}} \; \frac{1}{j+2} \left( n+ \frac{(j+2-n)(m-n)}{j+m-n+2} \right) \\
          &=& \frac{(m-1)_{n-1}} {(j+m+1)_{n-1}} \; \frac{1}{j+2} \left( \frac{ j(n+m-n) + n(m-n+2) + (2-n)(m-n)}{j+m-n+2} \right) \\
          &=&  \frac{(m-1)_{n-1}} {(j+m+1)_{n-1}} \; \frac{1}{j+2} \; \frac{(j+2)m} {j+m-n+2} \\
          &=& \frac{m (m-1)_{n-1}}{(j+m+1)_{n-1}(j+m-n+2)}\\
          &=& \frac{(m)_n} {(j+m+1)_n}.
\end{eqnarray*}
Now let $x \geq 0$ and let $k(n,x)$ denote the biggest integer which is strictly smaller than $n-\sqrt n x.$ We have
\begin{eqnarray*}
P(\frac{n-O_n}{\sqrt n} >x) &=& P(O_n < n- {\sqrt n} x) \\ 
                       &=& P(O_n \leq k(n,x)) \\
                       &=& P(\mbox{ all particles of $A_{n-k(n,x)}$ will get deleted within periods} \\
                       &  & \qquad n-k(n,x)+1, \ldots, n  ),
\end{eqnarray*}
and hence, using that $k(n,x)=n-\lfloor{{\sqrt n}x}\rfloor -1$ and $n-k(n,x)+1=\lfloor{{\sqrt n}x}\rfloor +2$ we get 
\begin{eqnarray*}
P(\frac{n-O_n}{\sqrt n} >x)   
                       &=& p(\lfloor{{\sqrt n}x}\rfloor +2,\lfloor{{\sqrt n}x}\rfloor +2,n-\lfloor{{\sqrt n}x}\rfloor -1) \\
                       &=&\frac{(n-\lfloor{{\sqrt n}x}\rfloor -1)_{\lfloor{{\sqrt n}x}\rfloor +2}}{(n+2)_{\lfloor{{\sqrt n}x}\rfloor +2} }  \\
                       &=& \prod_{i=1}^{\lfloor{{\sqrt n}x}\rfloor +2} \frac{n-\lfloor{{\sqrt n}x}\rfloor -i}{n+3-i}. 
\end{eqnarray*}
For $1 \leq i \leq  \lfloor{{\sqrt n}x}\rfloor +2$ we have
\[ \frac{n-\lfloor{{\sqrt n}x}\rfloor -i}{n+3-i} \leq \frac{n-{\sqrt n}x -(i-1)}{n-(i-1)} \leq \frac{n-{\sqrt n}x}{n} =1-\frac{x}{\sqrt n}, \]
and thus 
\[ P(\frac{n-O_n}{\sqrt n} >x) \leq (1-\frac{x}{\sqrt n})^{\sqrt n x}. \]
Using that the right hand side converges to $\exp(-x^2)$ as $n \rightarrow \infty$ we get
\begin{equation} \label{upper bound}
\limsup_{n \rightarrow \infty} P(\frac{n-O_n}{\sqrt n} >x) \leq \exp{(-x^2)}.
\end{equation}
Similarly we have
\begin{eqnarray*}
\frac{n-\lfloor{{\sqrt n}x}\rfloor -i}{n+3-i} &=&  \frac{n+3-i-(3+\lfloor{{\sqrt n}x}\rfloor )}{n+3-i} \geq 1-\frac{{\sqrt n}x+3}{n-i+3} \\
&=& 1-\frac{{\sqrt n}x}{n} \left(\frac{1+3/({\sqrt n}x)}{1-(i-3)/n} \right).     
\end{eqnarray*}
Thus for any $\varepsilon > 0$ there exists $n_0 \in {\bf N}$ such that  for all $n \geq n_0,$ $1 \leq i \leq {\sqrt n} x +2$, 
\[\frac{n-\lfloor{{\sqrt n}x}\rfloor -i}{n+3-i} \geq 1-\frac{x}{\sqrt n} (1+\varepsilon). \]
Hence for all sufficiently large $n$ we have
\[ P \left(\frac{n-O_n}{\sqrt n} >x \right) \geq \left(1-\frac{x(1+\varepsilon)}{\sqrt n} \right)^{\lfloor{\sqrt n x}\rfloor +2}. \]
Since the right hand side converges to $\exp(-x^2(1+\varepsilon))$ as $n \rightarrow \infty$ (and $\varepsilon > 0$ was arbitrarily chosen) we get
\begin{equation} \label{lower bound}
\liminf_{n \rightarrow \infty} P(\frac{n-O_n}{\sqrt n} >x) \geq \exp{(-x^2)}.
\end{equation}
Combining (\ref{lower bound}) and (\ref{upper bound}) the proof of the theorem thus gets completed. 
\newline
\vspace{7 mm}
\fine

\chapter{Holes in Diffusion Limited Aggregation}

\section{Description of the Model and Statement of Results}

We consider Diffusion Limited Aggregation (DLA) on ${\bf Z}^2,$ a cluster growth model which starts with a finite connected subset $A_0$ of ${\bf Z}^2,$ usually taken as $A_0=\{0\}.$ (A subset $b$ of ${\bf Z}^2$ is connected if any two points of $b$ can be connected by a nearest-neighbor lattice path all the points of which belong to $b$). Then, one by one, particles get released ``at infinity" and perform nearest-neighbor symmetric random walks until they attach at the first site visited which is adjacent to the cluster. More formally, let $A_n$ denote the cluster formed after $n$ steps and let $\partial A_n,$ the boundary of $A_n$, be the set of all lattice points which are adjacent to $A_n$ but not in $A_n.$ Then for any boundary site $x$ of $A_n$ we have 
\[ P(A_{n+1}=A_n \cup \{x\} \mid A_n)=\mu_{\partial A_n}(x), \]
where $\mu_{\partial A_n}(x)$ denotes harmonic measure on the boundary of $A_n$ (see \cite[Sec. 2.6] {Lawler:Intersections}).
We are interested in the
question whether the formed cluster will have ``holes'', i.e., whether its
complement will have maximal connected subsets which are finite. We will
show that with probability $1,$ $H(A_n)$, the number of
holes of the cluster $A_n$, tends to infinity as $n$ tends to infinity. 
This result will be derived as a corollary to the following

\begin{theorem} \label{main theorem} 
There exists a constant $c>0$ such that for any finite
connected subset $a$ of ${\bf Z}^2$ there exists an integer $L$ with
\[P(H(A_L) \geq H(A_0)+1 \mid A_0=a) \geq c .\]
\end{theorem}

\begin{corollary} \label{cor to main theorem} 
For any finite connected subset $a$ of ${\bf Z}^2$ we have
\[P(H(A_n) \mbox{ converges to infinity as } n \rightarrow \infty \mid
A_0=a) = 1. \]   
\end{corollary}

\section{Outline}

Let us first introduce some notation. Let ${\bf N}_0={\bf N} \cup \{0\}.$ For $x=(x_1,x_2) \in {\bf Z}^2$ define $|x|=|x_1|+|x_2|$ and for $x,y \in {\bf Z}^2$ define the lattice distance between $x$ and $y$ as $|x-y|.$ If $b$ is a nonempty subset of ${\bf Z}^2$ and $x \in {\bf Z}^2,$ define the lattice distance between $x$ and $b$ as $d(x,b)=\min \{|x-y|: y \in b\}.$ Let $S(b)=\sup\{|x|: x \in b \},$ and if $b$ is finite and $L$ is a positive integer, let $\Delta(b,L)=\{x \in {\bf Z}^2 : |x|=S(b)+L \}.$

A path in ${\bf Z}^2$ is a finite or infinite sequence $\rho$ of points of ${\bf Z}^2$ such that any two successive points $\rho_i,$ $\rho_{i+1}$ of the sequence are adjacent in ${\bf Z}^2.$ Moreover, the first point of the sequence will always be indexed by $0,$ i.e., if $\rho$ is finite then $\rho=(\rho_0, \ldots , \rho_l)$ for some $l \in {\bf N}_0$ and if $\rho$ is infinite then $\rho=(\rho_i)_{i \in {\bf N}_0}.$ If $\rho$ is finite and $l$ is chosen as above let us define $|\rho|,$ the length of $\rho,$ as $|\rho|=l$ and if $\rho$ is infinite, let $|\rho|=\infty.$ For $b \subset \Ztwo$ let $\tau(b,\rho)=\inf\{0 \leq i < |\rho|+1: \rho(i) \in
b \},$ ($\inf \emptyset = \infty$), denote the hitting time of the set $b$
for the (finite or infinite) path $\rho$ and for $x \in \Ztwo$ let $\tau(x,\rho)=\tau(\{x\},\rho);$ if $\tau(b,\rho) < \infty$ let  $\rho_{\tau(b)}=\rho_{\tau(b,\rho)}$ and similarly define $\rho_{\tau(x)}$ for $x \in \Ztwo$.

A component of $b,$ ($b \subset \Ztwo$), is a maximal
connected subset of $b$ and a hole of $b$ is a finite component of $\Ztwo
\setminus b.$ The number of holes of $b$ is denoted by $H(b).$ 

\vspace{7 mm}

Let us now consider DLA with initial set $A_0=a,$ where $a$ is a finite connected subset of ${\bf Z}^2.$ Let $L$ be a positive integer. We want to investigate the growth of the cluster up to time $L.$ Clearly, the first $L$ particles will hit $\Delta=\Delta(a,L)$ before (or upon) hitting the boundary of the growing cluster. Thus, by the strong
Markov property of simple random walk, instead of releasing particles ``at infinity'' we may as well have them start at $\Delta$
according to harmonic measure $\mu_{\Delta}$ on $\Delta.$ 

So let us work on the canonical probability space associated with $L$
independent random walks the starting positions of which are distributed
according to $\bigotimes_{i=1}^L \mu_{\Delta}.$
Let $ \Omega = \prod_{i=1}^L {\cal S},$ where ${\cal S}=\{\rho: \rho \mbox{ is an infinite path in
$\Ztwo$ with $\rho_0 \in \Delta$} \}.$  Let $\Phi_i: \Omega \rightarrow {\cal S}$ denote the projection onto
the $i$-th coordinate, i.e., $\Phi_i(\omega)=\omega_i,$ and define
$\Phi_{ij}: \Omega \rightarrow \Ztwo$ by $\Phi_{ij}(\omega)=\omega_{ij}$
and $\Phi_i^k: \Omega \rightarrow (\Ztwo)^{k+1}$ by
$\Phi_i^k(\omega)=(\omega_{ij})_{0 \leq j \leq k}.$
Let 
$\cal F$ be the $\sigma$-field in $\Omega$ generated by the maps
$\Phi_{ij},$ $1 \leq i \leq L,$ $j \in \Nzero$ and 
define the probability measure $P$ on $(\Omega, \cal F)$ by
\[P(\Phi_i^{l_i}=\rho_i, i=1, \ldots ,L) = \prod_{i=1}^L
\left(\mu_{\Delta}(\rho_{i0}) 4^{-l_i}\right) \]
for any finite paths $\rho_1, \ldots ,\rho_L$ with $\rho_{i0} \in \Delta$
and $|\rho_i|=l_i,$ $i=1, \ldots ,L.$ 
For any paths $\rho_1, \ldots, \rho_L$ in $\Ztwo$ 
with $\rho_{i0} \not\in a \cup {\partial a}$ 
define the clusters
${\cal C}(\rho_1, \ldots, \rho_n),$ $1 \leq n \leq L,$ recursively as follows: If $\tau(a,\rho_1) < \infty$ let ${\cal C}(\rho_1)= a \cup \{ \rho_{1 \tau(\partial a)} \}  $ and if $\tau(a,\rho_1) = \infty$ let ${\cal C}(\rho_1)={\tilde {\cal C}}_1,$ where ${\tilde {\cal C}}_1$ is a fixed cluster of size $|a|+1$ with $a \subset {\tilde {\cal C}}_1.$ Thus, if $\rho_1$ hits the boundary of the set $a,$ then ${\cal C}(\rho_1)$ is the set obtained from $a$ by adding the point where $\rho_1$ hits $\partial a.$
Now suppose that ${\cal C}(\rho_1, \ldots, \rho_n)$ has already been defined, $n < L$. If $\tau(a,\rho_i) < \infty$ for $i=1 \ldots n+1$ let $ {\cal C}(\rho_1, \ldots \rho_{n+1}) ={\cal C}(\rho_1, \ldots \rho_{n})
\cup \{ \rho_{n+1, \tau(\partial {\cal C}(\rho_1, \ldots \rho_{n}))} \}$ and if $\tau(a,\rho_i) = \infty$ for some index $1 \leq i \leq n+1,$ let ${\cal C}_{n+1}(\omega)={\tilde {\cal C}}_{n+1},$ where ${\tilde {\cal C}}_{n+1}$ is a fixed cluster of
size $|a|+n+1$ with $a \subset {\tilde {\cal C}}_{n+1}.$ In other words, if the previous paths $\rho_1 \ldots \rho_n$ hit $\partial a$ and if the next following path $\rho_{n+1}$ does so as well, then the new cluster ${\cal C}(\rho_1, \ldots, \rho_{n+1})$ is obtained from the previously constructed cluster ${\cal C}(\rho_1, \ldots, \rho_n)$ by adding the point where $\rho_{n+1}$ hits the boundary of ${\cal C}(\rho_1, \ldots, \rho_n).$ Define
 \newline
$\bullet \quad
A_0=a;$ \newline
$\bullet \quad A_n(\omega)  =  {\cal C}(\omega_1, \ldots ,\omega_n),$ $\; 
n=1, \ldots, L.$ \newline
Then $(A_n)_{n=0 \ldots L}$ is a Markov Chain on
$(\Omega,{\cal F},P)$ the transition probabilities of which are exactly
those of the DLA-model. (Here we use that random walk in $\Ztwo$
is recurrent and thus each random walk will hit $a$ with probability 1). 
If we want to emphasize the dependence of $(\Omega, \cal F, P)$ and
$A_n$ on the starting cluster $a$ and the integer $L,$
we write $(\Omega^{a,L}, {\cal F}^{a,L}, P^{a,L})$ and $A_n^{a,L}.$

\vspace{7 mm}

Using this particular construction, the hole problem will be solved by showing that for any sufficiently large finite cluster $a$ we can choose an integer $L$
depending on $a$ and we can define a map $\varphi:
\Omega^{a,L} \rightarrow \Omega^{a,L}$ satisfying

\begin{equation} \label{first property of varphi}
H(A_L^{a,L}(\varphi(\omega)) \geq H(a) + 1 
\end{equation}
for all $\omega \in \Omega^{a,L}$ with $\tau(a,\omega_i) < \infty,$ 
$i=1, \ldots, L,$  \newline
and
\begin{equation} \label{second property of varphi}
P^{a,L}(\varphi(\Omega^{a,L})) \geq c,
\end{equation}
where $c$ is a fixed positive constant independent of $a$ and $L.$

Note that by (\ref{first property of varphi}), 
the map $\varphi$
has to be constructed such that the changed paths $\varphi(\omega)_1,
\ldots ,
\varphi(\omega)_L$ will build a new cluster $A_L(\varphi(\omega))$ which has at least one more hole
than the cluster $a$ we start with. 
We will achieve this as follows: First
we surround the initial cluster $a$ by a simple closed curve $\Gamma$ such
that each lattice point on $\Gamma$ has distance 40 from $a.$ We will then
divide the region ``between'' $\Gamma$ and the cluster $a$ into (connected)
``patches'' such that the ``border line'' between any two adjacent patches
corresponds to a lattice path (of length 40) connecting $a$ to $\Gamma.$ 
We will then wait until one of these patches gets hit by a certain minimum number of
random walks strictly before these walks add onto the growing cluster. Those random walks which are the first to hit this ``lucky''
patch will get changed by inserting a new loop. Since we also control the ``chain reaction'' by changing a few additional paths (due to these additional changes the evolution of $A_n(\omega)$ and of $A_n(\varphi(\omega))$ will be very close to each other until those random walks hit the lucky patch), these new loops can be chosen such that the particles add onto the cluster at lattice points ``well in the middle'' of the ``lucky'' patch. These new positions will be such that a lattice point which belongs to the unbounded component of $\Ztwo \setminus a$ will get disconnected from the ``outside''. Thus (\ref{first property of varphi}) will be satisfied.

In order to obtain 
(\ref{second property of varphi}),
 the map $\varphi$ will satisfy the following: It will suffice to know each of the paths $\omega_1, \ldots , \omega_L$ up to the hitting time of $a$ in order to decide which paths will get changed and what these changes are. All new loops will be chosen such that their length is bounded above by a fixed constant and they will be inserted before the paths to be changed hit $a.$ Moreover, the total number of new loops will be bounded above by a fixed constant and $\varphi$ will be one-to-one.

Properties (\ref{first property of varphi}) and (\ref{second property of varphi})  of the map $\varphi$ then imply Theorem \ref{main theorem} as follows. 
By  (\ref{first property of varphi}) we can conclude that for any $\omega \in \Omega^{a,L}$ satisfying $\omega=\varphi(\tilde {\omega})$ for some $\tilde {\omega} \in \Omega^{a,L}$ we have $H(A_L^{a,L}(\omega))=H(A_L^{a,L}(\varphi(\tilde{\omega}))) \geq H(a)+1$ and hence
\[ \{H(A_L^{a,L}) \geq H(a)+1\} \supset \{\varphi(\Omega^{a,L})\} . \]
Combining this with (\ref{second property of varphi})
 we get 
\begin{eqnarray*}
P(H(A_L) \geq H(A_0)+1 \mid A_0=a) &=& P^{a,L}(H(A_L^{a,L}) \geq H(a)+1)\\
&\geq& P^{a,L}(\varphi(\Omega^{a,L})) \\
&\geq& c
\end{eqnarray*}
(for all sufficiently large finite clusters $a$). Thus Theorem \ref{main theorem} now follows by the Markov property of $A.$

The construction of $\varphi$ and the verification of properties (\ref{first property of varphi}) and (\ref{second property of varphi}) is our remaining task.
We will formally define the map $\varphi$ in section 3.6 after introducing most construction parts (patches, "lucky patch", new positions of attaching etc.) in sections 
3.3 - 3.5. In section 3.6 we will also prove that $\varphi$ satisfies  (\ref{first property of varphi}) and in section~3.7 we will verify that (\ref{second property of varphi}) holds, thus completing the proof of the theorem.  

\section{Construction of Patches}

Let $a$ be a finite cluster. For any two points $u,v$ in $\Ztwo$ let $u
\la v$
denote the line segment with endpoints $u,v.$ Points $z \in u \la v$ with $z
\not \in \{u,v\}$ are called inner points of $u \longleftrightarrow v.$ A
line segment $u \longleftrightarrow v$
with $u,v \in a$ and $|u-v|=1$ is called an edge of~$a.$ Let $G(a)$ denote the
subset of ${\bf R}^2$ obtained by taking the union over all edges of $a.$
The construction of patches will
be done in the following steps: \newline
 \newline
i) We surround $a$ by an oriented curve $\Gamma,$ which can be obtained by
joining together line segments of the form $u \longleftrightarrow u \pm
e_1,$ $u \la u \pm e_2$ and 
$u \la u \pm e_1 \pm e_2,$ $u \in \Ztwo,$ such that $a$
lies to the right of $\Gamma$ and  each lattice point on $\Gamma$ has
lattice distance 40 from $a.$\newline
 \newline
ii) For each lattice point $\xi$ on $\Gamma$ we choose a lattice path
$\gamma(\xi)$ of length 40 from $\xi$ to $a.$ Moreover, we choose these paths
such that whenever two paths touch they stay together.\newline
 \newline
iii) We select an integer $I$ (depending on $a$) and lattice points $u_0, u_1, \ldots ,u_I=u_0$ on $\Gamma$ (where
the enumeration corresponds to the order in which these points are hit by
$\Gamma$) such that $\gamma(u_0), \gamma(u_1), \ldots ,\gamma(u_{I-1})$ are
nonintersecting and, vaguely speaking, neither too close nor too far from
each other. (In this part as well as in (iv) we assume that the cluster $a$ is sufficiently big). \newline
 \newline
iv) Letting $\Gamma(u_{i-1},u_i)$ denote the part of $\Gamma$ which connects
$u_{i-1}$ to $u_i$, we define ${\cal D}_i$ to be the closure of the open region
bounded by $\gamma(u_{i-1}),$ $\Gamma(u_{i-1},u_i),$ 
$\gamma(u_{i})$ and those
edges of $a$ whose interior points can be connected by a curve in ${\bf
R}^2$ to $\Gamma(u_{i-1},u_i)$
without hitting any other point of $\Gamma(u_{i-1},u_i),$
$\gamma(u_{i-1}),$
$\gamma(u_{i})$ or any other edge of $a.$ 

\vspace{7 mm} 
 
Let us now start with step (i). We will construct a polygon whose contour $\Gamma$ (after
suitable orientation) has the desired properties.
Enumerate the points of $a$ as $y_1, y_2, \ldots ,y_{|a|}$ such that for $i=1,
\ldots, |a|,$ the set of points $\{y_1, y_2, \ldots ,y_i \}$ is a connected subset of $\Ztwo.$
Define $\mbox{Pol}(y_i)$ as the closed polygon with
corners $y_i \pm 40 e_1$ and $y_i \pm 40 e_2$
and let its contour (which forms a simple closed curve in ${\bf R}^2$) be denoted by $\mbox{Con}(y_i).$ 
Now suppose that the closed polygon $\mbox{Pol}(y_1, \ldots
,y_{i-1})$ has already been constructed and it satisfies

\renewcommand{\labelenumi}{\theenumi)} 
\begin{enumerate}
\item the contour $\mbox{Con}(y_1, \ldots ,y_{i-1})$ of    $\mbox{Pol}(y_1, \ldots ,y_{i-1})$ forms a simple closed       curve which can be obtained by joining line segments of    the form described in (i);  
\item $\mbox{Pol}(y_j) \subset \mbox{Pol}(y_1, \ldots ,y_{i-      1})$ for $j=1, \ldots, i-1;$ 
\item any lattice point on $\mbox{Con}(y_1, \ldots ,y_{i-               1})$ belongs
    to $\mbox{Con}(y_j)$ for some index $1 \leq j \leq i-1.$ 
\end{enumerate}

Let $\mbox{Pol}^*(y_1, \ldots ,y_{i})$ denote the union of $\mbox{Pol}(y_1, \ldots
,y_{i-1})$ and $\mbox{Pol}(y_i)$ and let the outer contour of $\mbox{Pol}^*(y_1, \ldots ,y_{i})$ be denoted by $\mbox{Con}^*(y_1, \ldots ,y_{i}).$
Thus, if $U$ denotes the open unbounded component of ${\bf R}^2
\setminus  \mbox{Pol}^*(y_1, \ldots ,y_{i}),$ then $x \in \mbox{Con}^*(y_1, \ldots
,y_{i})$ if and only if $x \in \bar{U} \setminus U.$ Since the interiors of 
$\mbox{Pol}(y_1, \ldots ,y_{i-1})$ and $\mbox{Pol}(y_i)$ intersect in a nonempty set
(here we use that $y_i$ is a neighbor in $\Ztwo$ of one of the points $y_1,
\ldots ,y_{i-1},$ and property (2) above),
$\mbox{Con}^*(y_1, \ldots ,y_{i})$ forms a simple closed curve which surrounds a
polygon containing $\mbox{Pol}^*(y_1, \ldots ,y_{i})$ as a (possibly proper) subset.
Note that $\mbox{Con}^*(y_1, \ldots ,y_{i})$ may have corners which are not in
$\Ztwo,$ these corners being of the form $u+1/2(e_1+e_2),$ $u \in \Ztwo.$ By property (1) such a corner can only be obtained if a line
segment $u \la u+e_1+e_2$ of $\mbox{Con}(y_1, \ldots ,y_{i-1}),$ $u \in {\bf Z}^2,$
crosses the line segment $u+e_1 \la u+e_2$ of $\mbox{Con}(y_{i})$ or vice
versa.

Let us assume that $\mbox{Con}^*(y_1, \ldots ,y_{i}),$  $\mbox{Con}(y_1, \ldots
,y_{i-1})$ and $\mbox{Con}(y_j)$ are  oriented such that the corresponding polygon
lies to the right and let us first consider the special case that the directed line segment 
$u \longrightarrow u+e_1+e_2$ of $\mbox{Con}(y_1, \ldots ,y_{i-1})$ crosses the
directed line segment $u+e_1 \longrightarrow u+e_2$ of $\mbox{Con}(y_{i}).$ Then,
since $\mbox{Pol}(y_1, \ldots ,y_{i-1})$ lies to the right of $\mbox{Con}(y_1, \ldots
,y_{i-1})$ and $\mbox{Pol}(y_{i})$ lies to the right of $\mbox{Con}(y_i),$  $\mbox{Con}^*(y_1, \ldots ,y_{i})$ uses the directed
line segments $u \longrightarrow u+1/2(e_1+e_2)$ and $u+1/2(e_1+e_2)
\longrightarrow u+e_2.$ In fact, whatever orientation we choose for the two
crossing line segments, $\mbox{Con}^*(y_1, \ldots ,y_{i})$ will form a corner at $u+ 
1/2(e_1+e_2)$ by using one of the following four pairs of directed line segments 

\begin{tabbing}
yes \= \kill
\> $\cdot$ $u \longrightarrow u+1/2(e_1+e_2)$ and        $u+1/2(e_1+e_2) \longrightarrow u+e_2,$ \\
\> $\cdot$ $u+e_2 \longrightarrow u+1/2(e_1+e_2)$ and     $u+1/2(e_1+e_2) \longrightarrow u+e_1+e_2,$ \\
\> $\cdot$ $u+e_1 \longrightarrow u+1/2(e_1+e_2)$ and $u+1/2(e_1+e_2) \longrightarrow u, $ \\
\> $\cdot$ $u+e_1+e_2 \longrightarrow u+1/2(e_1+e_2)$ and $u+1/2(e_1+e_2) \longrightarrow u+e_1.$
\end{tabbing}

Note that if in any of these cases we replace the two line segments by the
single line segment joining the lattice point on $\mbox{Con}^*(y_1, \ldots
,y_{i})$ immediately preceding the corner $u+1/2(e_1+e_2)$ to the lattice point
immediately following the corner, we obtain a closed curve which surrounds
an even bigger polygon.
 
Thus let $\mbox{Con}(y_1, \ldots ,y_{i})$ denote the closed curve obtained in this
way by removing all corners which are not in $\Ztwo$ and let $\mbox{Pol}(y_1,
\ldots ,y_{i})$ denote the closed polygon with contour $\mbox{Con}(y_1, \ldots
,y_{i}).$ By construction $\mbox{Pol}(y_1, \ldots ,y_{i})$ satisfies
properties (1) - (3). Finally define $\Gamma$ as the contour of $\mbox{Pol}(y_1, \ldots ,y_{|a|})$ and let $\Gamma$ be oriented such that $\mbox{Pol}(y_1, \ldots ,y_{|a|})$ (and hence $a$) lies to the right of $\Gamma.$

Note that the special form of $\Gamma$ (as described in (i)) implies
in particular that for $L$ sufficiently large any path in $\Ztwo$ from $\Delta(a,L)$ to $a$ hits the curve $\Gamma$
in a point belonging to $\Ztwo.$

\vspace{7 mm}

{\bf Notation.} Let $\rho$ be a path in $\Ztwo.$
For any index $i$ with $i<|\rho|+1$ define the paths $\rho^i=(\rho_j)_{0
\leq j \leq i}$ and $\rho^{(\alpha,i)}=(\rho_j)_{i \leq j < |\rho|+1}.$ (Note
that $\rho^i$ ends where $\rho^{(\alpha,i)}$ starts and by joining the two parts
$\rho^i$ and $\rho^{(\alpha,i)}$ of $\rho$ we get back the path $\rho$).   
For any
lattice point $z$ with $\tau(z,\rho)<\infty$ define
$\rho^z=\rho^{\tau(z,\rho)}$ and $\rho^{(\alpha,z)}= \rho^{(\alpha,\tau(z,\rho))}.$ (Thus $\rho^z$ is obtained by stopping $\rho$ at the hitting time of $z$ and $\rho^{(\alpha,z)}$ is obtained by restarting at the hitting time of $z$).  
If $y$ is a neighbor of $\rho_0$, let $(y,\rho)=(y,\rho_0, \ldots, \rho_{|\rho|})$. Similarly, if $\rho$ is finite and $y$ is a neighbor of $\rho_{|\rho|}$, let $(\rho,y)=(\rho_0, \ldots, \rho_{|\rho|},y)$.  
For lattice paths $\rho,$ $\pi$ with $|\rho|<\infty$ and
$\rho_{|\rho|}=\pi_0$ let $(\rho,\pi)$ denote the path obtained by
appending the path $\pi$ to the path $\rho$, i.e., 
$(\rho,\pi)_i=\rho_i$, if $i \leq |\rho|$ and $(\rho,\pi)_i=\pi_{i-|\rho|}$ if $|\rho| < i < |\rho|+|\pi|+1$. If the path $\pi$ starts at a neighbor of $\rho_{|\rho|}$, (still assuming that $\rho$ is finite), we define $(\rho,\pi)=((\rho,\pi_0),\pi).$
For any finite path $\rho$ let   
$R(\rho)$ denote the corresponding
reversed  path, i.e., $R(\rho)=(\rho_{|\rho|-j})_{0 \leq j \leq |\rho|}.$

\vspace{7 mm}

We are now ready to work on step (ii). Fix a point $z$ in $a$ with maximal second component and let $\xi_0=z+40e_2.$ Then $\xi_0 \in \Gamma.$ Starting at $\xi_0$ and moving along $\Gamma$ let us enumerate the lattice points on $\Gamma$ as $\xi_0, \xi_1, \ldots , \xi_l.$ Since each lattice point on $\Gamma$ has lattice distance 40 from $a$ we can select for each point $\xi_i$ a lattice path ${\tilde \gamma}(\xi_i)$
of length 40 connecting $\xi_i$ to $a.$ By choice of $z$ and $\xi_0,$ the point $z$ is the only point of $a$ which has lattice distance 40 to $\xi_0$ and thus ${\tilde \gamma}(\xi_0)=(\xi_0,\xi_0-e_2, \ldots , \xi_0 - 40 e_2).$ 

Let us now
recursively construct the paths $\gamma(\xi_i).$ 
Let $\gamma(\xi_0)={\tilde \gamma}(\xi_0)$ and suppose that $\gamma(\xi_0),
\ldots, \gamma(\xi_i)$ have already been constructed and satisfy the 
properties in (ii). If
${\tilde \gamma}(\xi_{i+1})$ does not intersect with any of the paths $\gamma(\xi_0), \ldots, \gamma(\xi_i),$ let
$\gamma(\xi_{i+1}) = {\tilde \gamma}(\xi_{i+1}).$ Otherwise let $w$ be the
first point of intersection. Note that $w$ belongs to $\gamma(\xi_i)$ or to $\gamma(\xi_0)$ since by induction hypothesis none of the paths $\gamma(\xi_1), \ldots, \gamma(\xi_i)$ crosses $\gamma(\xi_i)$ and none of the paths $\gamma(\xi_0), \ldots, \gamma(\xi_{i-1})$ crosses $\gamma(\xi_0).$  
Recall that $\xi_0, \ldots, \xi_{i+1}$ have lattice distance 40 from $a.$ Thus ${\tilde \gamma}(\xi_{i+1})$ is a shortest connection from $\xi_{i+1}$ to $a$ and $\gamma(\xi_i)$ (resp. $\gamma(\xi_0)$) is a shortest connection from $\xi_i$ (resp. $\xi_0$) to $a.$ Since $w$ belongs to ${\tilde \gamma}(\xi_{i+1})$ and also to $\gamma(\xi_i)$ (resp. $\gamma(\xi_0)$) we can conclude that ${\tilde \gamma}(\xi_{i+1})^{(\alpha,u)}$ as well as $\gamma(\xi_i)^{(\alpha,u)}$ (resp. $\gamma(\xi_0 )^{(\alpha,u)}$) are shortest connections from $w$ to $a$ and so these two paths have the same length. Defining $\gamma(\xi_{i+1}) = ({\tilde
\gamma}(\xi_{i+1})^w, \gamma(\xi_i)^{(\alpha,w)})$, if $w \in \gamma(\xi_i),$ and $\gamma(\xi_{i+1}) = ({\tilde
\gamma}(\xi_{i+1})^w, \gamma(\xi_0)^{(\alpha,w)}),$ if $w \in \gamma(\xi_0),$ we obtain a lattice path of length 40 from $\xi_{i+1}$ to $a$ which satisfies all the requirements. Thus we have now constructed the paths $\gamma(\xi),$ $\xi \in  \Gamma \cap \Ztwo$. 

\vspace{7 mm}

Let us now proceed to steps (iii) and (iv). Define $\xi_{l+1}=\xi_0.$ (This definition corresponds to the fact that, when moving along $\Gamma$ in the prescribed direction, the point $\xi_0$ is the next lattice point encountered after $\xi_l$). For any $\xi \in \Gamma$ let $y(\xi)$ denote the endpoint of $\gamma(\xi).$ We will next define an integer $I \leq l$ and we will define two subsequences $(u_i)_{0 \leq i \leq I}$ and $(v_i)_{1 \leq i \leq I}$ of
$(\xi_i)_{0 \leq i \leq l+1}$ as follows:
Let $u_0=\xi_0.$ Now suppose that $u_0, \ldots, u_i$ have already been
constructed, $u_i=\xi_{n(i)}$ for some index $n(i) \leq l.$ 
If there exists
a point $v$ in $(\xi_{n(i)+1}, \ldots ,\xi_l )$ such that any path from
$y(v)$ to $\gamma(u_i),$ which intersects $a$ at most in its two endpoints,
has length at least 40, let $v_{i+1}$ be the first such point in the sequence 
$(\xi_{n(i)+1}, \ldots ,\xi_l ).$ (Note that for any sufficiently large set $a$ there exists a point $v_1 \in \Gamma$ which has this property with respect to $\gamma(u_0)$~)~. 
If no such point exists, we decide to drop the previously
selected point $u_i$ from our sequence and we define $I=i$ and $u_I=\xi_{l+1},$ thus
terminating the procedure. Nevertheless, for future reference, let us
denote this previously selected point by $u^I.$ 

Now suppose that $v_{i+1}$ could be
defined, say $v_{i+1}=\xi_{m(i+1)}.$ If there exists
a point $u$ in $(\xi_{m(i+1)+1}, \ldots ,\xi_l )$ such that any path from
$y(v_{i+1})$ to $\gamma(u),$ which intersects $a$ at most in its two
endpoints, has length at least 40, let $u_{i+1}$ be the first such point in
the sequence $(\xi_{m(i+1)+1}, \ldots ,\xi_l ).$ (Once again let us point out that for any sufficiently large set $a$ there exists a point $u$ in $\{\xi_{m(i+1)+1}, \ldots ,\xi_l\}$  satisfying this property with respect to $y(v_1)$). Otherwise we decide
to drop the previously selected points $u_{i}$ and $v_{i+1}$ from our
sequences and - as above - we define $I=i$ and $u_I=\xi_{l+1},$ thus
terminating the procedure. For future reference the point $u^I$ gets defined as above. 

If $I \geq 2$, let us define the patch ${\cal D}_i,$ $1 \leq i \leq I,$ as the closure of the open
connected component of ${\bf R}^2 \setminus \left (\gamma(u_{i-1}) \cup
\Gamma(u_{i-1},u_i) \cup \gamma(u_i) \cup G(a) \right )$ which contains
the point $\gamma(v_i)_1.$ (Here we obviously identify a lattice path with
the subset of ${\bf R}^2$ obtained by taking the union over all line
segments connecting two successive lattice points and we identify $\Gamma(u_{i-1},u_i)$ with the corresponding subset of ${\bf R}^2$).

\begin{remark} \label{rem:size} There exists an integer $N_0$ such that for any finite cluster $a$ with $|a| \geq N_0$ the above construction yields sequences $(u_i)_{0 \leq i \leq I}$ and $(v_i)_{1 \leq i \leq I}$ with $I \geq 2.$ 
Thus for any such cluster at least two patches get defined.
\end{remark}

{\bf Notation.} For a patch ${\cal D}={\cal D}_i$ let $\xi^l({\cal D})=u_{i-1},$
$\xi^m({\cal D})=v_{i},$ $\xi^r({\cal D})=u_{i}$ and let $\gamma^-({\cal
D})=\gamma(u_{i-1}),$ $\gamma^*({\cal D})=\gamma(v_{i})$ and 
$\gamma^+({\cal D})=\gamma(u_{i}).$ Let $y^*({\cal D})$ denote the endpoint
of $\gamma^*({\cal D}).$ 
Let ${\cal D}^*$ denote the set of all lattice points which belong to the
interior of ${\cal D}.$

\section{Properties of Patches and Definition of the "Lucky Patch"}

The following propositions will show that the patches constructed are neither too narrow nor too wide for our purpose. Let us fix an integer $N_0$ satisfying Remark \ref{rem:size} and let the construction in the previous section be done for a finite cluster $a$ with $|a| \geq N_0.$

\begin{proposition} \label{first prop of patches} 
There exists a constant $C$ (independent of the cluster $a$) such that for any patch ${\cal D}$ the number of lattice points on $\Gamma(\xi^l({\cal D}),\xi^r({\cal D}))$ is bounded above by $C.$
\end{proposition}
{\bf Proof.} Fix a patch ${\cal D}$ and let $\xi^l=\xi^l(\pd),$
$\xi^m=\xi^m(\pd)$ and $\xi^r=\xi^r(\pd).$ We will show that any vertex
$\xi$ on $\Gamma(\xi^l,\xi^m)$ satisfies $|\xi-\xi^l| \leq 121$ and any
vertex $\xi$ on $\Gamma(\xi^m,\xi^r)$ satisfies $|\xi-\xi^m| \leq 363.$
Thus, letting $C_1$ (resp. $C_2$) denote the number of lattice points
$\psi$ satisfying $|\psi| \leq 121$ (resp. $|\psi| \leq 363$), it follows
(since $\Gamma$ is nonselfintersecting) that $|\Gamma(\xi^l,\xi^m) \cap
\Ztwo| \leq C_1$ and $|\Gamma(\xi^m,\xi^r) \cap \Ztwo| \leq C_2.$ Hence
$|\Gamma(\xi^l,\xi^r) \cap \Ztwo| \leq C$,where $C=C_1+C_2.$

So fix $\xi$ on $\Gamma(\xi^l,\xi^m),$ $\xi \neq \xi^m.$ Since $\Gamma$ hits $\xi$ before hitting $\xi^m$ there exists a path $\rho$ from $y(\xi)$ to $\gamma^-(\pd)$ with $|\rho| <40.$ Let $u \in \gamma^-(\pd)$ be the endpoint of $\rho.$ Then $(\gamma(\xi),\rho,R(\gamma^-(\pd))^{(\alpha,u)})$ is a lattice path of length at most $40+39+40=119$ connecting $\xi$ and $\xi^l$ and thus $|\xi-\xi^l| \leq 119.$ Moreover, letting $\psi$ denote the lattice point on $\Gamma$ immediately preceding $\xi^m,$  we get $|\xi^m-\xi^l| \leq |\xi^m-\psi| + |\psi-\xi^l| \leq 119+2=121.$ Thus the first part of our claim follows. 

Similarly, if $\pd=\pd_i$ with $i<I,$ it follows that $|\xi-\xi^m| \leq 121$ for any lattice point $\xi$ on $\Gamma(\xi^m,\xi^r) $ and if $\pd=\pd_I,$ we get 
$|\xi-\xi^m| \leq 121$ for any lattice point $\xi$ on $\Gamma(\xi^m(\pd_I),u^I). $ (Recall that $u^I$ got defined when working on step (iii) of the construction of patches). The same argument also
shows that $|\xi-u^I| \leq 121+121 =242$ for any lattice point $\xi$ on $\Gamma(u^I,\xi^r(\pd_I)). $ Thus $|\xi-\xi^m| \leq 121+242=363$ for any lattice point $\xi$ on $\Gamma(\xi^m(\pd_I),\xi^r(\pd_I))$ and our proof gets completed. \newline
\vspace{0 mm}
\fine

\begin{proposition} \label{second prop of patches}
Let $\pd$ be a patch of $a$ and let $\rho=(\rho_0, \ldots ,\rho_l)$ be a lattice path satisfying $\rho_0 \in a,$ $\rho_i \in \gamma^-(\pd)$ or $\rho_i \in \gamma^+(\pd)$ for some $0 \leq i < l,$ $\rho_l \in \gamma^*(\pd)$ and $\rho_j \not \in a$ for $j=1, \ldots, l-1.$ Then $l \geq 20.$   
\end{proposition}
{\bf Proof.} Let $\gamma^*= \gamma^*(\pd),$ let $y \in a$ denote the endpoint of $\gamma^*,$ let $u=\rho_l$ and let $r$ denote the length of $(\gamma^*)^{(\alpha,u)}.$
Fix an index $i$ with $\rho_i \in \gamma^-(\pd)$ or $\rho_i \in \gamma^+(\pd).$ If $l-i \geq 20$ the result follows immediately. 
So assume that $l-i < 20.$ Note that $(\rho^{(\alpha,i)},(\gamma^*)^{(\alpha,u)})$ is
a path from $\gamma^-(\pd)$ or $\gamma^+(\pd)$ to $y$ which intersects
$a$ at most in its two endpoints. Thus, by definition of $y$, its length is
at least 40. Hence $|\rho^{(\alpha,i)}|+|(\gamma^*)^{(\alpha,u)}| = l-i+r \geq 40$ or
$r \geq 40 - (l-i) \geq 20.$ Now note that both $R(\rho)$ and
$(\gamma^*)^{(\alpha,u)}$ connect $u$ to $a$. Thus, since $(\gamma^*)^{(\alpha,u)} $
is a shortest path from $u$ to $a,$
we finally get 
\[l=|R(\rho)| \geq |(\gamma^*)^{(\alpha,u)}|=r \geq 20.  \]
\fine

For the rest of this section let the integer $L$ be chosen as $L=6|a|+1$. 
Given lattice paths $\rho_1, \ldots, \rho_L$ we will now define the ``Lucky Patch''
$\Theta(\rho_1, \ldots, \rho_L)$ as the first one to be hit by 7 paths
before these hit the boundary of the growing cluster. To be more precise, let $\rho_1, \ldots, \rho_L$ be lattice paths starting at
$\Delta=\Delta(a,L),$ each of which hits $a.$ For each patch $\pd$ and
index $j \leq L$ we count those paths among $\rho_1, \ldots, \rho_j,$ which
hit $\pd$ (strictly) before hitting the growing cluster, i.e., we define
\[N_{j,{\cal D}}(\rho_1, \ldots, \rho_L) = \sum_{k=1}^{j} {\bf 1}_{
\{\tau({\cal D},\rho_k) < \tau(\partial {\cal C}(\rho_1, \ldots,
\rho_{k-1}),\rho_k)\}} \]
and we let 
\[ V_7(\rho_1, \ldots, \rho_L) = \inf \{j: \mbox{ there exists a patch
${\cal D}$  
with $N_{j,{\cal D}}(\rho_1, \ldots, \rho_L)=7$} \}. \]
Let $\rho=\rho_{V_7(\rho_1, \ldots, \rho_L) }.$ Let us look at the first point $\rho_i$ on $\rho$ which belongs to a patch whose count of visitors increases to $7$ at time $V_7(\rho_1, \ldots, \rho_L)$.  Note that $\rho_i$ may belong to two different patches $\pd$ and $\tilde{\pd}$ and each of these may have its count of visitors increase to $7$. However, if this is the case, then exactly one of these patches, say $\pd,$ satisfies $\rho_i \in \gamma^+(\pd),$ and the other patch satisfies $\rho_i \in \gamma^-(\tilde{\pd}).$ We can thus define $\Theta(\rho_1, \ldots, \rho_L)$ as follows: If $\rho_i$ belongs to two patches and both of these have their count of visitors increase to $7$, we define $\Theta(\rho_1, \ldots, \rho_L)$ to be the one patch $\pd$ among these two, which satisfies $\rho_i \in \gamma^+(\pd).$ Otherwise, (i.e., $\rho_i$ belongs to exactly one patch, or $\rho_i$ belongs to two patches, but only one of these has its count of visitors increase to $7$), we define $\Theta(\rho_1, \ldots, \rho_L)$ to be the patch which contains $\rho_i$ and which has its count of visitors increase to $7$.
Next we define $V_1(\rho_1, \ldots, \rho_L) < \ldots <V_6(\rho_1, \ldots,
\rho_L)$ to be the indices of the first six paths which hit $\Theta(\rho_1,
\ldots, \rho_L)$ before hitting the growing cluster.

\vspace{7 mm}

The next proposition shows that up to time $V_7(\rho_1, \ldots, \rho_L)$
the newly added points can not form big clusters.

\begin{proposition} \label{third prop of patches}
Let $j < V_7(\rho_1, \ldots, \rho_L). $ Then any component of ${\cal
C}(\rho_1, \ldots, \rho_j) \setminus a$ has at most 12 points.
\end{proposition}
{\bf Proof:} Let $C$ be a component of ${\cal C}(\rho_1,
\ldots, \rho_j) \setminus a.$ Choose $u \in C$ with $u \in \partial a$ and
choose a neighbor $v$ of $u$ with $v \in a.$ Let $\pd_i$ be the
patch satisfying $u \in \pd^*_i$ or $u \in \gamma^-(\pd_i).$

Case i) Suppose that $u$ is in the ``left'' part of $\pd_i,$ i.e., $u \in \gamma^-(\pd_i)$ or $u$ is in the open region bounded by
the curve $\Gamma(\xi^l(\pd_i), \xi^m(\pd_i)),$ the lattice paths
$\gamma^-(\pd_i)$ and $\gamma^*(\pd_i)$ and the corresponding part
of $G(a).$ We now claim that in this case $C$
contains no path connecting $u$ to the curve $\Lambda$ obtained by joining
$R(\gamma^*(\pd_{i-1})), $ $\Gamma(\xi^m(\pd_{i-1}),\xi^r(\pd_i))$
and $\gamma^+(\pd_i).$ 

So suppose that $\pi$ is such a path. Without loss of generality assume
that $\pi$ is nonselfintersecting. Let $l$ be the first index with
$\pi_l \in \Lambda.$ Note that any vertex $z$ on $\pi^l,$ $z \neq \pi_l,$
satisfies $z \in \pd_{i-1}^* \cup \pd_i^*$ or $z \in \gamma^-(\pd_i) \setminus 
\xi^l(\pd_i).$ Thus $z$ got included into the growing cluster by a path
$\rho$ which increases $N_{. {\cal D}(i-1)}(\rho_1, \ldots \rho_L)$ or
$N_{. {\cal D}(i)}(\rho_1, \ldots \rho_L).$  Since $j < V_7(\rho_1, \ldots
\rho_L)$ at most 12 of the paths $\rho_1, \ldots ,\rho_j$ can have this
property. Hence
\begin{equation} \label{eq:l leq 13}
l \leq 12+1 =13.
\end{equation}

Now, if $\pi_l \in \gamma^*(\pd_{i-1}),$ then $(v, \pi^l)$ is a path
connecting $a$ and $\gamma^*(\pd_{i-1})$ with $\pi_k \not \in a$ for $k=1,
\ldots ,l-1$ and such that $\pi_n \in \gamma^-(\pd_i)$ for some $1 \leq n
\leq l-1.$ Hence by Proposition \ref{second prop of patches} $|(v,\pi^l)| \geq 20$ which
contradicts (\ref{eq:l leq 13}). 

Next, if $\pi_l \in \gamma^+(\pd_{i}),$ then there exists an index $k<l$
with $\pi_k \in \gamma^*(\pd_i).$ Let $w=\pi_k.$ Since the path $(v, \pi^k)$
connects $a$ and $w \in \gamma^*(\pd_i)$ we have $|(v, \pi^k)| \geq
|(\gamma^*(\pd_i))^{(\alpha,w)}|.$ Hence, letting ${\tilde \pi}$ denote the path
obtained by joining $R(\gamma^*(\pd_i))^w$ and $(\pi^l)^{(\alpha,k)}$ we
get
$|{\tilde \pi}| = |R(\gamma^*(\pd_i))^w| + |(\pi^l)^{(\alpha,k)}| \leq
|(v, \pi^k)| + |(\pi^l)^{(\alpha,k)}| = |(v, \pi^l)| = l+1 \leq 13+1 =14.$
Yet, ${\tilde \pi}$ connects $y^*(\pd_i)$ and $\gamma^+(\pd_i)$ and thus by
definition of $\gamma^+(\pd_i),$ its length is at least 40. Once again we get a
contradiction. 

Finally, if $\pi_l \in \Gamma$ then $(v,\pi^l)$ connects $a$ and $\Gamma$
and hence $|(v, \pi^l)| \geq 40$ or $l \geq 39.$ 
Once again we get a contradiction and so our claim follows.

Case ii)  Suppose that $u$ is in the ``right'' part of $\pd_i,$ i.e.,  $u \in \gamma^*(\pd_i)$ or $u$ is in
the open region bounded by the curve $\Gamma(\xi^m(\pd_i), \xi^+(\pd_i)),$
the paths $\gamma^*(\pd_i)$ and $\gamma^r(\pd_i)$ and the corresponding part
of $G(a).$ Then it can be shown in a similar way that $C$ contains no
path connecting $u$ to the curve ${\tilde \Lambda}$ obtained by joining
$R(\gamma^-(\pd_i)),$ $\Gamma(\xi^l(\pd_i),\xi^m(\pd_{i+1}))$ and
$\gamma^m(\pd_{i+1}).$

Thus in case (i) all points of $C$ belong to ${\cal D}_{i-1}^* \cup {\cal
D}_{i}^*$ or to $\gamma^-(\pd_i).$
In particular, any path $\rho$
which leads to the inclusion of such a point into the growing cluster,
increases $N_{. {\cal D}(i-1)}$ or $N_{. {\cal D}(i)}.$ Since $j <
V_7(\rho_1, \ldots \rho_L),$ we get $|C| \leq 12.$ Similarly in case~(ii)
all points of $C$ belong to  ${\cal D}_{i}^* \cup {\cal
D}_{i+1}^*$ or to $\gamma^+(\pd_i)$ and once again we can conclude
$|C| \leq 12.$  \newline
\vspace{0 mm}
\fine

\begin{proposition} \label{fourth prop of patches}
Let $j < V_7(\rho_1, \ldots, \rho_L)$ and let $\pd$ be a patch of $a$. Let $z \in
{\cal C}(\rho_1, \ldots, \rho_j) \setminus a$ be a point with $z \in {\cal B}$ for
some patch ${\cal B} \neq \pd.$ Then any path $\lambda$ from $z$ to $\gamma^*(\pd)$
with $\lambda_i \not \in a$ for $i < |\lambda|$ has length at least 14.
\end{proposition}
(Note that we include in particular the case $z \in \gamma^{\pm}(\pd),$
since any lattice point on $\gamma^{\pm}(\pd)$ belongs to two patches.)

{\bf Proof.} Suppose $\lambda$ is a path from $z$ to $\gamma^*(\pd)$ with
$|\lambda| \leq 13$ and $\lambda_i \not \in a$ for $i < |\lambda|.$ Let $w \in
\gamma^*(\pd)$ be the endpoint of $\lambda.$ By the previous proposition there
exists a path $\lambda^*$ from $a$ to $z$ with $|\lambda^*| \leq 12$ and $\lambda^*_i \not
\in a$ for $i>0.$ Thus $(\lambda^*,\lambda)$ connects $a$ and $w \in
\gamma^*(\pd).$ Hence $|\gamma^*(\pd)^{(\alpha,w)}| \leq |(\lambda^*,\lambda)| \leq
12+13=25.$ So ${\hat \lambda}=(R(\gamma^*(\pd))^w,R(\lambda))$ is a path of length at
most $25+13=38$ connecting $y^*(\pd)$ and $z.$ Since $z$ belongs to a patch
${\cal B}$ with ${\cal B} \neq \pd$ and ${\hat \lambda}$ avoids $\Gamma$ (recall that any lattice
point on $\Gamma$ has distance 40 from $a$) ${\hat \lambda}$ hits $\gamma^-(\pd)$ or
$\gamma^+(\pd).$ We thus get a contradiction since by construction any path
which connects $y^*(\pd)$ to $\gamma^-(\pd)$ or $\gamma^+(\pd)$ and all
interior points of which belong to $\Ztwo \setminus a$ has length at least
40. \newline
\vspace{14 mm}
\fine

\section{Choice of the New Positions of Attaching}
 
For each patch $\pd$ we will define lattice points $x_1, \ldots, x_7$
close to $\gamma^*(\pd)$ and paths $\beta_1, \ldots, \beta_7$ starting at
$\xi^m(\pd)$ such that the following holds:
\begin{enumerate}
\item ${\cal C}(\beta_1, \ldots ,\beta_j)= a \cup \{ x_1, \ldots, x_j \}$ for $j=1, \ldots, 7.$ 
\item $a \cup \{ x_1, \ldots, x_7 \}$ has at least one more hole than $a.$
\end{enumerate}

Later on we will insert a loop into the path $\omega_{V_j(\omega)}$ at the hitting time of
$\Theta(\omega)$ which contains
$(\beta_j(\Theta(\omega)),R(\beta_j(\Theta(\omega))))$ and we will get
the point
$x_j(\Theta(\omega))$ as new point of attaching for the path $\varphi(\omega)_{V_j(\omega)}.$ 

So let $$\gamma=\gamma^*(\pd), y=\gamma_{40}  \mbox{ and } w=\gamma_{39}. $$ 
Note that $y$ belongs to $a$ and $w$ belongs to the unbounded component of $\Ztwo \setminus a.$  Let $f_2=w-y$ and let $f_1$ be a unit vector orthogonal to $f_2$ satisfying that, when embedding ${\bf R}^2$ into ${\bf R}^3$ in the usual way, (i.e., $x=(x_1,x_2)$ corresponds to ${\hat x}=(x_1,x_2,0)$ for all $x \in {\bf R}^2$), then ${\hat f}_1 \times {\hat f}_2 = (0,0,1).$ 
Let $w_1=w-f_1-f_2$, $w_2=w-f_1$, $w_3=w-f_1+f_2$, $w_4=w+f_2$, $w_5=w+f_1+f_2$, $w_6=w+f_1$ and $w_7=w+f_1-f_2$.
So, if we put $w_0=y,$ then $\Lambda=(w_0,w_1, \ldots, w_7,w_0)$ forms a simple loop surrounding $w.$ 
Now choose $k_1$ such that $w_{k_1}$ is the first lattice point on $\gamma$ which belongs to $
\{w_1, \ldots ,w_7 \}$, i.e., $w_{k_1}=\gamma_{\tau(\{w_1, \ldots ,w_7 \})}.$ Let 
${\tilde \beta}(w_{k_1})=\gamma^{\{w_1, \ldots ,w_7 \}}.$

Depending on $\gamma$ we will next choose a neighbor $w_{k_2}$ of $w_{k_1}$
and a path ${\tilde \beta}(w_{k_2})$ from $\xi^m(\pd)$ to $w_{k_2}.$
Moreover we will have to decide which of the points $w_{k_1}$ and $w_{k_2}$
we want to have included first into the cluster. 
For instance, if $(\gamma_{36},\gamma_{37},\gamma_{38})=(w_5+f_2,w_5,w_4),$
then $w_{k_1}=w_5$ and we choose $k_2=4,$  ${\tilde
\beta}(w_4)=(\gamma^{36}, w_4+f_2,w_4)$ and we decide to have $w_4$
included into the cluster before including $w_5.$ 
Once we have $w_{k_2}$ and ${\tilde \beta}(w_{k_2})$ selected and we know
which point to include first, we proceed as follows:
Let $k_{-} = k_1 \wedge k_2 \mbox{ , } k_{+} = k_1 \vee k_2 ,$ $m_1=\max \{ j \leq k_{-}: w_j \in \partial a \},$ $m_2=\min \{ j \geq k_{+}: w_j \in \partial a \}$ and let $w_{-}=w(k_{-})$ and $w_{+}=w(k_{+}).$ 
Let $l_1=k_{-}-m_1$ and for $1 \leq j \leq l_1$ define 
\[ x_j=w_{m_1+(j-1)} \mbox{  and  } \beta_j=\left({\tilde \beta}(w_{-}),
\left(R(\Lambda)^{x_j} \right) ^{(\alpha,w_{-})} \right). \]
Thus the path $\beta_j$ moves from $\xi^m(\pd)$ to $w_{-}$ along ${\tilde
\beta}(w_{-})$ and then moves counterclockwise around the loop $\Lambda$
until it reaches $x_j.$ 
Similarly, let $l_2=m_2-k_{+}$ and for $1 \leq j \leq l_2$ define 
\[ x_{l_1+j}=w_{m_2-(j-1)} \mbox{  and  } \beta_{l_1+j}=\left({\tilde \beta}(w_{+}),
(\Lambda^{x_{l_1+j}} ) ^{(\alpha,w_{+})} \right). \]
If we decide to have $w_{k_1}$ included before $w_{k_2},$ we put 
\[x_{l_1+l_2+1}=w_{k_1} \mbox{ and } \beta_{l_1+l_2+1}={\tilde \beta}(w_{k_1}), \]
\[x_{l_1+l_2+2}=w_{k_2} \mbox{ and } \beta_{l_1+l_2+2}={\tilde \beta}(w_{k_2}), \]
otherwise we exchange $k_1$ and $k_2.$
Finally, we let $n=\tau(\partial\{w_1, \ldots, w_7\},\gamma)$ and $l_3=7-(l_1+l_2+2),$
 and we define for $1 \leq j \leq l_3,$
\[ x_{l_1+l_2+2+j}= \gamma_{n-(j-1)}, \mbox{ and }
\beta_{l_1+l_2+2+j}=\gamma^{n-(j-1)}.\]
Thus the last points $x_j$ are chosen on the path $\gamma$ and  we obtain the corresponding path $\beta_j$ by stopping $\gamma$ at the hitting time of $x_j.$
In the previous example, if, say, $m_1=3$ and $m_2=7,$ then
$ x_1=w_3$, $x_2=w_7$, $x_3=w_6$, $x_4=w_4$, $x_5=w_5$, $x_6=\gamma_{36}$, $x_7=\gamma_{35}$, $\beta_1=({\tilde \beta}(w_4),w_3)$, $\beta_2=(\gamma^{w_5},w_6,w_7)$, $ \beta_3=(\gamma^{w_5},w_6)$, $\beta_4={\tilde \beta}(w_4)$, $\beta_5=\gamma^{w_5}$, $\beta_6=\gamma^{36}$ and  $\beta_7=\gamma^{35}$.

Let us now indicate our choices depending on $\gamma^{(\alpha,36)}:$ 
\vspace{7 mm}

\begin{tabular}{lll|l|l}
  && $(\gamma_{36},\gamma_{37},\gamma_{38})$ & $k_2$ & ${\tilde \beta}(w(k_2))$
\\ \hline
  &&$(w_6+2f_1,w_6+f_1,w_6)$ & 5 & $(\gamma^{37},w_5+f_1,w_5)$  \\
  &&$(w_6+f_1+f_2,w_6+f_1,w_6)$ & 5 & $(\gamma^{36},w_5)$       \\
  &&$(w_5+f_1,w_5,w_6)$ & 6 & $(\gamma^{36},w_6+f_1,w_6)$  \\
  &&$(w_5+f_2,w_5,w_6)$ & 4 & $(\gamma^{36},w_4+f_2,w_4)$  \\
  &&$(w_5+f_1,w_5,w_4)$ & 6 & $(\gamma^{36},w_6+f_1,w_6)$  \\
  &&$(w_5+f_2,w_5,w_4)$ & 4 & $(\gamma^{36},w_4+f_2,w_4)$  \\
  &&$(w_4+f_1+f_2,w_4+f_2,w_4)$ & 5 & $(\gamma^{36},w_5)$  \\
  &&$(w_4+2f_2,w_4+f_2,w_4)$ & 5 & $(\gamma^{37},w_5+f_2,w_5)$  \\
\end{tabular}
 
\vspace{7 mm}

\noindent If ${\tilde \beta}(w_{k_2})$ uses a lattice point $z \neq w_{k_2}$ with $z
\not \in \gamma,$ let $w_{k_2}$ get included into the cluster before
$w_{k_1},$ otherwise let $w_{k_1}$ get included before $w_{k_2}.$ 

Recalling that $|\gamma_j-y|=40-j$ for $0 \leq j \leq 40,$ we can easily
see that all other possible cases get obtained by reflecting one of the
above paths about the straight line through $y$ and $w.$ 

Let us note that $\gamma^{35}$ does not hit $\partial \{w_1, \ldots w_7 \}$
since $|w_i-y| \leq 3$ for $i=1 \ldots 7$ and $|\gamma_j - y| =40 - j \geq
5$ for $j \leq 35.$
Moreover, if ${\tilde \beta}(w_{k_2})$ uses a lattice point $z \neq
w_{k_2}$ with $z \not \in \gamma,$ then $z$ is a neighbor of $\gamma^{37}$
or $\gamma^{36}.$ Thus, $\mbox{dist}(z,a) \geq \mbox{dist}(\gamma^{37},a)-1 = 2 $ and $z
\not \in a \cup \partial a.$ Hence, for $j \leq l_1+l_2+2,$ $\beta_j$ does
not hit $\partial (a \cup \{ x_1, \ldots ,x_{j-1} \} )$ (strictly) before
hitting $w_{k_1}$ or $w_{k_2}.$ Thus it is straightforward to check that
property 1) holds for all indices $1 \leq j \leq 7.$

\begin{remark} \label{rem:new positions}
i) For each point $x_i(\pd),$ $1 \leq i \leq 7,$ there exists a path of
length at most 7 connecting $y^*(\pd)$ and $x_i,$ all the inner lattice 
points of which belong to $\Ztwo \setminus a.$ In particular, each of the
points $x_i(\pd)$ belongs to $\pd^*.$ \newline
ii) Let $w = \gamma^*(\pd)_{39}.$ For each lattice point $z$ on $\beta_i(\pd),$ $1 \leq i \leq 7,$ there
exists a path ${\tilde \rho}$ in $\Ztwo \setminus a$ connecting $z$ to
$\gamma^*(\pd)^w$ with $|{\tilde \rho}| \leq 2.$
\end{remark}

The next proposition will give us the second property of the set $\{x_1(\pd),
\ldots , x_7(\pd)\}.$

\begin{proposition} \label{proposition in "New Positions"}
Let $\pd$ be a patch and let $w = \gamma^*(\pd)_{39}.$ Then $w$ belongs to a
hole of $a \cup \{x_1(\pd) \ldots , x_7(\pd)\}.$
\end{proposition} 
{\bf Proof.} Let $x_i=x_i(\pd)$ and let $w_1, \ldots, w_7,$ $m_1, m_2$ and $k_1$
be defined as before. Let $\hat a = a \cup \{x_1 \ldots , x_7\}.$ Note that
$\{x_i: i=1 \ldots 7 \} \cap \{w_i: i=1 \ldots 7 \} = \{ w(m_1),w(m_1+1),
\ldots, w(m_2-1), w(m_2)\}$ and $m_1 < k_1 < m_2.$ By definition of $m_1$
and $m_2$ we have $w(m_i) \in \partial a,$ $i=1,2.$ So there exist
neighbors $v_1$ of $w(m_1)$ and $v_2$ of $w(m_2)$ with $v_i \in a,$
$i=1,2.$ Since $a$ is connected there exists a nonselfintersecting path
$\rho$ in $a$ from $v_1$ to $v_2$. Note that none of the points $w(m_1),
\ldots, w(m_2)$ belongs to $\rho,$ the path $\rho$ avoids $\Gamma$ and
$\rho$ intersects with $\gamma$ at most in the point $y.$ Thus
$(\rho,w(m_2),w(m_2-1), \ldots, w(m_1+1), w(m_1), v_1)$ is a simple loop in
$\hat a$ which surrounds $w$ and $w$ belongs to a hole of $\hat a.$ 
\newline
\vspace{0 mm}
\fine

\section{Definition of the Map $\varphi$ and Verification of Property \ref{first property of varphi}}

As before, let us fix an integer $N_0$ satisfying Remark  
\ref{rem:size}.
Let $a$ be a finite cluster with $|a| \geq N_0$ and let $L=6 |a| +1$. Let $\Xi=\{ \omega \in \Omega^{L,a}: \tau(\omega_i,a) < \infty, \mbox{ $i=1, \ldots, L$ } \}$ 
and for $\omega \in \Xi$ define $\omega^a=(\omega_1^a, \ldots, \omega_L^a).$
 
For $\omega \not \in \Xi$ we will define $\varphi(\omega)=\omega$ and for $\omega \in \Xi$ we will define $\varphi(\omega)$ by inserting extra loops
into some of the $L$ given paths $\omega_1, \ldots, \omega_L.$ The insertion of an extra loop in, say, $\omega_k,$ will be done at the hitting time of  
$\partial A_{k-1}(\omega) \cup \Theta(\omega).$ If $k \not \in  \{V_i(\omega), 1 \leq i \leq 7\},$ we will choose the
loop such that it hits $\partial A_{k-1}(\varphi(\omega))$ at a point $u$ which already belongs to the cluster
$A_{k-1}(\omega)$ and if $k=V_i(\omega)$ for some index $1 \leq i \leq 7,$ then we will choose the loop such that it hits $\partial A_{k-1}(\varphi(\omega))$ at the point $x_i(\Theta(\omega)).$ Thus, letting $n(\omega,j)=\max
\{i: V_i(\omega) \leq j \}$ we will get that for any index $j \leq V_7(\omega)$
\begin{eqnarray*}
\{x_1(\Theta(\omega)), \ldots ,x_{n(\omega,j)}(\Theta(\omega))\}  & \subset & A_{j}(\varphi(\omega)) \\  
      & \subset & A_{j}(\omega) \cup \{x_1(\Theta(\omega)),\ldots ,x_{n(\omega,j)}(\Theta(\omega))\} .
\end{eqnarray*}
Now, say, a loop has to be inserted into the path $\omega_k,$ where $k=V_i(\omega), 1 \leq i \leq 7,$ and let us assume that the point where $\omega_k$ hits $\partial A_{k-1}(\omega) \cup \Theta(\omega)$ belongs to $\gamma^{-}(\Theta(\omega)).$ Then we need to construct a path from this hitting point to $x_i(\Theta(\omega)).$ We will do this by first moving along
$R(\gamma^{-}(\Theta(\omega)))$ up to $\xi^l(\Theta(\omega))$ (a slight
modification is needed in this part of the construction since we may hit
$\partial A_{k-1}(\varphi(\omega))$
along this path), then following a lattice path to $\xi^m(\Theta(\omega))$ (this path will get easily obtained from the curve $\Gamma$) and finally descending along $\beta_i(\Theta(\omega))$ to $x_i.$ 
 
So, let us first obtain from the curve $\Gamma$ a simple lattice loop $\hat
\Gamma=(\psi_0, \psi_1, \ldots , \psi_{\hat l})$ as follows:
Let the lattice points on $\Gamma$ be enumerated as $\xi_0, \ldots ,\xi_l$
and put $\xi_{l+1}=\xi_0.$ Let $\psi_0=\xi_0$ and suppose that $\psi_0,
\ldots \psi_k$ have already been constructed with $\psi_k=\xi_j$ for some
$j \leq l.$ 
If $|\xi_{j+1}-\xi_j| = 1$ let $\psi_{k+1}=\xi_{j+1},$ otherwise let
$\psi_{k+2}=\xi_{j+1}$ and define
\[ \psi_{k+1}= \left \{ \begin{array}{lll}
                       \xi_j+e_2  & \mbox{if} &  \xi_{j+1}=\xi_j+e_1+e_2 \\
                       \xi_j+e_1  & \mbox{if} &  \xi_{j+1}=\xi_j+e_1-e_2 \\
                       \xi_j-e_1  & \mbox{if} &  \xi_{j+1}=\xi_j-e_1+e_2 \\
                       \xi_j-e_2  & \mbox{if} &  \xi_{j+1}=\xi_j-e_1-e_2
                        \end{array}
               \right. .\]
Note that by construction each vertex $\psi_j$ on ${\hat \Gamma}$ satisfies
$\mbox{ dist}(\psi_j,a) \geq 40$ and Proposition \ref{first prop of patches} is also valid for
$\hat \Gamma.$

\vspace{7 mm}

We will now formally define the map $\varphi.$ We will do this such that the decisions whether to insert a loop in, say, $\omega_k$, and which loop to insert will depend only on $\omega^a$. 
 
If $\omega \not \in \Xi$ let $\varphi(\omega)=\omega.$ Otherwise let us assume that $\f(\omega)_1, \ldots ,\f(\omega)_{k-1}$ have already been constructed, this construction being such that $A_{k-1}(\varphi(\omega))=A_{k-1}(\varphi({\tilde \omega}))$ for all $\tilde \omega$ satisfying ${\tilde \omega}^a=\omega^a$. Let  
$\pd=\Theta(\omega)$. (Before proceeding further let us note that $\Theta(\omega)=\Theta({\tilde \omega})$ for all $\tilde \omega$ satisfying ${\tilde \omega}^a=\omega^a$ and for all such ${\tilde \omega}$ we also have   $V_i(\omega)=V_i({\tilde \omega})$, $1 \leq i \leq 7$).   
 
If $k>V_7(\omega)$ define $\f(\omega)_k=\omega_k.$ \newline

If $k \leq V_7(\omega)$ we consider the following three cases: 

Case i) $\tau(\partial A_{k-1}(\omega),\omega_k)
= \tau(\partial A_{k-1}(\f(\omega)),\omega_k)
\leq \tau(\pd,\omega_k),$ i.e., 
the path $\omega_k$ hits $\partial A_{k-1}(\omega)$ before
or upon hitting $\pd$ and at the same time $\partial A_{k-1}(\f(\omega))$ gets hit. Define
\[ \f(\omega)_k=\omega_k. \]

Case ii) $\tau(\partial A_{k-1}(\omega),\omega_k)
\leq \tau(\pd,\omega_k)$ and $\tau(\partial A_{k-1}(\omega),\omega_k) < \tau(\partial A_{k-1}(\f(\omega)),\omega_k),$ i.e., 
the path $\omega_k$ hits $\partial A_{k-1}(\omega)$ before
or upon hitting $\pd$ and up to (including) this time $\partial A_{k-1}(\varphi(\omega))$ has not yet been hit. 
Fix a path $\lambda$ from $\omega_{k,\tau(\partial A_{k-1}(\omega))}$ to
$\partial a$ with $\lambda_j \in A_{k-1}(\omega)$ for $j>0$, this choice being the same for all ${\tilde \omega}$ satisfying ${\tilde \omega}^a=\omega^a$.   
Let $j_0= \min \{j:
\lambda_j \in \partial A_{k-1}(\f(\omega))\}$ and define
\[\f(\omega)_k=(\omega_k^i, \lambda^{j_0},
R(\lambda^{j_0}),\omega_k^{(\alpha,i)}), \]
where $i=\tau(\partial A_{k-1}(\omega), \omega_k).$ 

Case iii) $\tau(\pd,\omega_k) < \tau(\partial A_{k-1}(\omega),\omega_k),$ i.e.,
the path $\omega_k$ hits $\pd$ before hitting $\partial A_{k-1}(\omega).$
Then $k=V_i(\omega)$ for some index $1 \leq i \leq 7.$
Let $j=\tau_{\pd}(\omega_k).$ We will now construct a path
$\alpha_i(\omega)$ from $\omega_{kj}$ to $x_i(\pd)$ which hits $\partial
A_{k-1}(\f(\omega))$ at $x_i(\pd)$, this path, once again, being the same for all ${\tilde \omega}$ satisfying ${\tilde \omega}^a=\omega^a$. \newline
If $\omega_{kj} \in
\Gamma$ let ${\hat \Gamma}(i,\omega)$ denote that part of the loop $\hat
\Gamma$ or $R({\hat \Gamma})$ which starts at $\omega_{kj},$ ends at
$\xi^m(\pd)$ and which satisfies that all of its vertices which belong to
$\Gamma$ also belong to $\pd.$ Let $$\alpha_i(\omega)=({\hat
\Gamma}(i,\omega), \beta_i(\pd)).$$ 
If $\omega_{kj} \in \gamma^-(\pd)$ (resp. $\gamma^+(\pd)$) let ${\cal
R}_i(\omega)$ denote the collection of all paths $\lambda=(\lambda_0,
\ldots , \lambda_l),$ $l\in {\bf N}_0,$ satisfying $\lambda_0 =
\omega_{kj},$ $\lambda_l=\xi^l(\pd)$ (resp. $\lambda_l=\xi^r(\pd)$) and
$\lambda_m \not \in  \pd^*$ $ \cup $ $ \partial A_{k-1}(\f(\omega))$ for 
$m=0, \ldots, l.$ Thus, the
paths in ${\cal R}_i(\omega)$ connect $\omega_{kj}$ to $\xi^l(\pd)$ or
$\xi^r(\pd)$ without hitting the cluster built by the paths $\f(\omega)_1,
\ldots \f(\omega)_{k-1}$ or the interior of ${\cal D}.$ Among all the paths
in ${\cal R}_i(\omega)$ fix a path of minimum length, the same path being chosen for all ${\tilde \omega}$ satisfying ${\tilde \omega}^a=\omega^a$. (This is possible since ${\cal R}_i(\omega)={\cal R}_i({\tilde \omega})$ for all such ${\tilde \omega}$). Call this path
$r(i,\omega).$ If $\omega_{V_i \tau({\cal D})}$ belongs to $\gamma^-({\cal D})$ let ${\hat \Gamma}(i,\omega)$ denote
the part of the loop $\hat
\Gamma$ which starts at $\xi^l({\cal D})$ and ends at $\xi^m({\cal D})$ and if $\omega_{V_i \tau({\cal D})}$  belongs to $\gamma^+({\cal D})$ let ${\hat \Gamma}(i,\omega)$ denote the part of the loop $R(\hat \Gamma)$ which starts at $\xi^r({\cal D})$ and ends at $\xi^m({\cal D}).$  
Let 
$$\alpha_i(\omega)=(r(i,\omega),{\hat \Gamma}(i,\omega),
\beta_i(\pd)).$$ 
For all $\omega \in \Xi$ satisfying the assumptions of case (iii) define 
\[\f(\omega)_k=(\omega_k^{\tau({\cal D})}, \alpha_i(\omega),
R(\alpha_i(\omega)), \omega_k^{(\alpha,\tau({\cal D}))}). \] 

Any loop which gets inserted due to the occurence of case (ii) is called a catch-up loop.
Note that a catch-up loop may lead to the
inclusion of a point into $A_.(\varphi(\omega))$ which belongs to the interior of $\Theta(\omega)$ but which is different from any of the points
$x_1(\Theta(\omega)), \ldots, x_7(\Theta(\omega)).$ 
However, as we will show in the next proposition, any such point is
sufficiently far away from $\gamma^*(\Theta(\omega)).$ 

Let $z_i(\omega)=\omega_{V_i(\omega),\tau(\partial A_{V_i(\omega)-1}(\omega))}$
be the new lattice point for the cluster
$A_{V_i}(\omega)$ and let $x_i(\omega)=x_i(\Theta(\omega)),$ $i=1,
\ldots, 7.$

\begin{proposition} \label{simple properties of varphi}
The map $\f: \Omega^{L,a} \rightarrow \Omega^{L,a}$ is well-defined and for all $\omega \in \Xi$ and
indices $k \leq V_7(\omega)$ the map $\f$ satisfies
\renewcommand{\labelenumi}{\theenumi)} 
\begin{enumerate}
\item   
$\{x_1(\omega), \ldots ,x_{n(\omega,k)}(\omega)\}  
 \subset A_k(\f(\omega)) 
 \subset A_k(\omega) \cup \{x_1(\omega),\ldots ,x_{n(\omega,k)}(\omega)\}.$ \label{P1}
\item If $x \in A_k(\omega) \setminus A_k(\f(\omega))$ 
then there exists a path
$\rho=(\rho_0, \ldots, \rho_l)$ with $\rho_0=x,$ $\rho_j \in
A_k(\omega) 
\setminus a$ for $j=0, \ldots, l,$ and
$\rho_l=z_i(\omega)$ for some index $1 \leq i \leq n(\omega,k).$ \label{P2}
\item If $z \in A_k(\f(\omega)) \setminus
(\{x_1(\omega), \ldots ,x_{n(\omega,k)}(\omega)\} \cup a)$
then any path $\rho$ from $z$ to $\gamma^*(\Theta(\omega))$
with $\rho_j \not \in a$ for $j < |\rho|$ has length at least 8. \label{P3}
\end{enumerate}
\end{proposition}

\vspace{3.5 mm}
Before proving the proposition let us state the following two corollaries

\begin{corollary}
For all $\omega \in \Xi$ we have 
\[ H(A_L(\f(\omega))) \geq H(a)+1 . \]
\end{corollary}
{\bf Proof.} Let $\omega \in \Xi$. First note that if $B \subset \Ztwo \setminus a$ is a hole of
$a$ then $B$ is also a hole of $A_L(\varphi(\omega))$ since any random walk will hit $a$
before hitting $B.$ Thus it suffices to show that there exists a point $w$
which belongs to the unbounded component of $\Ztwo \setminus a,$ but which
belongs to a finite component of $\Ztwo \setminus A_L(\varphi(\omega)).$ Let
$\pd=\Theta(\omega)$ and let $w=\gamma^*(\pd)_{39}.$ Since in $\Ztwo
\setminus a$ the point $w$ is connected to $\Gamma$ by the path
$\gamma^*(\pd),$ it clearly belongs to the unbounded component of $\Ztwo
\setminus a.$ Next, noting that $w \in \partial a$ and $w \not \in
\{x_1(\pd), \ldots, x_7(\pd) \},$ we get by part (3) of the previous
proposition that $w$ does not
belong to $A_{V_7(\omega)}(\f(\omega)).$ By Proposition \ref{proposition in "New Positions"} the point $w$
belongs to a hole of $a \cup \{x_1(\pd), \ldots, x_7(\pd) \}$ and by
property (1), $a \cup \{x_1(\pd), \ldots, x_7(\pd) \} \subset
A_{V_7(\omega)}(\f(\omega)).$ Thus $w$ belongs to a hole of
$A_{V_7(\omega)}(\f(\omega)).$ But any hole of
$A_{V_7(\omega)}(\f(\omega))$ is also a hole of $A_L(\f(\omega))$ and our
claim follows.
\newline
\vspace{5 mm}
\fine

\begin{corollary} \label{second corollary of simple properties}
Let $\omega \in \Xi.$ Then
\[
\begin{array}{rllc} 
|A_{V_7(\omega)}(\f(\omega)) \cap \Theta^*(\omega)| & \geq & 7  & \mbox { and}\\
|A_{V_7(\omega)}(\f(\omega)) \cap B^*| &  \leq &  6  & \mbox{  for any patch $B
\neq \Theta(\omega)$.} 
\end{array} \]
\end{corollary}
{\bf Proof.}  Let $\omega \in \Xi$. Recalling that for any patch $\pd$ we have $\{x_1(\pd), \ldots
,x_7(\pd)\} \subset \pd^*,$ the first inequality follows immediately from part (1) of Proposition \ref{simple properties of varphi}.
 Now let $z \in A_{V_7(\omega)}(\f(\omega)) \cap B^*$ where
$B$ is a patch with $B \neq \Theta(\omega).$ Since $x_7(\Theta(\omega))$ is the last point which
got included into $A_{V_7(\omega)}(\f(\omega)),$ the point $z$ already
belonged to the previous cluster $A_{V_7-1}(\f(\omega))$ and  using the second set inclusion in part (1) of Proposition \ref{simple properties of varphi} we can conclude that $z$ belongs to $A_{V_7-1}(\omega).$ Thus $z=\omega_{k \tau(\partial A_{k-1}(\omega))}$ for some index $1 \leq k < V_7(\omega).$ In particular
it follows that $\omega_k$ hits $B$ strictly before hitting
$\partial A_{k-1}(\omega)$ and hence $N_{.,B}(\omega)$ increases by 1 at time $k.$
Now $N_{V_7-1,B}(\omega) \leq 6.$ Thus at most six of the paths $\omega_1,
\ldots, \omega_{V_7-1}$ can have the above property and this implies the
second inequality.
\newline
\vspace{7 mm}
\fine

{\bf Proof of Proposition \ref{simple properties of varphi}.} Clearly $\varphi$ is well-defined on $\Omega^{L,a} \setminus  \Xi.$ So let $\omega \in \Xi.$ We will show by induction on $k$ that each of the paths $\varphi(\omega)_k$ is well-defined and that $A_k(\varphi(\omega))$ satisfies properties (1) - (3) for $k \leq V_7(\omega).$ 

First note that $A_0(\omega)=A_0(\varphi(\omega))=a$ and thus properties (1) - (3) hold for the index $k=0.$ So let $k \geq 1$ and assume that $\varphi(\omega)_1, \ldots, \varphi(\omega)_{k-1}$ are well-defined and that $A_{(k-1) \wedge V_7(\omega)}(\varphi(\omega))$ satisfies (1) - (3). If $k > V_7(\omega)$ then $\varphi(\omega)_k=\omega_k$ and thus the path $\varphi(\omega)_k$ is well-defined. So let $k \leq V_7(\omega).$  Let ${\cal D}=\Theta(\omega).$ 

If $\tau(\partial A_{k-1}(\omega),\omega_k)$  $ = $  $\tau(\partial A_{k-1}(\f(\omega)),\omega_k)$
$\leq$  $\tau(\pd,\omega_k)$ then the same new point (which belongs to a patch
$\cal{B}$ with ${\cal B} \neq \Theta(\omega)$) is added to 
$A_{k-1}(\f(\omega))$ and to $A_{k-1}(\omega).$ Thus properties (1) and (2)
hold by induction hypothesis   and (3) follows from Proposition \ref{fourth prop of patches}. 

Let us now consider the case that a catch-up loop gets inserted into $\omega_k.$ Let 
$x=\omega_{k, \tau(\partial A_{k-1}(\omega))}$ be the new node for the
cluster $A_{k}(\omega)$ and let $u=\lambda_{j_0}$ with the notation used
in the definition of $\f.$ Recall that $\lambda^{j_0}$ is
a path connecting $x$ to $u$ all the vertices of which (with the exception
of $\lambda_0$) belong to $A_{k-1}(\omega) \setminus A_{k-1}(\f(\omega)).$
Moreover, $u$ is the first vertex on $\lambda$ which belongs to $\partial
A_{k-1}(\f(\omega)).$ Since $\f(\omega)_k$ is obtained by inserting the loop
$(\lambda^{j_0},R(\lambda^{j_0}))$ into $\omega_k$ we can conclude
that $u$ (which belongs to $A_{k-1}(\omega) \setminus A_{k-1}(\f(\omega))$ and thus to $A_k(\omega)$) is the new point for the
cluster $A_{k}(\f(\omega)).$ Combining this with our induction hypothesis,  property (1) follows. 
 
As for (2), we only need to check this property for the vertex
$x.$  Now $x \in \partial A_{k-1}(\omega) \setminus \partial
A_{k-1}(\f(\omega)).$ So there exists a neighbor $z$ of $x$ with $z \in
A_{k-1}(\omega) \setminus A_{k-1}(\f(\omega)).$ By induction hypothesis $z$ is connected to a vertex $z_i(\omega),$
$1 \leq i \leq n(k,\omega),$ by a path ${\tilde \rho},$ all the vertices of which
belong to $A_{k-1}(\omega) \setminus a.$ Thus all the vertices of the path
$(x, \tilde \rho),$ which connects x and $z_i(\omega),$ belong to
$A_k(\omega) \setminus a$ and (2) follows. 

Let us now verify that property
(3) holds. By induction hypothesis we only need to check this property
for the point $u=\lambda_{j_0}.$ 
If $u$ belongs to a patch ${\cal B}$ with ${\cal B} \neq \pd$ then (3) follows
from Proposition \ref{fourth prop of patches}. Otherwise $u \in \pd^*.$  Recalling
that $x \not \in \pd^*$ let us define $j_1=\min\{j \leq j_0: \lambda_j \in
\pd \setminus \pd^* \}.$ Without loss of generality let us assume that
$\lambda_{j_1} \in \gamma^-(\pd).$ Define $j_2=\max\{j \leq j_0: \lambda_j \in
\gamma^-(\pd) \}$ and let \[
j_3= \left\{ \begin{array}{ll}
                j_0  &  \mbox{ if $(\lambda^{j_0})^{(\alpha,j_2)}$ does not hit
$\gamma^*(\pd)$}\\
                \min\{j_2 \leq j \leq j_0: \lambda_j \in \gamma^*(\pd) \} &
\mbox{otherwise.}
             \end{array}
     \right.  \]
Then all vertices on $(\lambda^{j_3})^{(\alpha,j_2+1)}$ are vertices of
$A_{k-1}(\omega)$ which belong to $\pd^*.$ Thus, if $v$ is such a vertex and if
the path $\omega_j$ leads to the inclusion of $v$ into the cluster
$A_j(\omega),$ then $\omega_j$ hits $\pd$ strictly before hitting $\partial
A_{j-1}(\omega).$ Hence $j=V_i(\omega)$ for some index $1 \leq i < 7.$ In
particular, we can conclude that 
\begin{equation} \label{*}
|(\lambda^{j_3})^{(\alpha,j_2)}| \leq 6. 
\end{equation}
Thus by Proposition \ref{fourth prop of patches} the path $(\lambda^{j_0})^{(\alpha,j_2)}$ does not hit $\gamma^*(\pd)$ and
$\lambda_{j_3}=\lambda_{j_0}.$ So, if $\rho$ is a path from
$u=\lambda_{j_0}$ to $\gamma^*(\pd)$ all the vertices of which  (possibly
with the exception of the last vertex) belong to $\Ztwo \setminus a,$ then Proposition \ref{fourth prop of patches} implies that 
 $|(\lambda^{j_0})^{(\alpha,j_2)},\rho)| \geq 14$ 
and combining this with  (\ref{*}) we get
$|\rho| \geq 8.$ 

Finally, let us consider the case $\tau(\pd,\omega_k) < \tau(\partial
A_{k-1}(\omega),\omega_k).$ Thus $k=V_i(\omega)$ for some index $1 \leq i
\leq 7.$ 
We want to show that $x_i(\pd)$ is the new point for the cluster $A_k(\varphi(\omega)).$ In order to do this we
first claim that $x_i(\pd)$ is the first vertex on
$\beta_i(\pd)$ which belongs to $\partial A_{k-1}(\f(\omega)).$ By construction
$x_i(\pd)$ is the first vertex on $\beta_i(\pd)$ which belongs to $\partial
(a \cup \{x_1(\pd), \ldots, x_{i-1}(\pd) \})$ and thus by property (1) of our
induction hypothesis the point $x_i(\pd)$ also belongs to $A_{k-1}(\f(\omega)) \cup
\partial A_{k-1}(\f(\omega)).$ Now let $v$ be a neighbor of a vertex $z$ of
$\beta_i(\pd),$ $z \neq x_i(\pd)$. Then by the above $v \not \in a$. By construction of $\beta_i(\pd)$ there exists a path $\tilde \rho$ in
$\Ztwo \setminus a$ with $|\tilde\rho| \leq 2$ connecting $z$ to
$\gamma^*(\pd).$Thus $v$ is connected to $\gamma^*(\pd)$  by a path in
$\Ztwo \setminus a$ of length at most 3. Hence by property (3) of our
induction hypothesis $v \not \in A_{k-1}(\f(\omega)) \setminus (\{x_1(\pd),
\ldots, x_{i-1}(\pd) \} \cup a)$ and our claim follows.

Next observe that each of the points $x_1(\pd), \ldots, x_7(\pd)$ satisfies
$d(x_i(\pd),a) \leq 7$ and each of the points of
$A_{V_7(\omega)}(\omega)$ has at least distance $40-13=27$ from $\Gamma.$ Thus by property (1) of our induction
hypothesis we can conclude that no neighbor of a lattice point
on $\hat \Gamma$ belongs to $A_{k-1}(\f(\omega)).$ Hence, if $\omega_{k,
\tau({\cal D})} \in \Gamma$ then $x_i(\pd)$ is the first vertex on $\alpha(i, \omega)$
which belongs to $\partial A_{k-1}(\f(\omega)).$
If $\omega_{k,\tau({\cal D})} \not \in \Gamma,$ assume without loss of
generality that $\omega_{k,\tau({\cal D})} \in \gamma^-(\pd).$ Let $l=\max \{j
\leq \tau(\pd,\omega_{k}): \omega_{kj} \in \Gamma\}.$ Then none of the vertices
on $(\omega_k^{\tau({\cal D})})^{(\alpha,l)}$ belongs to $\pd^* \cup \partial
A_{k-1}(\omega)$ and thus by property (1) of our induction hypothesis none
of these belongs to $\pd^* \cup \partial A_{k-1}(\f(\omega)).$ So, letting
$\Lambda$ denote that part of the loop $\hat \Gamma$ which starts at
$\omega_{kl}$ and ends at $\xi^l(\pd)$ it follows that
$(R((\omega_k^{\tau({\cal D})})^{(\alpha,l)}),\Lambda)$ belongs to ${\cal
R}_i(\omega).$ In particular, ${\cal R}_i(\omega)$ is nonempty and thus
$r(i,\omega)$ is well-defined. Moreover, using the definition
of ${\cal R}_i(\omega)$ and arguing as before it follows that $x_i(\pd)$ is the
first vertex on $\alpha(i, \omega)$ which belongs to $\partial
A_{k-1}(\f(\omega)).$ Thus $x_i(\pd)$ is in fact the new point for the cluster
$A_{k}(\f(\omega))$ and combining this with our induction
hypothesis, properties (1) - (3) follow. 

Finally note that by property (1) the definition of $\f(\omega)_k$ is
complete, i.e., cases (i)~-~(iii) exhaust all of the possible cases.
\newline
\vspace{7 mm}
\fine 

\section{Verification of Property \ref{second property of varphi} for the Map $\f$}

Let us now work on proving property \ref{second property of varphi} for the map $\varphi$, (the initial set $a$ and the integer $L$ being chosen as in the previous section).

\begin{proposition} \label{length of a V-loop}
There exists a constant $K$ independent of $a$ such
that for all $\omega \in \Xi$ and for all indices $k$ with $k=V_i(\omega)$
for some index $1 \leq i \leq 7$ the loop inserted into $\omega_k$ has
length at most $K.$
\end{proposition}
{\bf Proof.}  Let $\omega \in \Xi$ and let $k=V_i(\omega).$ Let
$\pd=\Theta(\omega).$ By Proposition \ref{first prop of patches} and our construction
of $\beta_i(\pd)$ we only need to consider the case $\omega_{k, \tau({\cal
D})} \in \gamma^-({\cal D}) \setminus \xi^l({\cal D})$ (resp. $\omega_{k, \tau({\cal
D})} \in \gamma^+({\cal D}) \setminus \xi^r({\cal D})$)  and it suffices to show that there exists a constant $K_1$
independent of $a$ and $\omega$ such that 
\[|r(i,\omega)| \leq K_1.\]
Without loss of generality assume that $x=\omega_{k, \tau({\cal
D})}$ belongs to  $\gamma^-=\gamma^-({\cal D}).$
It seems to be quite obvious that we can construct a path
from $x$ to $\xi^l(\pd)$ by following
$\gamma^-(\pd)$ or ``walking around'' $\partial A_{k-1}(\f(\omega))$ and
that the length of such a path is bounded above by a constant which is
independent of $a$ and $\omega.$ Nevertheless, in order to obtain a formal proof, we will follow a slightly different approach. 

We need to find a constant $K_1$ (independent of $a$ and $\omega$) and a path $\rho$ connecting $x$ and $\xi^l(\pd)$ such that $\rho$ does not hit $\partial A_{k-1}(\varphi(\omega))$ and such that $|\rho| \leq K_1.$ Now observe that $\omega_k$ does not hit $\partial A_{k-1}(\varphi(\omega))$ up to the hitting time of $\pd$ and since $\omega_k$ starts on $\Delta(a,L),$ it hits $\Gamma$ before hitting $\pd.$ Let $z$ denote the last point on $\omega_k^x$ which belongs to $\Gamma$ and let $\pi$ denote the non-selfintersecting path obtained from $(\omega_k^x)^{(\alpha,z)}$ by removing all loops. Let $\hat{\Gamma}(z,\xi_l)$ denote the lattice path obtained from $\hat{\Gamma}$ by starting at $z$ and then following the loop $\hat{\Gamma}$ until it hits the point $\xi_l$. Let  ${\cal M}_1$ denote the closure of the open
region bounded by  $R(\pi),$ $\hat{\Gamma}(z,\xi_l)$ and $(\gamma^-)^x$ and let $y^l$ denote the endpoint of $\gamma^-({\cal D})$.
 
We now claim that all points of ${\cal M}_1$ which belong to $A_{k-1}(\varphi(\omega)) \cup \partial A_{k-1}(\varphi(\omega))$ are within lattice distance 13 to $y^l.$ 
Once we have shown this claim we construct $\rho$ as follows: 

If $|x-y^l| > 13$ let $\rho=R(\gamma^-)^{(\alpha,x)}.$ (Note that in this case the length of $\rho$ is obviously bounded by the length of $\gamma$ and thus by $40$). 

If $ |x-y^l| \leq 13$ let us start at $x$ and let us follow $R(\pi)$ until we reach the first point $u$ which is at distance 14 to $y^l.$ Let us connect $u$ to $R(\gamma^-)_{14}$ by a non-selfintersecting path all the points of which are either at distance 14 or 15 to $y^l$, or they are within distance 13 to $y^l$ and belong to $\pi$ (see the formal definition below). Let us finally follow $R(\gamma^-)$ until we reach $\xi^l(\pd).$ 
 
More precisely, let ${\tilde \rho}$ denote the lattice loop
surrounding $\{u \in {\bf R}^2: |u-y^l| \leq 14 \}.$ 
 (${\tilde \rho}$ can be
constructed from the boundary of the above set in exactly the same way the
lattice loop ${\hat \Gamma}$ got constructed from the curve $\Gamma$). Let ${\cal M}_2$ denote the closure of
the bounded component of ${\bf R}^2 \setminus \tilde \rho$ and let
${\cal M}={\cal M}_1 \cap {\cal M}_2.$ Let $\hat{\rho}$ denote the lattice loop
surrounding ${\cal M}$,  
oriented such that
${\cal M}$ lies to the right of $\hat{\rho}.$ Thus the points $x$ and $w=R(\gamma^-)_{14}$ both lie on $\hat{\rho}.$
Define \[\rho=((\hat{\rho}^{(\alpha,x)})^w,
R(\gamma^-)^{(\alpha,w)}).\]
Then $\rho$ is a path from $x$ to $\xi^l(\pd)$ such that each
vertex $v$ on $\rho$ satisfies $v \in \pi$ (and thus $v \in \omega_k^{\tau({\cal D})}$) or $v \in {\cal M}_1$ and
$|v-y^l| \geq 14.$ Therefore $\rho$ does not hit $\partial
A_{k-1}(\varphi(\omega)).$ Moreover,
since $(\hat{\rho}^{(\alpha,x)})^w$ is a non-selfintersecting path, each vertex of
which belongs to ${\cal M}_2,$ we can conclude that $|(\hat{\rho}^{(\alpha,x)})^w|
\leq |{\cal M}_2 \cap \Ztwo|.$ Thus, $|\rho| \leq |(\hat{\rho}^{(\alpha,x)})^w| + |\gamma^-| \leq |{\cal M}_2 \cap \Ztwo| + 40.$ Hence, in both cases above, the constant $K_1=|{\cal M}_2 \cap \Ztwo| + 40$ is an upper bound for the length of $\rho.$ 

To finish the proof of the proposition, it remains to show that all points of ${\cal M}_1,$ which belong to $A_{k-1}(\varphi(\omega)) \cup \partial A_{k-1}(\varphi(\omega)),$ are within lattice distance 13 to $y^l.$ So let $u \in {\cal M}_1 \cap A_{k-1}(\varphi(\omega)).$ By part (1) of Proposition \ref{simple properties of varphi} we have $u \in A_{k-1}(\omega).$ Hence there exists a path $\lambda$ from $u$ to $a$ with $\lambda_i \in A_{k-1}(\omega) \setminus a$ for $i < |\lambda|.$ Yet, none of the points of ${\cal M}_1$ belongs to $a.$ (Here we use that $a$ is connected and none of the vertices of $(R(\pi)^z, 
{\hat \Gamma}(z,\xi_l),
(\gamma^-)^x)$ belongs to $a$). Thus $\lambda$ exits ${\cal M}_1$ and so it hits the lattice loop $(R(\pi)^z, 
{\hat \Gamma}(z,\xi_l),(\gamma^-)^x).$
 Moreover, since none of the points of $((R(\pi)^z, 
{\hat \Gamma}(z,\xi_l))$
 belongs to $A_{k-1}(\omega),$ we can conclude that $\lambda$ hits $(\gamma^-)^x.$ Let $v$ be such a point of intersection. Using that $(\gamma^-)^{(\alpha,v)}$ is a shortest path from $v$ to $a,$ we get $|\lambda| \geq |(\lambda^v, (\gamma^-)^{(\alpha,v)})| \geq |u-y^l|$ and thus by Proposition \ref{third prop of patches}
\[|u-y^l| \leq 12. \]
Now let $u \in {\cal M}_1 \cap \partial A_{k-1}(\varphi(\omega))$ and let $v$ be a neighbor of $u$ which belongs to $A_{k-1}(\varphi(\omega)).$ If $v \in {\cal M}_1$ we can conclude by the above that
\[|u-y^l| \leq 1+|v-y^l| \leq 13.\]
So assume that $v \not \in {\cal M}_1.$ Then $u \in ((R(\pi)^z, ({\hat \Gamma}^{\xi_l})^{(\alpha,z)},(\gamma^-)^x)$ and since none of the lattice points of $((R(\pi)^z, ({\hat \Gamma}^{\xi_l})^{(\alpha,z)})$ belongs to $A_{k-1}(\omega) \cup \partial A_{k-1}(\varphi(\omega))$ we have that $u \in (\gamma^-)^x.$ Thus $y^l$ is a point in $a$ which has minimal lattice distance to $u$ and so
\[ |u-y^l| =d(u,a) \leq 1+d(v,a). \]
Since $v \in A_{k-1}(\varphi(\omega))$ we have $v \in  A_{k-1}(\omega) \cup \{u_1(\pd),\ldots ,u_7(\pd)\}$ by part (1) of Proposition \ref{simple properties of varphi}. Now, if $v \in  A_{k-1}(\omega)$ then $d(v,a) \leq 12$ by Proposition \ref{third prop of patches} and if $v \in \{u_1(\pd),\ldots ,u_7(\pd)\}$ then $d(v,a) \leq d(v,y^*(\pd)) \leq 7$ by Remark \ref{rem:new positions}. Thus in either case the lattice distance between $v$ and $a$ is at most $12$ and combining this with the above inequality our claim follows. 
\newline
\vspace{7 mm}
\fine

\begin{proposition} \label{prop:catch-up loops}
i) Any catch-up loop has length at most 12. \newline
ii) For any $\omega \in \Xi$ at most 72 catch-up loops get inserted.
\end{proposition}
{\bf Proof.} i) Let $k$ be an index such that a catch-up loop $({\tilde
\lambda},R({\tilde \lambda}))$ with ${\tilde
\lambda}_0=\omega_{k,\tau(\partial A_{k-1}(\omega))}$ gets inserted
into $\omega_k.$ Then ${\tilde \lambda}_j \in A_{k-1}(\omega) \setminus
A_{k-1}(\f(\omega))$ for $j>0.$
Now 
$|A_{k-1}(\omega)| = |A_{k-1}(\f(\omega))| $ and
$A_{k-1}(\f(\omega)) \subset A_{k-1}(\omega) \cup \{u_1, \ldots, u_6\}.$
Thus $|{\tilde \lambda}| \leq |A_{k-1}(\omega) \setminus A_{k-1}(\f(\omega))|$ $ \leq 6.$ \newline
ii) Let ${\cal K}(\omega) = \{1 \leq k < V_7(\omega): \mbox{a
catch-up loop gets inserted into $\omega_k$} \}.$ For each $k \in {\cal K}(\omega)$ let
$u_k(\omega)$ denote the new point for the cluster $A_k(\f(\omega))$ and
let $z_i(\omega),$ $1 \leq i \leq 6,$ denote the new point for the cluster
$A_{V_i(\omega)}(\omega).$ Let ${\cal E}_i(\omega)$ denote the component of
$A_{V_7(\omega)-1}(\omega) \setminus a$ which contains $z_i(\omega).$ Note
that by Proposition \ref{third prop of patches}  $|{\cal E}_i(\omega)| \leq 12$ for $i=1, \ldots, 6.$ Now by part (2) of Proposition \ref{simple properties of varphi} for each $k \in {\cal K}(\omega)$ there exists an index $1 \leq i \leq 6$ with $u_k(\omega) \in {\cal E}_i(\omega).$ Hence we get
\[ |{\cal K}(\omega)| \leq |\bigcup_{i=1}^{6} {\cal E}_i(\omega)| \leq
\sum_{i=1}^{6} |{\cal E}_i(\omega)| \leq 6 \times 12 = 72. \]
\newline
\vspace{0 mm}
\fine

The next proposition will show that the map $\f$ is one-to-one. Let ${\hat
{\cal S}}=\{\rho: \rho \mbox{ is a finite path} $ in $\mbox{
$\Ztwo$ with $\rho_0 \in \Delta(a,L)$ and $\tau(a,\rho)=|\rho|$} \}.$
 
\begin{proposition} \label{prop:one-to-one}
If $(\rho_1, \ldots, \rho_L)$ , $({\tilde \rho}_1, \ldots, {\tilde \rho}_L) \in {\hat {\cal S}}^L$ with $(\rho_1, \ldots, \rho_L) \neq  ({\tilde \rho}_1, \ldots, {\tilde \rho}_L)$ then 
\[ \f(\omega_i^a = \rho_i \mbox{ , } i=1, \ldots, L) \cap  \f(\omega_i^a = {\tilde \rho}_i \mbox{ , } i=1, \ldots, L) = \emptyset. \]
\end{proposition}
{\bf Proof.} Suppose that $\omega^*$ lies in the intersection of the above
two sets. Then $\omega^*=\f(\omega)$ for some $\omega \in \Omega$ with
$\omega_i^a = \rho_i$ , $i=1, \ldots, L$ and
$\omega^*=\f({\tilde \omega})$ for some ${\tilde \omega} \in \Omega$ with
${\tilde \omega}_i^a = {\tilde \rho}_i$ , $i=1,
\ldots, L.$ Let $\pd=\Theta(\omega)$, ${\tilde {\cal D}}=\Theta({\tilde
\omega}),$ $n_7=V_7(\omega),$ ${\tilde n}_7=V_7({\tilde \omega}).$ 

Let us first consider the case $\pd \neq {\tilde {\cal D}}.$ Without loss of
generality assume that $n_7 \leq {\tilde n}_7.$ Since $\omega^*=\f(\omega)$
and $V_7(\omega)=n_7,$ Corollary \ref{second corollary of simple properties} applied to $\omega$ implies that at
least 7 points of $A_{n_7}(\omega^*)$ belong to the interior of $\pd.$ On
the other hand we also have $\omega^*=\f({\tilde \omega})$ and $V_7({\tilde
\omega})={\tilde n}_7.$ Thus, since $\Theta({\tilde \omega})= {\tilde {\cal D}}
\neq \pd,$ the same corollary applied to ${\tilde \omega}$ implies that at
most 6 points of $A_{{\tilde n}_7}(\omega^*)$ belong to the interior of $\pd$.
Noting that $A_{n_7}(\omega^*) \subset
A_{{\tilde n}_7}(\omega^*)$ since $n_7 \leq {\tilde n}_7,$ we arrive at a
contradiction. 

Let us now consider the case $\pd={\tilde {\cal D}}.$ Let $i_0=\min\{i:\rho_i
\neq {\tilde \rho}_i \}$ and $j_0=\min\{j:\rho_{i_0,j} \neq {\tilde
\rho}_{i_0,j} \}.$ Since for any $1 \leq i \leq L,$ $\tau(\partial A_{i-1}(\omega),
\omega_i) \leq \tau(\partial a, \omega_i)$ we have $A_i(\omega)={\cal
C}(\rho_1, \ldots, \rho_i)$ and similarly $A_i({\tilde \omega})={\cal
C}({\tilde \rho}_1, \ldots, {\tilde \rho}_i).$  Now by definition of $i_0,$ 
$(\rho_1, \ldots, \rho_{i_0-1})=({\tilde \rho}_1, \ldots, {\tilde
\rho}_{i_0-1})$ and hence
\begin{equation} \label{eqn 1 in one-to-one}
A_{i_0-1}(\omega)=A_{i_0-1}({\tilde \omega}).
\end{equation}
Moreover, by choice of $\omega$ and ${\tilde \omega}$ we have
$\f(\omega)=\omega^*=\f({\tilde \omega})$ and so, in particular,   
\begin{equation} \label{eqn 2 in one-to-one}
A_{i_0-1}(\f(\omega))=A_{i_0-1}(\f({\tilde \omega})).
\end{equation}
Thus, using that $\Theta(\omega)=\Theta({\tilde \omega}),$ we apply exactly
the same procedure to obtain $\f(\omega)_{i_0}$ from $\omega_{i_0}$ as the
one we use to obtain  $\f({\tilde \omega})_{i_0}$ from ${\tilde
\omega}_{i_0}.$ Loosely speaking, this rule is a step-by-step rule and
hence we pick up the first step where $\rho_{i_0}$ and ${\tilde
\rho}_{i_0}$ differ. More formally, let us observe that for $j \leq j_0-1$
a loop gets inserted into $\omega_{i_0}$ at time $j$ (in order to form
$\f(\omega)_{i_0}$ ) if and only if this loop gets inserted into ${\tilde
\rho}_{i_0}$ at time $j_0$ and vice versa. This follows by definition of
$\f$ from (\ref{eqn 1 in one-to-one}), (\ref{eqn 2 in one-to-one}) and the fact that $\omega_i^j={\tilde \omega}_i^j$ for
any $j \leq j_0-1.$ Thus, if no loop gets inserted into $\omega_{i_0}$ up
to (including) time $j_0-1$ then no loop gets inserted into ${\tilde
\omega}_{i_0}$ up to (including) time $j_0-1$ and hence 
\[ \f(\omega)_{i_0j_0}=\omega_{i_0j_0} \]
and
\[ \f({\tilde \omega})_{i_0j_0}={\tilde \omega}_{i_0j_0}. \]
Yet, the two left hand sides agree since $\f(\omega)=\omega^*=\f({\tilde
\omega}),$ while the two right hand sides are different by choice of $i_0$
and $j_0.$
Thus we arrive at a contradiction. On the other hand, if a loop gets
inserted into $\omega_{i_0}$ at time $j \leq j_0-1,$ then 
\[ \f(\omega)_{i_0}=(\omega_{i_0}^j, \lambda, R(\lambda),\omega_{i_0}^{(\alpha,j)}) \]
and
\[ \f({\tilde \omega})_{i_0}=({\tilde \omega}_{i_0}^j, \lambda,
R(\lambda),{\tilde \omega}_{i_0}^{(\alpha,j)}) \]
for some path $\lambda$ with $\lambda_0=\omega_{i_0j}={\tilde \omega}_{i_0j}.$
Yet, $\omega_{i_0}^j={\tilde \omega}_{i_0}^j$ and $\omega_{i_0}^{(\alpha,j)}
\neq {\tilde \omega}_{i_0}^{(\alpha,j)}$ by choice of $i_0$ and $j_0.$ Thus,
once again, the two left hand sides agree while the two right hand sides
are different and we arrive at a contradiction.  
\newline
\vspace{7 mm}
\fine

We are now ready to prove property \ref{second property of varphi} of the map $\f.$ As before let $N_0$ be an integer satisfying the conclusion of Remark \ref{rem:size}.
\begin{lemma}
There exists a constant $c>0$ such that for any finite cluster $a$ with $|a| \geq N_0$ 
\[P(\f(\Omega^{a,L})) \geq c,\]
where $L=6|a|+1.$
\end{lemma}
{\bf Proof.} Let ${\hat {\cal S}}$ 
be as previously defined.
For paths $\rho_1, \ldots, \rho_L$ in ${\hat {\cal S}}$ define the event
${\cal G}_{(\rho_1, \ldots, \rho_L)}=\{\omega_1^a=\rho_1, \ldots,
\omega_L^a=\rho_L \}.$ Note that by definition of $\f$ for any
paths $\rho_1, \ldots \rho_L \in {\hat {\cal S}}$ there
exist paths ${\tilde \rho}_1, \ldots, {\tilde \rho_L} \in {\hat {\cal S}}$ with
$\f({\cal G}_{(\rho_1, \ldots, \rho_L)})={\cal G}_{({\tilde \rho}_1, \ldots,
{\tilde \rho}_L)}.$ Moreover, ${\tilde \rho}_k$ is different from $\rho_k$
if and only if 
$k=V_i(\rho_1, \ldots, \rho_L)$ for some index $1 \leq i \leq 7$ or if
the path $\rho_k$ hits the boundary of the cluster built by the paths $\rho_1, \ldots, \rho_{k-1}$ strictly before it hits the boundary of the cluster built by the paths ${\tilde \rho}_1, \ldots, {\tilde \rho}_{k-1}$ and if 
the patch $\Theta(\rho_1, \ldots, \rho_{L})$ has not been hit prior to this time.     
In this last case a catch-up loop is inserted into $\rho_k.$
Thus by the second part of Proposition \ref{prop:catch-up loops}, ${\tilde \rho}_k \neq \rho_k$ for at most $72+7=79$ of these
paths. Now by Proposition \ref{length of a V-loop} for $k=V_i(\rho_1, \ldots, \rho_L)$ the
length of the loop inserted into $\rho_k$ is bounded above by a constant
$K$ which is independent of $a,L$ and $\rho_1, \ldots, \rho_L$ and by the first part of 
Proposition \ref{prop:catch-up loops} the length of a catch-up loop is bounded above by 12. Thus
$\sum_{k=1}^L |{\tilde \rho_k}| \leq \sum_{k=1}^L |{\rho_k}| + 79 K_1,$ where 
$K_1=K \vee 12.$ Letting $c=4^{-79 K_1}$ we get
\begin{eqnarray*}
P(\f({\cal G}_{(\rho_1, \ldots, \rho_L)})) 
  & = & P({\cal G}_{({ \tilde \rho}_1, \ldots, {\tilde \rho}_L)}) \\
  & = & \left( \prod_{k=1}^L \mu_{\Delta(a,L)}({\tilde \rho}_{k0}) \right)
4^{-\sum_{k=1}^L |{\tilde \rho_k}|} \\
  & \geq & c \left( \prod_{k=1}^L \mu_{\Delta(a,L)}({\rho}_{k0}) \right)
4^{-\sum_{k=1}^L |{\rho_k}|} \\ 
  & = & c \, P({\cal G}_{(\rho_1, \ldots, \rho_L)}). 
\end{eqnarray*}

Combining this last inequality with Proposition \ref{prop:one-to-one} we finally get
\begin{eqnarray*}
P(\f(\Omega^{a,L})) 
  & = & P(\f(\Xi)) \\
  & = & P \left(\f( \bigcup_{(\rho_1, \ldots, \rho_L) \in {\hat {\cal S}}^L} {\cal
G}_{(\rho_1, \ldots \rho_L)}) \right) \\
  & = & \sum_{(\rho_1, \ldots, \rho_L) \in {\hat {\cal S}}^L} P \left( \f({\cal
G}_{(\rho_1, \ldots \rho_L)}) \right) \\
  & \geq & c \sum_{(\rho_1, \ldots, \rho_L) \in {\hat {\cal S}}^L} P \left( {\cal
G}_{(\rho_1, \ldots \rho_L)} \right) \\
  & = & c \, P(\Xi) \\
  & = & c.
\end{eqnarray*}

\vspace{0 mm}
\fine

\chapter*{Glossary for Chapter 2}      
\addcontentsline{toc}{chapter}{Glossary for Chapter 2}

\thispagestyle{myheadings}{}{}
\begin{description}
\item[$K$:] number of additions in a period, $K \in {\bf N},$ $ K \geq 2$.
\item[$\kappa$:] number of steps in a period, $\kappa=2K+1$. 
\item[$h$:] the function $h(n)= \log n /{\log { \log n}}$. 
\item[$f$:] the function describing the transition from counting in the form 
$0$, $(1,1)$,$(1,2)$,$ \ldots$, $(1,\kappa)$, $(2,1)$, $(2,2)$, $\ldots$,$ (2,\kappa)$,$ \ldots$  to counting
in the form  $0,1,2, \ldots$. 
\item[$A_n$:] the (random) set of sites after $n$ periods are completed .
\item[$A_{n,j}$:] the (random) set of sites after step $j$ in period $n$ is completed .
\item[$X_{n,j}$:] the site added to (resp. deleted from) the set $A_{n,j-1}$ in step $j$ of period $n$.
\end{description}

\vspace{3.5 mm}
 
In the following $a$ is a finite subset of ${\bf Z}$.
  
\begin{description}
\item[$G(a)$:] $\sharp$ of gaps of $a$. 
\item[$G^{(2)}(a)$:] $\sharp$ of gaps of $a$ of size at least 2 .
\item[$L(a)$:] $\sharp$ of gap sites of $a$].
\item[$C(a)$:] set of gap sites of $a$ which are within distance 4 to some other gap site of $a$.
\item[$C^{(2)}(a)$:] set of gap sites of $a$ which are within distance 4 to
at least 2 other gap sites of $a$ 
\item[$\partial a$:] boundary of $a$, i.e., 
$\{y \in {\bf Z} \setminus a: \exists x \in a \mbox{ with } |y-x| =1 \}$.
\item[$\bar{\partial} a$:] outer boundary of $a$, i.e., $\{\min a-1, \max a +1 \}$.
\end{description}

$G_n,$ $G_n^{(2)},$ $L_n$ etc. stand for $G(A_n),$ $G^{(2)}(A_n)$ $L(A_n)$ etc.

\vspace{6 mm}

$(a, \lambda)$ - new and old gaps of $a$: \newline 
The gaps of $a$ are enumerated `from left to right' as $g_1, \ldots ,g_{G(a)}$. $\lambda$ is a subset of the index set $\{1, \ldots, G(a)\}$. The gap $g_i$ of $a$ is called {\it new} if $i \in \lambda,$ otherwise {\it old}.

\vspace{6 mm}

$G(a,\lambda),$  $G^{(2)}(a,\lambda)$ etc. defined as $G(a),$ $G^{(2)}(a)$ etc. above with `gap' replaced by `new gap' at each occurrence, thus, for instance, $C^{(2)}(a,\lambda)= \sharp $ of new gap sites  of $a$ which are within distance 4 to at least 2 other new gap sites of $a$.

\begin{description}
\item[$\partial^o(a,\lambda)$:] subset of $\partial a$ consisting of all boundary points  belonging to old gaps. 
\item[$\partial^n(a,\lambda)$:] subset of $\partial a$ consisting of all boundary points  belonging to new gaps .
\item[$V^o(a,\lambda)$:] subset of $a$ consisting of all points which are adjacent to an old gap 
but not an endpoint of $a$. 
\item[$V^n(a,\lambda)$:] subset of $a$ consisting of all points which are within distance 4 to a new gap
site but which are neither adjacent to an old gap nor an endpoint of $a$. 
\item[$I_n$:] a (random) subset of $\{1, \ldots, G(A_n)\}$ describing the index set of new gaps of $A_n$.
\end{description}

$G^I_n,$ $G_n^{(2),I},$ $L^I_n$ etc. stand for $G(A_n,I_n),$ $G^{(2)}(A_n,I_n),$ $L(A_n,I_n)$ etc.

\begin{description}
\item[$\Gamma_m$:] $\sharp$ of new gap sites which are within distance 4 to some other new gap site
after a total of $m$ steps is completed, $\Gamma_m=|C^I|(A_{f^{-1}(m)},I_{f^{-1}(m)})$.
\item[$M$:] a positive integer .
\item[${\cal G}$:] $\{ a \subset {\bf Z}: |a| < \infty, G(a) \leq 10K \mbox{ and
} L(a) \leq (K-1) |a| \}$ .
\item[${\cal S}_M$:] $\{ a \subset {\bf Z} : |a| < \infty, L(a) \leq h(|a|)/M,
|C(a)| \leq 3 \mbox{ and } |C^{(2)}(a)| \leq 1 \} $ . \\
\item[$\tau^1$:] $\inf \{n: L^I_n \geq (K^{-1}(K-1) + 1)h(|A_0|) \},$ i.e., $\tau^1$ denotes the first period such that  
at the end of this period the number of new gap sites  exceeds 
  $K^{-1}(K-1) h(|A_0|)$
 by at least $h(|A_0|)$ . 
\item[$T^{M,2}$:] $\inf \{n: L_n \geq (K^{-1}(K-1)+ 3 M^{-1} )h(|A_0|) \} ,$ i.e., $T^{M,2}$ denotes the first period 
such that 
 at the end of this period the number of gap sites exceeds $K^{-1}(K-1) h(|A_0|)$ 
 by at least $3 M^{-1} h(|A_0|)$.
\item[$\tau^{M,2}$:] $\inf \{n: L^I_n \geq (K^{-1}(K-1)+ 3 M^{-1} )h(|A_0|) \} ,$ i.e., $\tau^{M,2}$ denotes the first period such that 
 at the end of this period the number of new gap sites exceeds  $K^{-1}(K-1) h(|A_0|)$ 
 by at least $3  M^{-1} h(|A_0|)$. 
\item[$\tau_M$:] $\inf \{ n \geq |A_0|^{11/10}: L^I_{n} \leq M^{-1} h(|A_0|),$ i.e., $\tau_M$ denotes the first
period after period $|A_0|^{11/10}$
 such that 
at the end of this period the
number of new gap sites is `small' compared 
 to $h(|A_0|)$ . 
\item[$ T_{M,2}$:] $\inf \{ n \geq |A_0|: L(A_n) \leq M^{-1} h(|A_0|) \},$ i.e., 
$T_{M,2}$ denotes the first
period after period  $|A_0|$ 
 such that 
at the end of this period the
number of gap sites is `small' compared to $h(|A_0|)$. 
\item[$\tau_{M,2}$:] $\inf \{ n \geq |A_0|: L^I(A_n,I_n) \leq h(|A_0|)/M
\},$ i.e.,  $\tau_{M,2}$ denotes the first
period after period  $|A_0|$  such that 
at the end of this period the
number of new gap sites is `small' compared to $h(|A_0|)$.   
\item[$S$:] $\inf \{(n,i): i>K, |C_{n,i-1}| \geq 2 \mbox{ and } X_{n,i} \in
C_{n,i} \},$ i.e., $S$ denotes  
the first deletion step in which the following conditions hold: \\ 
(i) $X_{n,i}$ is a gap site for $A_{n,i},$  \\
(ii) $X_{n,i}$ is within distance 4 to some other gap site for $A_{n,i},$ \\  
(iii) $A_{n,i-1}$ had at least 2 gap sites which are within distance 4 to each other. 
\item[$\bar{S}$:] the period containing $S$.
\item[$\sigma$:] $\inf \{(n,i): i>K, |C_{n,i-1}^I| \geq 2 \mbox{ and } X_{n,i} \in
C_{n,i}^I \}, $ i.e., $\sigma$ is defined as $S$ with  `gap sites' replaced by `new gap sites' at each occurrence.
\item[$\bar{\sigma}$:] the period containing $\sigma$ .
\item[$\sigma_0$:] $\inf \{n: n \geq |A_0|^{11/10} \mbox{ and } G_{n}^{I,2}=0
\} ,$ i.e., 
$\sigma_0$ denotes the first period after time $|A_0|^{11/10}$  such
that at the end of this period all new gaps have size 1.
\item[$\theta_0$:] $\inf\{n:L^I_n=L_n\},$
 i.e., $\theta_0$ denotes the first period such
that at the end of this period there are no more old gaps left.
\item[${\tilde \zeta}$:] $\sigma \wedge (\tau^1 ,\kappa)$.
\item[$\zeta$:] $f({\tilde \zeta})$.
\item[$U$:] a process of independent random variables which are uniformly distributed 
on $(0,1),$ 
 in short, an  i.u. process$^*$. 
\item[$U^{(1)}, \ldots, U^{(3)}$:] i.u. processes$^*$ satisfying that  $U,U^{(1)},U^{(2)}$ and $U^{(3)}$ are  independent .
\item[$U^{(4)}$:] an i.u. process$^*$ which is independent of $U$ ,
($U^{(4)}$ is only used in the proof of Proposition 2.10.2). 
\item[${\cal N}$:] (p. 33,34) $\sharp$ of excursions of $L^I$ which start strictly before 
period $\bar{\sigma}$ .
\end{description}

\vspace{3 mm}

For $i \leq {\cal N}(\omega)$ the following is defined (p. 33,34).

\begin{description}
\item[$f_i(\omega)$:] start of the $i$-th excursion of $L^I(\omega).$ 
\item[$\rho_i(\omega)$:] length of the $i$-th excursion of $L^I(\omega)$. 
\item[$e^i(\omega)$:] the $i$-th excursion of $L^I(\omega),$ $e^i(\omega)=(L^I_{f_{i-1}+j})_{j=0, \ldots, \rho_i}(\omega)$ . \\
\item[$U^L$:] (p. 34,35) an i.u. process$^*$ constructed pathwise by pasting together parts of $U(\omega)$ and of $U^{(1)}(\omega)$. 
\item[$\hat{L}^I$:] (p. 34) a random walk constructed from $U^L$. 
\item[$\hat{f}_i$:] end of the $i$-th successive $K$-decrease for the
process $\hat{L}^I,$ i.e., $\hat{f}_0  = 0,$ $\hat{f}_n(\omega)  = \inf \{j > \hat{f}_{n-1}(\omega):$  $\hat{L}_j(\omega)
-  \hat{L}_{\hat{f}_{n-1}(\omega)}(\omega) \leq -K \}$. 
\item[$\hat{\rho}_i$:] length of the $i$-th successive $K$-decrease for the
process $\hat{L}^I,$ i.e., $\hat{\rho}_i=\hat{f}_i-\hat{f}_{i-1}$.
\item[$V_i$:] $V_0(\omega)=\inf \{m: \Gamma_m(\omega) > 0\},$ i.e., $V_0$ denotes the first step at which there exist two new gap sites which are within distance 4 to each other, \\
 $i \geq 1$: $V_i(\omega) = \inf \{m>V_{i-1}(\omega): \Gamma_m(\omega) >
\Gamma_{m-1}(\omega) \},$i.e.,  $V_i$ denotes the $i$-th point of increase for $\Gamma($ after time $V_0$.
\item[${\cal M}$:] $\sharp$ of
points of increase of $\Gamma$ in the time interval $(V_0,\alpha),$ i.e.,   ${\cal M}(\omega)=\sup \{i:V_i(\omega) < \alpha(\omega) \},$  ($ \sup \emptyset =0$) .
\item[$W^{(1)}_i$:] the first time of decrease of $\Gamma$ after time $V_{i},$ i.e., $W^{(1)}_i(\omega)= \inf\{m > V_i(\omega): \Gamma_m(\omega)< \Gamma_{m-1}(\omega)\}$  .
\item[$W^{(2)}_i$:] the first time after $V_{i}$ at which $\Gamma$ equals 0, i.e., $W^{(2)}_i(\omega)= \inf\{m > V_i(\omega): \Gamma_m(\omega)=0 \}$.  
\item[$U^X$:] (p. 44,45) an i.u. process$^*$ constructed pathwise by pasting together parts of $U(\omega)$ and of $U^{(2)}(\omega)$. 
\item[$V_i^X$:] the $i-$th ``d-success" of $U^X$, i.e., $V_0^X=0$ and  for $i \geq 1,$ $V_i^X(\omega) = \inf \{n > V^X_{i-1}(\omega):U^X_{n}(\omega) \leq \ln N/N \}$. 
\item[$U^Y$:] (p. 46-48) an i.u. process$^*$ constructed pathwise by pasting together parts of $U(\omega)$ and of $U^{(3)}(\omega)$.
\item[$V_i^Y$:] the $i-$th a-success of $U^Y$ , i.e., $V_0^Y=0$ and for $i \geq 1,$ $V_i^Y(\omega) = \inf \{n > V^Y_{i-1}(\omega):  U^Y_{n}(\omega)  > 1-1/{\ln N} \}$ .
\item[$\alpha_i$:] $\alpha_0=0$ and for $i \geq 1,$ $\alpha_i(\omega)=\inf \{n \geq \alpha_{i-1}(\omega)+C: L_n(\omega) \leq C \},$  where $C$ is a fixed integer with $C \geq 10K+1$. 
\item[$\beta_i$:] $\alpha_i-\alpha_{i-1}-C$.
\item[$U^*$:] (p. 63) an i.u. process$^*$ constructed pathwise by pasting together parts of $U(\omega)$ and of $U^{(4)}(\omega)$. 
\item[$L^*$:] (p. 63) a random walk constructed from $U^*$.
\item[$f_n^*$:] end of the $n$-th CK-drop for $L^*$, i.e., $f^*_0(\omega)=0$ and for $n \geq 1,$ $f_n^*(\omega)=\inf \{j>f_{n-1}^*(\omega):(L^*_j-L^*_{f_{n-1}^*})(\omega) \leq -CK \}$.
\item[$\beta_n^*$:] length of the $n$-th CK-drop for $L^*,$ i.e.,  $\beta_n^*=f_n^*-f_{n-1}^*$.

\end{description}
\chapter*{Glossary for Chapter 3}      
\addcontentsline{toc}{chapter}{Glossary for Chapter 3}

\thispagestyle{myheadings}{}{}
\begin{description}
\item[{$\partial a$}:] the boundary of $a$, ($a \subset \Ztwo$), i.e., $x \in \partial a$  if $x$ is a neighbor of a point in $a$, 
but $x \not \in a$.
\item[$L$:] an integer usually chosen depending on the given initial set $a$;
in sections 3.6 and 3.7 and in most of section 3.4, $L$ is chosen as  
$L=6|a|+1$.
\item[$H(a)$:] the number of holes of the set $a$.
\item[$\mu_{b}$:] harmonic measure on $b$, $b \subset \Ztwo$. 
\item[$|x|$:] the lattice norm of $x$, $|x|=|x_1|+|x_2|$; for $x,y \in \Ztwo$, $|x-y|$ is called the 
lattice distance between $x$ and $y$. 
\item[$d(x,b)$:] $\min \{|x-y|: y \in b\}$, the lattice distance between $x \in \Ztwo$ and $b \subset \Ztwo$, ($b \neq \emptyset$).  
\item[$S(a)$:] $\sup\{|x|: x \in a \}$. 
\item[$\Delta(a,L)$:]  $\{x \in {\bf Z}^2 : |x|=S(a)+L \}$.  
\item[$|\rho|$:] (p. 70) the length of the path $\rho$ . 
\item[$\tau(b, \rho)$:] the hitting time of the set $b$ for the path $\rho$, i.e.,
$\tau(b,\rho)=\inf\{0 \leq i < |\rho|+1: \rho(i) \in
b \},$ ($\inf \emptyset = \infty$).
\item[$\tau(x, \rho)$:] $\tau(\{x\},\rho)$, $x \in \Ztwo$ .
\item[$\rho_{\tau(b)}$:] $\rho_{\tau(b,\rho)}$, $b \subset \Ztwo$,  
(assuming that $\tau(b,\rho) < \infty$);   
similarly $\rho_{\tau(x)}=\rho_{\tau(x,\rho)}$, $x \in \Ztwo$, (assuming that $\tau(x,\rho) < \infty$).
\item[${\cal S}$:] $\{\rho: \rho \mbox{ is an infinte path in
$\Ztwo$ with $\rho_0 \in \Delta(a,L)$} \}$ 
 ($a$ being the given initial set).
\item[${\hat {\cal S}}$:] $\{\rho: \rho \mbox{ is a finite path in $\Ztwo$ with $\rho_0 \in \Delta(a,L)$ and $\tau(a,\rho)=|\rho|$} 
\}$ .
\item[$\Omega^{a,L}$:] $\prod_{i=1}^L {\cal S}$ . 
\item[$\Phi_i$:] the projection onto
the $i$-th coordinate path, i.e., $\Phi_i(\omega)=\omega_i,$ $\omega \in \Omega^{a,L}$ .
\item[$\Phi_{ij}$:] $\Phi_{ij}(\omega)=\omega_{ij}$. 
\item[$\Phi_i^k$:]  
$\Phi_i^k(\omega)=(\omega_{ij})_{0 \leq j \leq k}$.
\item[${\cal F}^{a,L}$:] the $\sigma$-field in $\Omega$ generated by the maps $\Phi_{ij},$ $1 \leq i \leq L,$ $j \in \Nzero$ .
\item[$P^{a,L}$:] the probability measure  on $(\Omega^{a,L}, {\cal F}^{a,L})$ defined by 
$P(\Phi_i^{l_i}=\rho_i, i=1, \ldots ,L) = \prod_{i=1}^L
\left(\mu_{\Delta}(\rho_{i0}) 4^{-l_i}\right) $ 
for any finite paths $\rho_1, \ldots ,\rho_L$ with 
$\rho_{i0} \in \Delta(a,L)$
and $|\rho_i|=l_i,$ $i=1, \ldots ,L$ .
\item[${\cal C}(\rho_1, \ldots, \rho_n)$:] (p. 71) 
rougly speaking, ${\cal C}(\rho_1, \ldots, \rho_n)$ is the cluster obtained from the 
initial set $a$  
after $n$ particles which start at the points $\rho_{i0} 
\in \Delta(a,L)$,$i=1, \ldots, n$,
( $\rho_{i0} \not\in a \cup {\partial a}$)  
and follow 
the paths $\rho_i$, $i=1, \ldots, n,$ are attached .
\item[$A^{a,L}$:] the Markov process on $(\Omega^{a,L},{\cal F}^{a,L},P^{a,L})$ defined by  $A_n^{a,L}(\omega)={\cal C}(\omega_1, \ldots, \omega_n)$;  
$A_n^{a,L}$ describes the random cluster obtained from 
the initial set $a$ after $n$ new 
particles have been added .
\item[$\varphi$:] (Section 3.6) the map $\varphi: \Omega^{a,L} \rightarrow \Omega^{a,L}$ is constructed by inserting additional loops into  
some of the given paths $\omega_1, \ldots , \omega_L$; the changed paths $\varphi(\omega)_1,
\ldots ,
\varphi(\omega)_L$ build a new cluster $A_L(\varphi(\omega))$ which has at least one more hole than the cluster $a$ we start with.
\item[$u \la v$:] the line segment with endpoints $u,v$;  
a point $z \in u \la v$ with $z \not \in \{u,v\}$ is called an {\it inner point} of $u \longleftrightarrow v$;   
an {\it edge} of $a$ is a line segment  $u \la v$ with $u,v \in a$ and $|u-v|=1$.
\item[$G(a)$:] the union over all edges of $a$.
\item[$\Gamma$:] (p. 74-76) an oriented curve surrounding $a$; the set $a$ lies to the right of $\Gamma$ and each lattice point on $\Gamma$ has lattice distance 40 from $a$ . 
\item[$\gamma(\xi)$:] a lattice path of length 40 from $\xi$, ($\xi$ a lattice point on $\Gamma$), to $a$;  these  paths are chosen such that whenever two paths touch they stay together .
\item[$y(\xi)$:] the endpoint of $\gamma(\xi)$.
\item[$\xi_0$:] $\xi_0 \in \Gamma$ is chosen such that $\xi_0-40e_2$ is a point of $a$ with maximal second component, thus $\gamma(\xi_0)=(\xi_0, \xi_0-e_2, \ldots, \xi_0-40e_2)$ .
\item[$\Gamma(\xi, {\tilde \xi})$:] the part of $\Gamma$ which connects $\xi$ to $\tilde \xi$ .
\item[$I$:] (p. 78) an integer which is chosen in the construction of the sequence $u$ below  .
\item[$u$:] (p. 78) $(u_i)_{0 \leq i \leq I}$ is a sequence of lattice points on $\Gamma$ satisfying $u_I=u_0=\xi_0$, these points being chosen such that $\gamma(u_0), \gamma(u_1) \ldots \gamma(u_{I-1})$ are nonintersecting and, loosely speaking, neither too close nor too far from each other .
\item[$u^I$:] (p. 78) a lattice point on $\Gamma(u_{I-1},u_0)$ which is chosen when constructing the sequence $u$ above .
\item[$v$:] (p. 78) $(v_i)_{1 \leq i \leq I}$ is a sequence of lattice points on $\Gamma$ with $v_i \in \Gamma(u_{i-1},u_i)$, which is chosen such that any path from $y(v_i)$ to $\gamma(u_{i-1})$ or $\gamma(u_i)$ which intersects $a$ at most in its two endpoints has length at least 40 . 
\item[${\cal D}_i$:] (p. 78) the $i$-th {\it patch} of $a$, defined as the closure of the open region bounded by $\gamma(u_{i-1}),$ $\Gamma(u_{i-1},u_i),$ $\gamma(u_{i})$ and those edges of $a$ whose interior points can be connected by a curve in ${\bf R}^2$ to $\Gamma(u_{i-1},u_i)$ without hitting any other point of $\Gamma(u_{i-1},u_i),$ $\gamma(u_{i-1}),$
$\gamma(u_{i})$ or any other edge of $a$ .
\end{description}

\noindent For the following eight entries  ${\cal D}$ is a patch of $a$ and $i$ is 
chosen such that ${\cal D}={\cal D}_i$. 

\begin{description}
\item[$\xi^l({\cal D})$:] $u_{i-1}$. 
\item[$\xi^m({\cal D})$:] $v_{i}$ . 
\item[$\xi^r({\cal D})$:] $u_{i}$.

\item[$\gamma^-({\cal D})$:] $\gamma(\xi^l({\cal D}))$ .
\item[$\gamma^*({\cal D})$:] $\gamma(\xi^m({\cal D}))$ .
\item[$\gamma^+({\cal D})$:] $\gamma(\xi^r({\cal D}))$ .
\item[$y^*({\cal D})$:] the endpoint of $\gamma^*({\cal D})$. 
\item[${\cal D}^*$:] the set of all lattice points which belong to the interior of ${\cal D}$ .   
\end{description}

\vspace{7 mm}

\noindent $\mbox{Pol}(y_{i})$: (p. 74) the closed polygon with corners $y_i \pm 40 e_1$ and $y_i \pm 40 e_2$  .

\vspace{7 mm}

\noindent $\mbox{Con}(y_{i})$: the contour of $\mbox{Pol}(y_{i})$.

\vspace{7 mm}

\noindent $\mbox{Pol}(y_1, \ldots ,y_{i})$: (p.74-76) 
  a closed polygon which
contains $\, \mbox{Pol}(y_{j})$ for $1 \leq j \leq i$ and 

which satisfies that any lattice point on its contour belongs to $\mbox{Con}(y_{j})$ for some

index $1 \leq j \leq i$.

\begin{description}
\item[$\rho^i$:] the truncated path $(\rho_j)_{0 \leq j \leq i}$.
\item[$\rho^{(\alpha,i)}$:] the path obtained from $\rho$ by removing $\rho^{i-1}$,  i.e., $\rho^{(\alpha,i)}_j=\rho_{i+j}$, $0 \leq j \leq |\rho|-i$. 
\item[$\rho^b$:] $\rho^{\tau(b,\rho)}$,  i.e., $\rho$ is stopped when hitting the set $b$ (assuming that $\tau(b,\rho) < \infty$); for $z \in \Ztwo$ let  $\rho^z=\rho^{\{z\}}$.
\item[$\rho^{(\alpha,b)}$:]  $\rho^{(\alpha,\tau(b,\rho))}$; for $z \in \Ztwo$ let  $\rho^{(\alpha,z)}=\rho^{ (\alpha,\{z\}) }$.
\item[$(z,\rho)$:] the path obtained from $\rho$ by including $z$ as new starting point, i.e., $(z,\rho)_0=z$ and $(z,\rho)_i=\rho_{i-1}$ for $i<|\rho|+2$ (assuming that $z$ is a neighbor of $\rho_0$) .
\item[$(\rho,z)$:] the path obtained from $\rho$ by appending the point $z$, i.e.,  $(\rho,z)=(\rho_0, \ldots, \rho_{|\rho|},z)$ 
(assuming that $\rho$ is finite and $z$ is a neighbor of $\rho_{|\rho|}$).
\item[$(\rho,\pi)$:] (p. 76,77) the path obtained by appending the path $\pi$ to $\rho$ (assuming that $\rho$ is  
finite and  $\rho_{|\rho|}$ and $\pi_0$ are equal or neighbors of each other) .
\item[$R(\rho)$:] 
the path obtained by reversing $\rho$, i.e., $R(\rho)_i=\rho_{|\rho|-i}$, $1 \leq i \leq |\rho|$.
\item[$N_{j,{\cal D}}(\rho_1, \ldots, \rho_L)$:] the number of paths among $\rho_1, \ldots, \rho_j$, ($j \leq L$), which hit the patch $\cal D$ (strictly) before hitting the growing cluster, i.e., $N_{j,{\cal D}}(\rho_1, \ldots, \rho_L) =  
\sum_{k=1}^{j} {\bf 1}_{
\{\tau({\cal D},\rho_k) < \tau(\partial {\cal C}(\rho_1, \ldots,
\rho_{k-1}),\rho_k)\}}$.
\item[$\Theta(\rho_1, \ldots, \rho_L)$:] (p. 80,81) the ``lucky" patch, i.e., the first patch to be hit by 7 paths before 
these hit the boundary of the growing cluster  .
\item[$V_i(\rho_1, \ldots, \rho_L)$:] (p. 80,81) the index of the $i$-th path , ($1 \leq i \leq 7$), which hits $\Theta(\rho_1,
\ldots, \rho_L)$ before hitting the boundary of the growing cluster.
\item[$x_i({\cal D})$:] (p. 83-85) For each patch $\cal D$ of $a$, the lattice points $x_i({\cal D})$, $(1 \leq i \leq 7)$, are chosen such that the point $y^*({\cal D})$ (which by definition 
 belongs to the unbounded component of $\Ztwo \setminus a$) belongs to a hole of $a \cup \{x_1({\cal D}), \ldots, x_7({\cal D})\}$ . 
\item[$\beta_i({\cal D})$:] a path from $\xi^m({\cal D})$ to $x_i({\cal D})$; the paths  $\beta_1({\cal D}), \ldots ,\beta_7({\cal D})$ satisfy that ${\cal C}(\beta_1, \ldots, \beta_j)=a \cup \{x_1({\cal D}), \ldots, x_7({\cal D})\}$ for $j=1, \ldots, 7$ .
\item[$\Xi$:] the set  $\{ \omega \in \Omega^{a,L}: \tau(a,\omega_i) < \infty, \mbox{ $i=1, \ldots, L$ } \}$.
\item[$n(\omega,j)$:] For $\omega \in \Xi$, $n(\omega,j)$ denotes the number of paths among $\omega_1, \ldots, \omega_{j \wedge V_7(\omega)}$ which hit the ``lucky" patch  $\Theta(\omega)$ before hitting the boundary of the growing cluster,  i.e.,   $n(\omega,j)=\max \{i:V_i(\omega) \leq j \} $.
\item[${\cal R}_i(\omega)$:] (p. 89) (defined for $\omega \in \Xi$ with $\omega_{V_i \tau(\Theta(\omega))} \in \gamma^-(\Theta(\omega)) \cup \gamma^+(\Theta(\omega))$).  
If $\omega_{V_i \tau(\Theta(\omega))}$  belongs to $\gamma^-(\Theta(\omega))$ (resp. $\gamma^+(\Theta(\omega))$)  then ${\cal R}_i(\omega)$ is defined as the collection of all paths from $\omega_{V_i \tau(\Theta(\omega))}$ to $\xi^l(\Theta(\omega))$  (resp. $\xi^r(\Theta(\omega))$) 
which avoid the cluster built by the paths $\varphi(\omega)_1, \ldots, \varphi(\omega)_{V_i(\omega)-1}$ and the interior of $\Theta(\omega)$ .
\item[$\mbox{r}(i,\omega)$:]  For all  $\omega \in \Xi$ satisfying that $\omega_{V_i \tau(\Theta(\omega))}$ belongs to $\gamma^-(\Theta(\omega))$ or to $\gamma^+(\Theta(\omega)),$ 
$\mbox{r}(i,\omega)$ is defined as a fixed path of minimum length in ${\cal R}_i(\omega)$. 
\item[${\hat \Gamma}$:] (p. 88) a lattice loop constructed from $\Gamma$; (${\hat \Gamma}$ may have additional lattice points which 
do not belong to $\Gamma$; an additional point needs to be included whenever two successive lattice points on $\Gamma$ are not neighbors in $\Ztwo$).  
\item[${\hat \Gamma}(i,\omega)$:] (defined for  $\omega \in \Xi$). Let ${\cal D}= \Theta(\omega)$. If $\omega_{V_i,\tau({\cal D})}$ belongs to 
$\Gamma$  
then ${\hat \Gamma}(i,\omega)$  denotes 
that part of the loop $\hat \Gamma$ or $R({\hat \Gamma})$ which starts at $\omega_{V_i,\tau({\cal D})}$, ends at
$\xi^m(\cal D)$ and which satisfies that all of its vertices which belong to $\Gamma$ also belong to ${\cal D}.$
If $\omega_{V_i \tau({\cal D})}$ belongs to $\gamma^-({\cal D})$ then ${\hat \Gamma}(i,\omega)$ denotes the part of the loop $\hat \Gamma$ which
starts at $\xi^l({\cal D})$ and ends at $\xi^m({\cal D})$ and, similarly, if $\omega_{V_i \tau({\cal D})}$ belongs to $\gamma^+({\cal D})$ then ${\hat \Gamma}(i,\omega)$ denotes the part of the loop $R(\hat \Gamma)$ which starts at $\xi^r({\cal D})$ and ends at $\xi^m({\cal D})$.
\item[$\alpha_i(\omega)$:] (defined for  $\omega \in \Xi$).  If $\omega_{V_i,\tau(\Theta(\omega))}$ belongs to $\Gamma$
then $\alpha_i(\omega)= ({\hat
\Gamma}(i,\omega), \beta_i(\Theta(\omega)))$. 
Otherwise, (i.e., $\omega_{V_i,\tau(\Theta(\omega))}$ belongs to $\gamma^-(\Theta(\omega))$ or to $\gamma^+(\Theta(\omega))$), we have $\alpha_i(\omega)= 
 (\mbox{r}(i,\omega),{\hat \Gamma}(i,\omega),
\beta_i(\Theta(\omega)))$. In either case, $\alpha_i(\omega)$ is a path  from $\omega_{V_i \tau({\Theta(\omega)})}$ to $x_i(\Theta(\omega))$ which hits $\partial
A_{k-1}(\varphi(\omega))$ at $x_i(\Theta(\omega)).$
\end{description}

 }

\bibliographystyle{plain}

\end{document}